# STABILITY OF INTEGRAL DELAY EQUATIONS AND STABILIZATION OF AGE-STRUCTURED MODELS

**Iasson Karafyllis[*] and Miroslav Krstic[**]**

[*]Dept. of Mathematics, National Technical University of Athens,
Zografou Campus, 15780, Athens, Greece, email: iasonkar@central.ntua.gr

[**]Dept. of Mechanical and Aerospace Eng., University of California, San Diego, La Jolla, CA 92093-0411, U.S.A., email: krstic@ucsd.edu

**Abstract**

We present bounded dynamic (but observer-free) output feedback laws that achieve global stabilization of equilibrium profiles of the partial differential equation (PDE) model of a simplified, age-structured chemostat model. The chemostat PDE state is positive-valued, which means that our global stabilization is established in the positive orthant of a particular function space—a rather non-standard situation, for which we develop non-standard tools. Our feedback laws do not employ any of the (distributed) parametric knowledge of the model. Moreover, we provide a family of highly unconventional Control Lyapunov Functionals (CLFs) for the age-structured chemostat PDE model. Two kinds of feedback stabilizers are provided: stabilizers with continuously adjusted input and sampled-data stabilizers. The results are based on the transformation of the first-order hyperbolic partial differential equation to an ordinary differential equation (one-dimensional) and an integral delay equation (infinite-dimensional). Novel stability results for integral delay equations are also provided; the results are of independent interest and allow the explicit construction of the CLF for the age-structured chemostat model.

## 1. Introduction

Continuous-time age-structured models are described by the so-called McKendrick-von Foerster equation (see [3,4,5,26] and the references therein), which is a first order hyperbolic Partial Differential Equation (PDE) with a non-local boundary condition. Age-structured models are natural extensions of standard chemostat models (see [27]). Optimal control problems for age-structured models have been studied (see [3,8,28] and the references therein). The ergodic theorem (see [11,12] and [26] for similar results on asymptotic similarity) has been proved an important tool for the study of the dynamics of continuous-time age structured models (see also [29] for a study of the existence of limit cycles).

This work initiates the study of the global stabilization problem by means of feedback control for age-structured models. More specifically, the design of explicit output feedback stabilizers is sought for the global stabilization of an equilibrium age profile for an age-structured chemostat model. Just as in other chemostat feedback control problems described by Ordinary Differential Equations (ODEs; see [9,13,14,15,21,22]), the dilution rate is selected to be the control input while the output is a weighted integral of the age distribution function. The assumed output functional form is chosen because it is an appropriate form for the expression of the measurement of the total concentration of the microorganism in the bioreactor or for the expression of any other measured variable (e.g., light absorption) that depends on the amount (and its size distribution) of the microorganism in the bioreactor. The main idea for the solution of the feedback control problem is the transformation of the first order hyperbolic PDE to an Integral Delay Equation (IDE; see [16]) and the application of the strong ergodic theorem. This feature differentiates the present work from recent works on feedback control problems for first order hyperbolic PDEs (see [1,2,6,7,16,20]).

The present work studies the global stabilization problem of an equilibrium age profile for an age-structured chemostat model by means of two kinds of feedback stabilizers: (i) a continuously applied feedback stabilizer, and (ii) a sampled-data feedback stabilizer. The entire model is assumed to be unknown and two cases are considered for the equilibrium value of the dilution rate: the case where the equilibrium value of the dilution rate is unknown (i.e., absolutely nothing is known about the model), and the case where the equilibrium value of the dilution rate is a priori known. In the first case, a family of observer-based (dynamic), output feedback laws with continuously adjusted dilution rate is proposed: the equilibrium value of the dilution rate is estimated by the observer. In the second case, a sampled-data output feedback law is proposed for arbitrarily sparse sampling schedule. In all cases, the dilution rate (control input) takes values in a pre-specified bounded interval and consequently input constraints are taken into account. The main idea for the solution of the feedback control problem is the transformation of the PDE to an ODE and an IDE. Some preliminary results for the sampled-data case, which are extended in the present work, were given in [17,18].

However, instead of simply designing dynamic, output feedback laws which guarantee global asymptotic stability of an equilibrium age profile, the present work has an additional goal: the explicit construction of a family of Control Lyapunov Functionals (CLFs) for the age-structured chemostat model. In order to achieve this goal, the present work provides/uses novel stability results on linear IDEs, which are of independent interest. The newly developed results, provide a Lyapunov-like proof of the scalar, strong ergodic theorem for special cases of the integral kernel. Stability results for linear IDEs similar to those studied in this work have been also studied in [23].

Since the state of the chemostat model is the population density of a particular age at a given time, the state of the chemostat PDE is non-negative valued. Accordingly, the desired equilibrium profile (a function of the age variable) is positive-valued. So the state space of this PDE system is the positive orthant in a particular function space. We pursue *global* stabilization of the positive equilibrium profile in such a state space. This requires a novel approach and even a novel formulation of stability estimates in which the norm of the state at the desired equilibrium is zero but takes the infinite value not only when the population density (of some age) is infinite but also when it is zero, i.e., we infinitely penalize the population death (the so-called "washout"), as we should. Our main idea in this development is a particular logarithmic transformation of the state, which penalizes both the overpopulated and underpopulated conditions, with an infinite penalty on the washout condition.

The structure of the paper is described next. In Section 2, we describe the chemostat stabilization problem in a precise way and we provide the statements of the main results of the paper (Theorem 2.1 and Theorem 2.4). Section 3 provides useful existing results for the uncontrolled PDE, while Section 4 is devoted to the presentation of stability results on IDEs, which allow us to construct CLFs for the chemostat problem. The proofs of the main results are provided in Section 5. Section 6 presents a result, which is similar to Theorem 2.1, but uses a reduced order observer instead of a full-order observer. Simulations, which illustrate the application of the obtained results, are given in Section 7. The concluding remarks of the paper are given in Section 8. Finally, the Appendix provides the proofs of certain auxiliary results.

**Notation.** Throughout this paper we adopt the following notation.

∗ For a real number $x \in \Re$, $[x]$ denotes the integer part of $x \in \Re$. $\Re_+$ denotes the interval $[0,+\infty)$.

∗ Let $U$ be an open subset of a metric space and $\Omega \subseteq \Re^m$ be a set. By $C^0(U;\Omega)$, we denote the class of continuous mappings on $U$, which take values in $\Omega$. When $U \subseteq \Re^n$, by $C^1(U;\Omega)$, we denote the class of continuously differentiable functions on $U$, which take values in $\Omega$. When $U = [a,b) \subseteq \Re$ (or $U = [a,b] \subseteq \Re$) with $a < b$, $C^0([a,b);\Omega)$ (or $C^0([a,b];\Omega)$) denotes all functions $f:[a,b) \to \Omega$ (or $f:[a,b] \to \Omega$), which are continuous on $(a,b)$ and satisfy $\lim_{s \to a^+}(f(s)) = f(a)$ (or $\lim_{s \to a^+}(f(s)) = f(a)$ and $\lim_{s \to b^-}(f(s)) = f(b)$). When $U = [a,b) \subseteq \Re$, $C^1([a,b);\Omega)$ denotes all functions $f:[a,b) \to \Omega$ which are continuously differentiable on $(a,b)$ and satisfy $\lim_{s \to a^+}(f(s)) = f(a)$ and $\lim_{h \to 0^+} h^{-1}(f(a+h) - f(a)) = \lim_{s \to a^+} f'(s)$.

∗ $K_\infty$ is the class of all strictly increasing, unbounded functions $a \in C^0(\Re_+;\Re_+)$, with $a(0) = 0$ (see [19]).

∗ For any subset $S \subseteq \Re$ and for any $A > 0$, $PC^1([0,A];S)$ denotes the class of all functions $f \in C^0([0,A];S)$ for which there exists a finite (or empty) set $B \subset (0,A)$ such that: (i) the derivative $f'(a)$ exists at every $a \in (0,A) \setminus B$ and is a continuous function on $(0,A) \setminus B$, (ii) all meaningful right and left limits of $f'(a)$ when $a$ tends to a point in $B \cup \{0,A\}$ exist and are finite.

∗ Let a function $f \in C^0(\Re_+ \times [0,A])$ be given, where $A > 0$ is a constant. We use the notation $f[t]$ to denote the profile at certain $t \geq 0$, i.e., $(f[t])(a) = f(t,a)$ for all $a \in [0,A]$.

∗ Let a function $x \in C^0([-A,+\infty);\Re)$ be given, where $A > 0$ is a constant. We use the notation $x_t \in C^0([-A,0];\Re)$ to denote the "$A$–history" of $x$ at certain $t \geq 0$, i.e., $(x_t)(-a) = x(t-a)$ for all $a \in [0,A]$.

∗ Let $0 < D_{\min} < D_{\max}$ be given constants. The saturation function $sat(x)$ for the interval $[D_{\min}, D_{\max}]$ is defined by $sat(x) := \min(D_{\max}, \max(D_{\min}, x))$, for all $x \in \Re$.

## 2. Problem Description and Main Results

### 2. I. The model

Consider the age-structured chemostat model:

$$\frac{\partial f}{\partial t}(t,a) + \frac{\partial f}{\partial a}(t,a) = -(\mu(a) + D(t))f(t,a), \text{ for } t > 0,\ a \in (0,A) \tag{2.1}$$

$$f(t,0) = \int_0^A k(a)f(t,a)da, \text{ for } t \geq 0 \tag{2.2}$$

where $D(t) \in [D_{\min}, D_{\max}]$ is the dilution rate, $D_{\max} > D_{\min} > 0$ are constants, $A > 0$ is a constant and $\mu:[0,A] \to \Re_+$, $k:[0,A] \to \Re_+$ are continuous functions with $\int_0^A k(a)da > 0$. System (2.1), (2.2) is a continuous age-structured model of a microbial population in a chemostat. The function $\mu(a) \geq 0$ is called the mortality function, the function $f(t,a)$ denotes the density of the population of age $a \in [0, A]$ at time $t \geq 0$ and the function $k(a) \geq 0$ is the birth modulus of the population. The boundary condition (2.2) is the renewal condition, which determines the number of newborn individuals $f(t,0)$. Finally, $A > 0$ is the maximum reproductive age. Physically meaningful solutions of (2.1), (2.2) are only the non-negative solutions, i.e., solutions satisfying $f(t,a) \geq 0$, for all $(t,a) \in \Re_+ \times [0, A]$.

The chemostat model (2.1), (2.2) is derived by neglecting the dependence of the growth of the microorganism on the concentration of a limiting substrate. A more accurate model would involve an enlarged system that has one PDE for the age distribution, coupled with one ODE for the substrate (as proposed in [29], in the context of studying limit cycles with constant dilution rates). However, the approach of neglecting the nutrient's equation in the chemostat is not new (see for example [25]).

We assume that there exists $D^* \in (D_{\min}, D_{\max})$ such that

$$1 = \int_0^A k(a) \exp\left(-D^* a - \int_0^a \mu(s) ds\right) da \tag{2.3}$$

This assumption is necessary for the existence of an equilibrium point for the control system (2.1), (2.2), which is different from the identically zero function. Any function of the form:

$$f^*(a) = M \exp\left(-D^* a - \int_0^a \mu(s) ds\right), \text{ for } a \in [0, A] \tag{2.4}$$

with $M \geq 0$ being an arbitrary constant, is an equilibrium point for the control system (2.1), (2.2) with $D(t) \equiv D^*$. Notice that there is a continuum of equilibria.

The measured output of the control system (2.1), (2.2) is given by the equation:

$$y(t) = \int_0^A p(a) f(t,a) da, \text{ for } t \geq 0 \tag{2.5}$$

where $p : [0, A] \to \Re_+$ is a continuous function with $\int_0^A p(a) da > 0$. Notice that the case $p(a) \equiv 1$ corresponds to the total concentration of the microorganism in the chemostat.

## 2.II. Feedback Control with Continuously Adjusted Input

Let $y^* > 0$ be an arbitrary constant (the set point) and let $f^*(a)$ be the equilibrium age profile given by (2.4) with $M = y^* \left( \int_0^A p(a) \exp\left(-D^* a - \int_0^a \mu(s) ds\right) da \right)^{-1}$. Consider the dynamic feedback law given by

$$\dot{z}_1(t) = z_2(t) - D(t) - l_1\left(z_1(t) - \ln\left(\frac{y(t)}{y^*}\right)\right)$$

$$\dot{z}_2(t) = -l_2\left(z_1(t) - \ln\left(\frac{y(t)}{y^*}\right)\right) \tag{2.6}$$

$$z(t) = (z_1(t), z_2(t))' \in \Re^2$$

and

$$D(t) = sat\left(z_2(t) + \gamma \ln\left(\frac{y(t)}{y^*}\right)\right) \tag{2.7}$$

where $l_1, l_2, \gamma > 0$ are constants. Next consider solutions of the initial-value problem (2.1), (2.2), (2.5), (2.6), (2.7) with initial condition $(f_0, z_0) \in \tilde{X} \times \Re^2$, where $\tilde{X}$ is the set $\tilde{X} = \left\{ f \in PC^1([0,A];(0,+\infty)) : f(0) = \int_0^A k(a) f(a) da \right\}$. By a solution of the initial-value problem (2.1), (2.2), (2.5), (2.6), (2.7) with initial condition $(f_0, z_0) \in \tilde{X} \times \Re^2$, we mean a pair of mappings $f \in C^0([0,\tau) \times [0,A];(0,+\infty))$, $z \in C^1([0,\tau); \Re^2)$, where $\tau \in (0,+\infty]$, which satisfies the following properties:

(i) $f \in C^1(D_f; (0,+\infty))$, where $D_f = \{(t,a) \in (0,\tau) \times (0,A) : (a-t) \notin B \cup \{0,A\}\}$ and $B \subseteq (0,A)$ is the finite (possibly empty) set where the derivative of $f_0 \in \tilde{X}$ is not defined or is not continuous,

(ii) $f[t] \in \tilde{X}$ for all $t \in [0,\tau)$, where $(f[t])(a) = f(t,a)$ for $a \in [0,A]$ (see Notation),

(iii) equations (2.5), (2.6), (2.7) hold for all $t \in [0,\tau)$,

(iv) equation $\frac{\partial f}{\partial t}(t,a) + \frac{\partial f}{\partial a}(t,a) = -(\mu(a) + D(t)) f(t,a)$ holds for all $(t,a) \in D_f$, and

(v) $z(0) = z_0 = (z_{1,0}, z_{2,0})$, $f(0,a) = f_0(a)$ for all $a \in [0,A]$.

The mapping $[0,\tau) \ni t \to (f[t], z(t)) \in \tilde{X} \times \Re^2$ is called the *solution of the closed-loop system (2.1), (2.2), (2.5) with (2.6), (2.7) and initial condition* $(f_0, z_0) \in \tilde{X} \times \Re^2$ defined for $t \in [0,\tau)$.

Define the functional $\Pi : C^0([0,A]; \Re) \to \Re$ by means of the equation

$$\Pi(f) := \frac{\int_0^A f(a) \left( \int_a^A k(s) \exp\left( \int_s^a (\mu(l) + D^*) dl \right) ds \right) da}{\int_0^A a k(a) f^*(a) da} \tag{2.8}$$

and assume that the following technical assumption holds for the non-negative function

$$\tilde{k}(a) := k(a) \exp\left(-D^* a - \int_0^a \mu(s) ds\right), \text{ for } a \in [0,A] \tag{2.9}$$

that satisfies $\int_0^A \tilde{k}(a) da = 1$ (recall (2.3)):

**(A)** *There exists a constant $\lambda > 0$ such that $\int_0^A \left| \tilde{k}(a) - r\lambda \int_a^A \tilde{k}(s) ds \right| da < 1$, where $r := \left( \int_0^A a \tilde{k}(a) da \right)^{-1}$.*

We are now ready to state the first main result of the present work, which provides stabilizers with continuously adjusted input.

**Theorem 2.1 (continuously adjusted input and unknown equilibrium value of the dilution rate):** *Consider the age-structured chemostat model (2.1), (2.2) with $k \in PC^1([0,A]; \Re_+)$ under Assumption (A). Then for every $f_0 \in \tilde{X}$ and $z_0 \in \Re^2$ there exists a unique solution of the closed-loop (2.1), (2.2), (2.5) with (2.6), (2.7) and initial condition $(f_0, z_0) \in \tilde{X} \times \Re^2$. Furthermore, there exist a constant $L > 0$ and a function $\rho \in K_\infty$ such that for every $f_0 \in \tilde{X}$ and $z_0 \in \Re^2$ the unique solution of* the closed-loop (2.1), (2.2), (2.5) with (2.6), (2.7) and initial condition $(f_0, z_0) \in \tilde{X} \times \Re^2$ is defined for all $t \geq 0$ and satisfies the following estimate:

$$\max_{a \in [0,A]}\left(\left|\ln\left(\frac{f(t,a)}{f^*(a)}\right)\right|\right) + |z_1(t)| + |z_2(t) - D^*| \leq$$
$$\exp\left(-\frac{L}{4}t\right)\rho\left(\max_{a \in [0,A]}\left(\left|\ln\left(\frac{f_0(a)}{f^*(a)}\right)\right|\right) + |z_{1,0}| + |z_{2,0} - D^*|\right), \text{ for all } t \geq 0 \quad (2.10)$$

Moreover, let $p_1, p_2 > 0$ be a pair of constants satisfying $(2 + l_1 p_1 - 2l_2 p_2)^2 < 8l_1 p_1 - 4l_2 p_1^2$, $p_1^2 < 4p_2$. Then the continuous functional $W : \Re^2 \times C^0([0, A]; (0, +\infty)) \to \Re_+$ defined by:

$$W(z, f) := (\ln(\Pi(f)))^2 + G\sqrt{Q(z, f)} + \beta Q(z, f) \quad (2.11)$$

where $\beta \geq 0$ is an arbitrary constant,
$Q(z, f) :=$

$$(z_1 - \ln(\Pi(f)))^2 - p_1(z_1 - \ln(\Pi(f)))(z_2 - D^*) + p_2(z_2 - D^*)^2 + \frac{M}{2}\left(\frac{\max_{a \in [0,A]}\left(\exp(-\sigma a)\left|\frac{f(a) - \Pi(f)f^*(a)}{f^*(a)}\right|\right)}{\min\left(\Pi(f), \min_{a \in [0,A]}\left(\frac{f(a)}{f^*(a)}\right)\right)}\right)^2 \quad (2.12)$$

$\sigma > 0$ is a sufficiently small constant and $M, G > 0$ are sufficiently large constants, is a Lyapunov functional for the closed-loop system (2.1), (2.2), (2.5) with (2.6), (2.7), in the sense that every solution $(f[t], z(t)) \in \tilde{X} \times \Re^2$ of the closed-loop system (2.1), (2.2), (2.5) with (2.6), (2.7) satisfies the inequality:

$$\limsup_{h \to 0^+}\left(h^{-1}(W(z(t+h), f[t+h]) - W(z(t), f[t]))\right) \leq -L\frac{W(z(t), f[t])}{1 + \sqrt{W(z(t), f[t])}}, \text{ for all } t \geq 0 \quad (2.13)$$

As remarked in the Introduction, Theorem 2.1 does not only provide formulas for dynamic output feedback stabilizers that guarantee global asymptotic stability of the selected equilibrium age profile, but also provides explicit formulas for a family of CLFs for system (2.1), (2.2). Indeed, the continuous functional $W : \Re^2 \times C^0([0, A]; (0, +\infty)) \to \Re_+$ defined by (2.11), (2.12) is a CLF for system (2.1), (2.2).

**Remark 2.2:**
**i)** The family of feedback laws (2.6), (2.7) (parameterized by $l_1, l_2, \gamma > 0$) guarantees global asymptotic stabilization of every selected equilibrium age profile. Moreover, the feedback law (2.6), (2.7) achieves a global exponential convergence rate (see estimate (2.10)), in the sense that estimate (2.10) holds for all physically meaningful initial conditions ($f_0 \in \tilde{X}$). As indicated in the Introduction, the logarithmic penalty in (2.10) penalizes both the overpopulated and underpopulated conditions, with an infinite penalty on zero density for some age. The state converges to the desired equilibrium profiles from all positive initial conditions, but not from the zero-density initial condition, which itself is an equilibrium (population cannot develop from a "dead" initial state).
**ii)** The feedback law (2.6), (2.7) is a dynamic output feedback law. The subsystem (2.6) is an observer that primarily estimates the equilibrium value of the dilution rate $D^*$. The observer (2.6) is a highly reduced order, since it estimates only two variables, the afore-mentioned constant $D^*$ and the scalar functional of the infinite-dimensional state, $\Pi(f)$, introduced in (2.8). All the remaining infinitely many states are not estimated. This is the key achievement of our paper—attaining stabilization without the estimation of nearly the entire infinite-dimensional state and proving this result in an appropriately constructed transformed representation of that unmeasured infinite-dimensional state.

**iii)** The family of feedback laws (2.6), (2.7) does not require knowledge of the mortality function of the population, the birth modulus of the population and the maximum reproductive age of the population. Accordingly, it does not require the knowledge of the equilibrium value of the dilution rate $D^*$ either. Instead, $D^*$ is estimated by the observer state $z_2(t)$ (see estimate (2.10)).

**iv)** The feedback law (2.6), (2.7) can work with arbitrary input constraints. The only condition that needs to be satisfied is that the equilibrium value of the dilution rate $D^*$ must satisfy the input constraints, i.e., $D^* \in (D_{\min}, D_{\max})$, which is a reasonable requirement (otherwise the selected equilibrium age profile is not feasible).

**v)** The parameters $l_1, l_2, \gamma > 0$ can be used by the control practitioner for tuning the controller (2.6), (2.7): the selection of the values of these parameters affects the value of the constant $L > 0$ that determines the exponential convergence rate. Since the proof of Theorem 2.1 is constructive, useful formulas showing the dependence of the constant $L > 0$ on the parameters $l_1, l_2, \gamma > 0$ are established in the proof of Theorem 2.1.

**vi)** It should be noted that for every pair of constants $l_1, l_2 > 0$ it is possible to find constants $p_1, p_2 > 0$ satisfying $(2 + l_1 p_1 - 2l_2 p_2)^2 < 8l_1 p_1 - 4l_2 p_1^2$, $p_1^2 < 4p_2$. Indeed, for every $l_1, l_2 > 0$ the matrix $\begin{bmatrix} -l_1 & 1 \\ -l_2 & 0 \end{bmatrix}$ is a Hurwitz matrix. Consequently, there exists a positive definite matrix $\begin{bmatrix} 1 & -p_1/2 \\ -p_1/2 & p_2 \end{bmatrix}$ so that the matrix

$$\begin{bmatrix} -l_1 & -l_2 \\ 1 & 0 \end{bmatrix} \begin{bmatrix} 1 & -p_1/2 \\ -p_1/2 & p_2 \end{bmatrix} + \begin{bmatrix} 1 & -p_1/2 \\ -p_1/2 & p_2 \end{bmatrix} \begin{bmatrix} -l_1 & 1 \\ -l_2 & 0 \end{bmatrix} = \begin{bmatrix} -2l_1 + l_2 p_1 & 1 + l_1 p_1/2 - l_2 p_2 \\ 1 + l_1 p_1/2 - l_2 p_2 & -p_1 \end{bmatrix}$$

is negative definite. This implies the inequalities $p_1^2 < 4p_2$ and $(2 + l_1 p_1 - 2l_2 p_2)^2 < 8l_1 p_1 - 4l_2 p_1^2$.

**vii)** The main idea for the construction of the feedback law (2.6), (2.7) is the transformation of the PDE problem (2.1), (2.2) into a system that consists of an ODE and an IDE along with the logarithmic output transformation $Y(t) = \ln\left(\dfrac{y(t)}{y^*}\right)$. The transformations are presented in Figure 1 and are exploited rigorously in the proof of Theorem 2.1. Figure 1 also shows that the full-order observer (2.6) is actually an observer for the system $\dot{\eta}(t) = D^*(t) - D(t)$, $\dot{D}^*(t) = 0$.

## 2.III. Checking Assumption (A)

Theorem 2.1 assumes that the birth modulus of the population satisfies Assumption (A). This is not an assumption that is needed for the establishment of the exponential estimate (2.10). Estimate (2.10) could have been established without Assumption (A) by means of the strong ergodic theorem (see Section 3). The role of Assumption (A) is crucial for the establishment of the CLF, given by (2.11), (2.12). However, since Assumption (A) demands a specific property for the function $\tilde{k}(a) := k(a) \exp\left(-D^* a - \int_0^a \mu(s) ds\right)$ that involves the (unknown) equilibrium value of the dilution rate $D^*$, the verification of the validity of Assumption (A) becomes an issue. The following proposition provides useful sufficient conditions for Assumption (A). Its proof is provided in the Appendix.

**Proposition 2.3 (Means of checking Assumption (A)):** *Let $\tilde{k} \in C^0([0, A]; \Re)$ be a function that satisfies the following assumption:*

**(B)** *The function $\tilde{k} \in C^0([0,A]; \Re)$ satisfies $\tilde{k}(a) \geq 0$ for all $a \in [0,A]$ and $\int_0^A \tilde{k}(a)da = 1$. Moreover, there exists a constant $\varepsilon > 0$ such that the set $S_\varepsilon = \{a \in [0,T] : \tilde{k}(a) \leq \varepsilon\}$, where $T := \sup\{a \in [0,A] : \tilde{k}(a) > 0\}$, has Lebesgue measure $|S_\varepsilon| < (2r)^{-1}$, where $r := \left(\int_0^A a\tilde{k}(a)da\right)^{-1}$.*

*Then for every $\lambda \in [0, r^{-1}\varepsilon]$ it holds that $\int_0^A \left|\tilde{k}(a) - r\lambda \int_a^A \tilde{k}(s)ds\right|da \leq 1 - \lambda(1 - 2r|S_\varepsilon|)$.*

Proposition 2.3 shows that Assumption (A) is valid for a function that satisfies Assumption (B). On the other hand, we know that Assumption (B) holds for every function $\tilde{k} \in C^0([0,A]; \Re_+)$ satisfying $\int_0^A \tilde{k}(a)da = 1$ and having only a finite number of zeros in the interval $[0,A]$. Since $\tilde{k}(a) := k(a)\exp\left(-D^*a - \int_0^a \mu(s)ds\right)$, we can be sure that Assumption (A) necessarily holds *for all birth moduli $k \in C^0([0,A]; \Re_+)$ of the population with only a finite number of zeros in the interval $[0,A]$, no matter what the equilibrium value of the dilution rate $D^*$ is and no matter what the mortality function $\mu : [0,A] \to \Re_+$ is.*

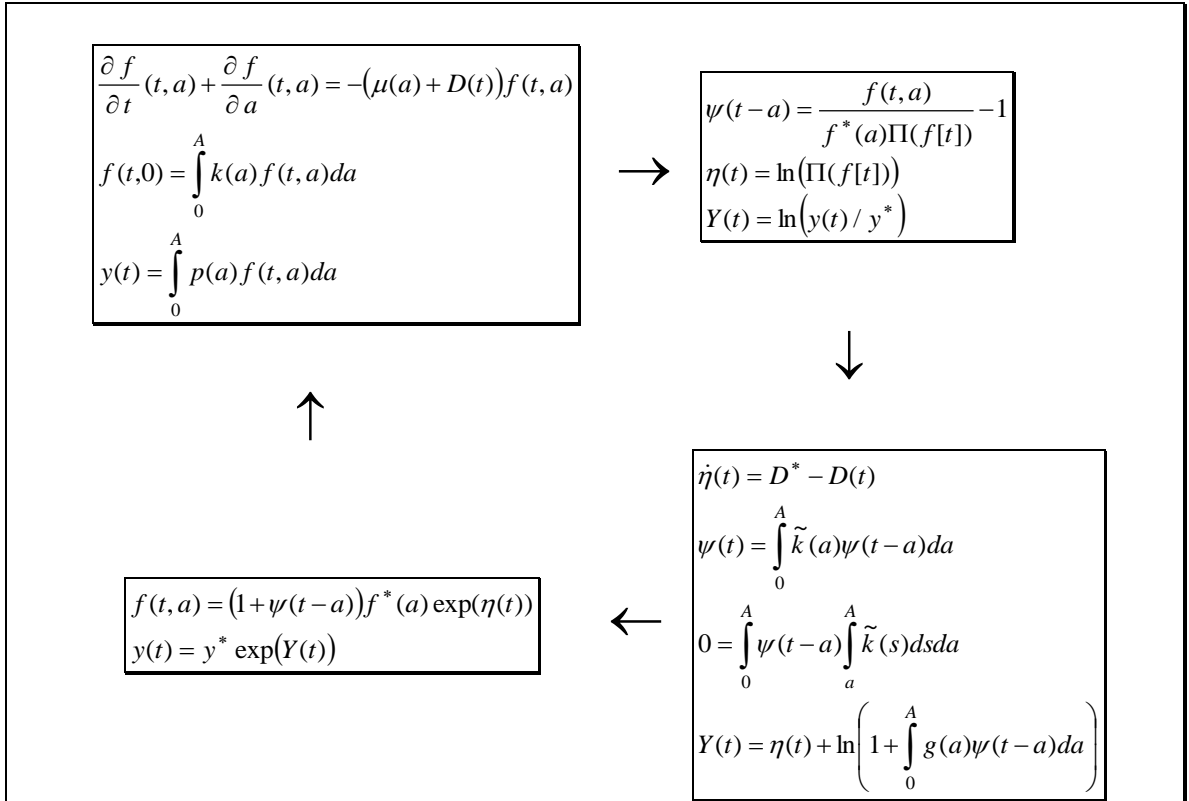

**Figure 1:** The transformation of the PDE (2.1) with boundary condition given by (2.2) to an IDE and an ODE and the inverse transformation.

## 2.IV. Sampled-Data Control

On the other hand, when the equilibrium value of the dilution rate $D^*$ is a priori known, then we are in a position to achieve sampled-data stabilization. Let $T > 0$ be the sampling period and consider the closed-loop system (2.1), (2.2), (2.5) with the sample-and-hold feedback law

$$D(t) = sat\left(D^* + T^{-1} \ln\left(\frac{y(iT)}{y^*}\right)\right), \text{ for all } t \in [iT,(i+1)T) \text{ and for all integers } i \geq 0 \quad (2.14)$$

By a solution of the initial-value problem (2.1), (2.2), (2.5), (2.14) with initial condition $f_0 \in \tilde{X}$, where $\tilde{X}$ is the set $\tilde{X} = \left\{ f \in PC^1([0,A];(0,+\infty)): f(0) = \int_0^A k(a)f(a)da \right\}$, we mean a mapping $f \in C^0([0,\tau) \times [0,A];(0,+\infty))$, where $\tau \in (0,+\infty]$, which satisfies the following properties:

(i) $f \in C^1(D_f;(0,+\infty))$, where $D_f = \left\{ (t,a) \in (0,\tau) \times (0,A): (a-t) \notin B \cup \{0,A\}, t \neq T\lfloor T^{-1}t \rfloor \right\}$ and $B \subseteq (0,A)$ is the finite (possibly empty) set where the derivative of $f_0 \in \tilde{X}$ is not defined or is not continuous,

(ii) $f[t] \in \tilde{X}$ for all $t \in [0,\tau)$,

(iii) equations (2.5), (2.14) hold for all $t \in [0,\tau)$,

(iv) equation $\frac{\partial f}{\partial t}(t,a) + \frac{\partial f}{\partial a}(t,a) = -(\mu(a) + D(t))f(t,a)$ holds for all $(t,a) \in D_f$, and

(v) $f(0,a) = f_0(a)$ for all $a \in [0,A]$.

The mapping $[0,\tau) \ni t \to f[t] \in \tilde{X}$ is called the *solution of the closed-loop system (2.1), (2.2), (2.5) with (2.14)* and initial condition $f_0 \in \tilde{X}$ defined for $t \in [0,\tau)$.

We are now ready to state the second main result of the present work.

**Theorem 2.4 (sampled-data feedback and known equilibrium value for the dilution rate):** *Consider the age-structured chemostat model (2.1), (2.2) with $k \in PC^1([0,A];\Re_+)$. Then for every $f_0 \in \tilde{X}$ there exists a unique solution of the closed-loop (2.1), (2.2), (2.5) with (2.14) and initial condition $f_0 \in \tilde{X}$. Furthermore, there exist a constant $L > 0$ and a function $\rho \in K_\infty$ such that for every $f_0 \in \tilde{X}$ the unique solution of the closed-loop (2.1), (2.2), (2.5) with (2.14) and initial condition $f_0 \in \tilde{X}$ is defined for all $t \geq 0$ and satisfies the following estimate:*

$$\max_{a \in [0,A]}\left(\left|\ln\left(\frac{f(t,a)}{f^*(a)}\right)\right|\right) \leq \exp(-Lt)\rho\left(\max_{a \in [0,A]}\left(\left|\ln\left(\frac{f_0(a)}{f^*(a)}\right)\right|\right)\right), \text{ for all } t \geq 0 \quad (2.15)$$

The differences of Theorem 2.4 with Theorem 2.1 are:
(i) Theorem 2.4 applies the sampled-data feedback (2.14) while Theorem 2.1 applies a continuously adjusted feedback,
(ii) Theorem 2.4 assumes knowledge of the equilibrium value of the dilution rate $D^*$,
(iii) Theorem 2.4 does not assume property (A) but does not provide a CLF for the system. This was explained above: assumption (A) is only needed for the explicit construction of a CLF.

Finally, the reader should notice that there is no constraint for the sampling period $T > 0$: arbitrarily large values for $T > 0$ are allowed (arbitrarily sparse sampling). In the case where the output is given by

$$y(t) = f(t,0), \text{ for } t \geq 0 \quad (2.16)$$

instead of (2.5), the proof of Theorem 2.4 works with only minor changes (the proof is omitted; this is the case considered in [18]).

*2.V. Ideas Behind the Proofs of the Main Results*

The basic tool for the proofs of the main results of the present work is the transformation shown in Figure 1. The main idea comes from the recent work [16]: the transformation of a first-order hyperbolic PDE to an IDE. However, if we applied the results of [16] in a straightforward way, then we would end up with the following IDE:

$$v(t) = \int_0^A k(a) \exp\left(-\int_0^a \mu(s)ds\right) \exp\left(-\int_{t-a}^t D(s)ds\right) v(t-a)da, \quad (2.17)$$

where $v(t) = f(t,0)$ and $f(t,a) = \exp\left(-\int_0^a \mu(s)ds\right) \exp\left(-\int_{t-a}^t D(s)ds\right) v(t-a)$. However, the IDE is input-dependent. Instead, we would like to describe the effect of the control input in a more convenient way: this is achieved by introducing one more state

$$\eta(t) = \ln(\Pi(f[t])), \quad (2.18)$$

where $\Pi$ is given by (2.8). The evolution of $\eta(t)$ is described by the ODE $\dot\eta(t) = D^* - D(t)$. Then we are in a position to obtain the transformation

$$\psi(t-a) = \frac{f(t,a)}{f^*(a)\Pi(f[t])} - 1, \text{ for all } (t,a) \in \Re_+ \times [0,A] \quad (2.19)$$

which decomposes the dynamics of (2.17) to the input-independent dynamics of the IDE $\psi(t) = \int_0^A \tilde k(a)\psi(t-a)da$ evolving on the subspace described by the equation $\int_0^A \psi(t-a)\int_a^A \tilde k(s)dsda = 0$ and the input-dependent ODE $\dot\eta(t) = D^* - D(t)$. After achieving this objective, the next step is the stability analysis of the zero solution of the IDE $\psi(t) = \int_0^A \tilde k(a)\psi(t-a)da$: this is exactly the point where the strong ergodic theorem or the results on linear IDEs are used.

# 3. The Uncontrolled PDE

The present section aims to give to the reader the background mathematical knowledge which is used for the study of age-structured PDEs. More specifically, we aim to make the reader familiar to the strong ergodic theorem for age-structured PDEs and to show the relation of age-structured PDEs to linear IDEs.

Let $A > 0$ be a constant and let $\mu:[0,A] \to \Re_+$, $k:[0,A] \to \Re_+$ be continuous functions with $\int_0^A k(a)da > 0$. Consider the initial value PDE problem:

$$\frac{\partial z}{\partial t}(t,a) + \frac{\partial z}{\partial a}(t,a) = -\mu(a)z(t,a), \text{ for } t > 0, \ a \in (0,A) \quad (3.1)$$

$$z(t,0) = \int_0^A k(a)z(t,a)da, \text{ for } t \geq 0 \quad (3.2)$$

with initial condition $z(0,a) = z_0(a)$ for all $a \in [0,A]$. The following existence and uniqueness result follows directly from Proposition 2.4 in [10] and Theorems 1.3-1.4 on pages 102-104 in [24]:

**Lemma 3.1 (existence/uniqueness):** *For each absolutely continuous function $z_0 \in C^0([0, A]; \Re)$ with $z_0(0) = \int_0^A k(a) z_0(a) da$, there exists a unique function $z:[0,+\infty) \times [0, A] \to \Re$ with $z(0, a) = z_0(a)$ for all $a \in [0, A]$ that satisfies: (a) For each $t \geq 0$, the function $z[t]$ defined by $(z[t])(a) = z(t, a)$ for $a \in [0, A]$ is absolutely continuous and satisfies $(z[t])(0) = \int_0^A k(a)(z[t])(a) da$ for all $t \geq 0$, (b) the mapping $\Re_+ \ni t \to z[t] \in L^1([0, A]; \Re)$ is continuously differentiable, and (c) equation (3.1) holds for almost all $t > 0$ and $a \in (0, A)$. Moreover, if $z_0(a) \geq 0$ for all $a \in [0, A]$ then $z(t, a) \geq 0$, for all $(t, a) \in \Re_+ \times [0, A]$.*

The function $z:[0,+\infty) \times [0, A] \to \Re$ is called the solution of (3.1), (3.2). When additional regularity properties hold then the solution of (3.1), (3.2) satisfies the properties shown by the following lemma.

**Lemma 3.2 (regularity/relation to IDEs):** *If $k \in PC^1([0, A]; \Re_+)$, then for every $z_0 \in PC^1([0, A]; \Re)$ satisfying $z_0(0) = \int_0^A k(a) z_0(a) da$ the function $z:[0,+\infty) \times [0, A] \to \Re$ from Lemma 3.1 is $C^1$ on*

$$S = \{ (t, a) \in (0,+\infty) \times (0, A) : (a - t) \notin B \cup \{0, A\} \}$$

*where $B$ is the finite (or empty) set where the derivative of $z_0$ is not defined or is not continuous, satisfies (3.1) on $S$ and equation (3.2) for all $t \geq 0$. Also,*

$$z(t, a) = \exp\left(-\int_0^a \mu(s) ds\right) v(t - a), \text{ for all } (t, a) \in \Re_+ \times [0, A] \quad (3.3)$$

*where $v \in C^0([-A, +\infty); \Re) \cap C^1((0, +\infty); \Re)$ is the unique solution of the Integral Delay Equation (IDE):*

$$v(t) = \int_0^A k(a) \exp\left(-\int_0^a \mu(s) ds\right) v(t - a) da, \text{ for } t \geq 0 \quad (3.4)$$

*with initial condition $v(-a) = \exp\left(\int_0^a \mu(s) ds\right) z_0(a)$, for all $a \in (0, A]$.*

Lemma 3.2 is obtained by integration on the characteristic lines of (3.1). The solution $v \in C^0([-A, +\infty); \Re) \cap C^1((0, +\infty); \Re)$ of the IDE (3.4) is obtained as the solution of the delay differential equation

$$\dot{v}(t) = \bar{k}(0) v(t) - \bar{k}(A) v(t - A) + \int_0^A \frac{d\bar{k}}{da}(a) v(t - a) da \quad (3.5)$$

where $\bar{k}(a) := k(a) \exp\left(-\int_0^a \mu(s) ds\right)$ for $a \in [0, A]$. The differential equation (3.5) is obtained by formal differentiation of the IDE (3.4) and its solution satisfies (3.4) (the verification requires integration by parts).

It is straightforward to show that the function $h(D) = \int_0^A k(a)\exp\left(-Da - \int_0^a \mu(s)ds\right)da$ is strictly decreasing with $\lim_{D\to+\infty} h(D) = 0$ and $\lim_{D\to-\infty} h(D) = +\infty$. Therefore, there exists a unique $D^* \in \Re$ such that (2.3) holds. Equation (2.3) is the Lotka-Sharpe condition [4]. The following strong ergodicity result follows from the results of Section 3 in [11] and Proposition 3.2 in [10]:

**Theorem 3.3 (scalar strong ergodic theorem):** *Let $D^* \in \Re$ be the unique solution of (2.3). Then, there exist constants $\varepsilon > 0$, $K \geq 1$ such that for every absolutely continuous function $z_0 \in C^0([0,A];\Re)$ with $z_0(0) = \int_0^A k(a)z_0(a)da$, the corresponding solution $z:[0,+\infty)\times[0,A] \to \Re$ of (3.1), (3.2) satisfies for all $t \geq 0$:*

$$\int_0^A \exp\left(\int_0^a \mu(s)ds\right)\left|z(t,a) - \exp\left(D^*(t-a) - \int_0^a \mu(s)ds\right)\Phi(z_0)\right|da \leq K\exp((D^*-\varepsilon)t)\int_0^A \exp\left(\int_0^a \mu(s)ds\right)|z_0(a)|da \quad (3.6)$$

*where $\Phi: L^1([0,A];\Re) \to \Re$ is the linear continuous functional defined by:*

$$\Phi(z_0) := \frac{\int_0^A z_0(a)\int_a^A k(s)\exp\left(\int_s^a (\mu(l) + D^*)dl\right)ds\,da}{\int_0^A ak(a)\exp\left(-\int_0^a (\mu(l) + D^*)dl\right)da} \quad (3.7)$$

## 4. Results on Linear Integral Delay Equations

Since the previous sections have demonstrated the relation of age-structured PDEs to linear IDEs, we next focus on the study of linear IDEs. The present section provides stability results for the system described by the following linear IDE:

$$x(t) = \int_0^A \varphi(a)x(t-a)da \quad (4.1)$$

where $x(t) \in \Re$, $A > 0$ is a constant and $\varphi \in C^0([0,A];\Re)$. The results of the present section allow the construction of Lyapunov functionals for linear IDEs, which provide formulas for Lyapunov functionals of age-structured PDEs (since the zero dynamics of the controlled age-structured model are described by linear IDEs). All proofs of the results of the present section are provided in the Appendix.

*4.I. The notion of the solution-existence/uniqueness*

For every $x_0 \in C^0([-A,0];\Re)$ with $x_0(0) = \int_0^A \varphi(a)x_0(-a)da$ there exists a unique function $x \in C^0([-A,+\infty);\Re)$ that satisfies (4.1) for $t \geq 0$ and $x(-a) = x_0(-a)$ for all $a \in [0,A]$. This function is called the solution of (4.1) with initial condition $x_0 \in C^0([-A,0];\Re)$. The solution is obtained as the solution of the neutral delay equation $\frac{d}{dt}\left(x(t) - \int_0^A \varphi(a)x(t-a)da\right) = 0$ (Theorem 1.1 on page 256 in [10] guarantees the existence of a unique function $x \in C^0([-A,+\infty);\Re) \cap C^1((0,+\infty);\Re)$ that satisfies $\frac{d}{dt}\left(x(t) - \int_0^A \varphi(a)x(t-a)da\right) = 0$ for $t \geq 0$ and $x(-a) = x_0(-a)$ for all $a \in [0, A]$).

Therefore, the IDE (4.1) defines a dynamical system on $X = \left\{ x \in C^0([-A,0];\Re) : x(0) = \int_0^A \varphi(a)x(-a)da \right\}$ with state $x_t \in X$, where $(x_t)(-a) = x(t-a)$ for all $a \in [0, A]$ (see Notation).

*4.II. A Basic Estimate and its Consequences*

The first result of this section provides useful bounds for the solution of (4.1) with non-negative kernel. Notice that the following lemma allows discontinuous solutions for (4.1) as well as discontinuous initial conditions.

**Lemma 4.1 (A Basic Estimate for the solution of linear IDEs):** *Let $\varphi \in C^0([0,A];\Re_+)$ be a given function with $\int_0^A \varphi(a)da \geq 1$ and consider the IDE (4.1). Let $\delta > 0$ be an arbitrary constant with $\int_0^\delta \varphi(a)da < 1$. Then for every $\xi \in L^\infty([-A,0);\Re)$ there exists a unique function $x \in L^\infty_{loc}([-A,+\infty);\Re)$ with $x(a) = \xi(a)$ for $a \in [-A,0)$ that satisfies (4.1) for $t \geq 0$ a.e.. Moreover, $x \in L^\infty_{loc}([-A,+\infty);\Re)$ satisfies for all $t \geq 0$ the following inequality:*

$$\min\left(\inf_{-A \leq a < 0}(\xi(a)), \left(\frac{L-c}{1-c}\right)^{1+h^{-1}t} \inf_{-A \leq a < 0}(\xi(a))\right) \leq \inf_{-A \leq a < 0}(x(t+a))$$
$$\leq \sup_{-A \leq a < 0}(x(t+a)) \leq \max\left(\left(\frac{L-c}{1-c}\right)^{1+h^{-1}t} \sup_{-A \leq a < 0}(\xi(a)), \sup_{-A \leq a < 0}(\xi(a))\right) \quad (4.2)$$

*where $h := \min(\delta, A - \delta)$, $L := \int_0^A \varphi(a)da \geq 1$, $c := \int_0^\delta \varphi(a)da < 1$.*

A direct consequence of Lemma 4.1 and Lemma 3.2 is that if $k \in PC^1([0,A];\Re_+)$, then for every $z_0 \in PC^1([0,A];\Re)$ satisfying $z_0(0) = \int_0^A k(a)z_0(a)da$ and $z_0(a) > 0$ for all $a \in [0,A]$, the corresponding solution of (3.1), (3.2) satisfies $z(t,a) > 0$, for all $(t,a) \in \Re_+ \times [0,A]$. To see this, notice that if $\int_0^A k(a)\exp\left(-\int_0^a \mu(s)ds\right)da \geq 1$ then we may apply Lemma 3.2 and Lemma 4.1 directly for the IDE (3.4). On the other hand, if $\int_0^A k(a)\exp\left(-\int_0^a \mu(s)ds\right)da < 1$ then we define $x(t) = \exp(pt)v(t)$ for all $t \geq -A$, where $p > 0$. It follows that $x(t) = \int_0^A k(a) \exp\left(pa - \int_0^a \mu(s)ds\right) x(t-a) da$ for $t \geq 0$ and that $\int_0^A k(a) \exp\left(pa - \int_0^a \mu(s)ds\right) da \geq 1$ for $p > 0$ sufficiently large.

Another direct consequence of Lemma 4.1 and Lemma 3.2 is that if $k \in PC^1([0,A]; \Re_+)$, then the quantity $\frac{f(t,a)}{f^*(a)\Pi(f[t])} - 1$ appearing in the right hand side of the transformation (2.19) is only a function of $t-a$ (and thus (2.19) is a valid transformation). Indeed, it is straightforward to verify that for every piecewise continuous function $D: \Re_+ \to [D_{\min}, D_{\max}]$ and for every $f_0 \in PC^1([0,A];(0,+\infty))$ with $f_0(0) = \int_0^A k(a) f_0(a) da$, the solution of (2.1), (2.2) with $f(0,a) = f_0(a)$ for $a \in [0, A]$, corresponding to input $D: \Re_+ \to [D_{\min}, D_{\max}]$ satisfies $f(t,a) = z(t,a) \exp\left(-\int_0^t D(s)ds\right)$ for all $(t,a) \in \Re_+ \times [0,A]$, where $z: [0,+\infty) \times [0,A] \to \Re$ is the solution of (3.1), (3.2) with same initial condition $z(0,a) = f_0(a)$ for $a \in [0,A]$. Using (2.4), (3.3) and equation $f(t,a) = z(t,a) \exp\left(-\int_0^t D(s)ds\right)$, we get:

$$\frac{f(t,a)}{f^*(a)} = M^{-1} \exp\left(D^* a - \int_0^t D(s)ds\right) v(t-a), \text{ for all } (t,a) \in \Re_+ \times [0,A]$$

Using (2.4), (2.8), (3.3), equation $f(t,a) = z(t,a) \exp\left(-\int_0^t D(s)ds\right)$ and definition (2.9), we get:

$$M \int_0^A w\tilde{k}(w) dw \, \Pi(f[t]) = \exp\left(-\int_0^t D(s)ds\right) \int_0^A v(t-a) \exp\left(D^* a\right) \left(\int_a^A \tilde{k}(s)ds\right) da$$

$$= \exp\left(D^* t - \int_0^t D(s)ds\right) \int_{t-A}^t v(l) \exp(-D^* l) \left(\int_{t-l}^A \tilde{k}(s)ds\right) dl$$

Since $v(t) > 0$ for all $t \geq -A$ (a consequence of (3.3) and the conclusion of the previous paragraph), the above equation implies that $\Pi(f[t]) > 0$ for all $t \geq 0$. Combining the two above equations, we get:

$$\frac{f(t,a)}{f^*(a)\Pi(f[t])} = \frac{v(t-a)\exp(-D^*(t-a)) \int_0^A w\tilde{k}(w) dw}{\int_{t-A}^t v(l) \exp(-D^* l) \left(\int_{t-l}^A \tilde{k}(s)ds\right) dl}, \text{ for all } (t,a) \in \Re_+ \times [0,A]$$

where $\tilde{k}(a) = k(a) \exp\left(-D^* a - \int_0^a \mu(s)ds\right)$ for $a \in [0,A]$. Notice that (3.4) implies that $\frac{d}{dt} \int_{t-A}^t v(l) \exp(-D^* l) \left(\int_{t-l}^A \tilde{k}(s)ds\right) dl = 0$ for all $t \geq 0$. Indeed, we have for all $t \geq 0$:

$$\frac{d}{dt} \int_{t-A}^t v(l) \exp(-D^* l) \left(\int_{t-l}^A \tilde{k}(s)ds\right) dl = v(t) \exp(-D^* t) \left(\int_0^A \tilde{k}(s)ds\right) - \int_{t-A}^t v(l) \tilde{k}(t-l) \exp(-D^* l) dl$$

$$= v(t) \exp(-D^* t) \left(\int_0^A \tilde{k}(s)ds\right) - \int_0^A v(t-a) \tilde{k}(a) \exp(-D^*(t-a)) da$$

Using definition $\tilde{k}(a) = k(a)\exp\left(-D^*a - \int_0^a \mu(s)ds\right)$ for $a \in [0, A]$ and the fact that $\int_0^A \tilde{k}(s)ds = 1$ (a consequence of (2.3)), we get for all $t \geq 0$:

$$\frac{d}{dt}\int_{t-A}^{t} v(l)\exp(-D^*l)\left(\int_{t-l}^{A}\tilde{k}(s)ds\right)dl = \exp(-D^*t)\left(v(t) - \int_0^A v(t-a)k(a)\exp\left(-\int_0^a \mu(s)ds\right)da\right)$$

which combined with (3.4) gives $\frac{d}{dt}\int_{t-A}^{t} v(l)\exp(-D^*l)\left(\int_{t-l}^{A}\tilde{k}(s)ds\right)dl = 0$ for all $t \geq 0$. Therefore, using the fact that $v(-a) = \exp\left(\int_0^a \mu(s)ds\right)z(0,a) = \exp\left(\int_0^a \mu(s)ds\right)f_0(a)$, for all $a \in (0, A]$, we get

$$\int_{t-A}^{t} v(l)\exp(-D^*l)\left(\int_{t-l}^{A}\tilde{k}(s)ds\right)dl = \int_{-A}^{0} v(l)\exp(-D^*l)\left(\int_{-l}^{A}\tilde{k}(s)ds\right)dl$$
$$= \int_0^A f_0(w)\exp\left(D^*w + \int_0^w \mu(s)ds\right)\left(\int_w^A \tilde{k}(s)ds\right)dw$$
, for all $t \geq 0$

Consequently, the quantity $\frac{f(t,a)}{f^*(a)\Pi(f[t])} - 1$ is a function only of $t-a$, since we have:

$$\frac{f(t,a)}{f^*(a)\Pi(f[t])} = \frac{\int_0^A w\tilde{k}(w)dw}{\int_0^A f_0(w)\exp\left(D^*w + \int_0^w \mu(s)ds\right)\left(\int_w^A \tilde{k}(s)ds\right)dw} v(t-a)\exp(-D^*(t-a)), \text{ for all } (t,a) \in \Re_+ \times [0, A]$$

## *4.III. The Strong Ergodic Theorem in terms of IDEs*

Next, we state the strong ergodic theorem (Theorem 3.3) in terms of the IDE (4.1). To this goal, we define the operator

$$G: C^0([-A, 0]; \Re) \to C^0([0, A]; \Re)$$

for every $v \in C^0([-A, 0]; \Re)$ by the relation $(Gv)(a) = v(-a)$ for all $a \in [0, A]$.

If $\mu \in C^0([0, A]; \Re_+)$, $k \in PC^1([0, A]; \Re_+)$ satisfy (2.3) for certain $D^* \in \Re$, then it follows from Lemma 3.2 and Theorem 3.3 that there exist constants $\varepsilon > 0$, $K \geq 1$ such that for every $z_0 \in PC^1([0, A]; \Re)$ satisfying $z_0(0) = \int_0^A k(a)z_0(a)da$, the unique solution of the IDE (3.4) with initial condition $v(-a) = \exp\left(\int_0^a \mu(s)ds\right)z_0(a)$ for all $a \in [0, A]$ satisfies for all $t \geq 0$ the following estimate:

$$\int_0^A |v(t-a) - \exp(D^*(t-a))\Phi(z_0)|da \leq K\exp((D^* - \varepsilon)t)\int_0^A |v(-a)|da \tag{4.3}$$

The above property can be rephrased without any reference to the PDE: for every $\bar{k} \in PC^1([0, A]; \Re_+)$ with $1 = \int_0^A \bar{k}(a)\exp(-D^*a)da$ there exist constants $\varepsilon > 0$, $K \geq 1$ such that for every $v_0 \in C^0([-A, 0]; \Re)$ with $v_0(0) = \int_0^A \bar{k}(a)v_0(-a)da$ and $(Gv_0) \in PC^1([0, A]; \Re)$, the unique solution of the IDE $v(t) = \int_0^A \bar{k}(a) v(t-a) da$ with initial condition $v(-a) = v_0(-a)$, for all $a \in [0, A]$ satisfies (4.3) for all $t \geq 0$, with $z_0(a) = v(-a) \exp\left(-\int_0^a \mu(s) ds\right)$ for $a \in [0, A]$.

Using the transformation $x(t) = \exp(-D^* t) v(t)$, for all $t \geq -A$, we obtain a "one-to-one" mapping of solutions of the IDE $v(t) = \int_0^A \bar{k}(a) v(t-a) da$ to the solutions of the IDE (4.1) with $\varphi(a) := \bar{k}(a) \exp(-D^* a)$ for all $a \in [0, A]$. Moreover, estimate (4.3) implies the following estimate for all $t \geq 0$:

$$\int_0^A |x(t-a) - P(x_0)| da \leq K \exp(-\varepsilon t) \exp(D^* A) \int_0^A |x(-a)| da$$

where the functional $P : C^0([-A, 0]; \Re) \to \Re$, is defined by means of the equation

$$P(x) = r \int_0^A x(-a) \int_a^A \varphi(s) ds \, da \tag{4.4}$$

and $r := \left(\int_0^A a \varphi(a) da\right)^{-1}$. The functional $P : C^0([-A, 0]; \Re) \to \Re$ is found by substituting $z_0(a) = x(-a) \exp\left(-D^* a - \int_0^a \mu(s) ds\right)$ (for $a \in [0, A]$) in the functional $\Phi : L^1([0, A]; \Re) \to \Re$ defined by (3.7).

Therefore, we are in a position to conclude that the following property holds: for every $\varphi \in PC^1([0, A]; \Re_+)$ with $1 = \int_0^A \varphi(a) da$ there exist constants $\tilde{K}, \varepsilon > 0$ such that for every $x_0 \in C^0([-A, 0]; \Re)$ with $x_0(0) = \int_0^A \varphi(a) x_0(-a) da$ and $(Gx_0) \in PC^1([0, A]; \Re)$, the unique solution of the IDE (4.1) with initial condition $x(-a) = x_0(-a)$, for all $a \in [0, A]$ satisfies the following estimate for all $t \geq 0$

$$\int_0^A |x(t-a) - P(x_0)| da \leq \tilde{K} \exp(-\varepsilon t) \int_0^A |x(-a)| da \tag{4.5}$$

Using this property, we obtain the following corollary, which is a restatement of the strong ergodic theorem (Theorem 3.3) in terms of IDEs and the $L^\infty$ norm (instead of the $L^1$ norm). Recall that $X = \left\{ x \in C^0([-A, 0]; \Re) : x(0) = \int_0^A \varphi(a) x(-a) da \right\}$.

**Corollary 4.2 (The strong ergodic theorem in terms of IDEs):** *Suppose that* $\varphi \in PC^1([0, A]; \Re_+)$ *with* $1 = \int_0^A \varphi(a) da$. *Then there exist constants* $M, \sigma > 0$ *such that for every* $x_0 \in X$ *with* $(Gx_0) \in PC^1([0, A]; \Re)$, *the unique solution of the IDE (4.1) with initial condition* $x(-a) = x_0(-a)$ *for all* $a \in [0, A]$ *satisfies the following estimate for all* $t \geq 0$:

$$\max_{-A \leq \theta \leq 0} \left(|x(t+\theta) - P(x_0)|\right) \leq M \exp(-\sigma t) \max_{-A \leq a \leq 0} \left(|x_0(a)|\right) \tag{4.6}$$

## 4.IV. The Construction of Lyapunov Functionals

The problem with Corollary 4.2 is that it does not provide a Lyapunov-like functional which can allow the derivation of the important property (4.6). Moreover, it does not provide information about the magnitude of the constant $\sigma > 0$. In order to construct a Lyapunov-like functional and obtain information about the magnitude of the constant $\sigma > 0$, we need some technical results. The first result deals with the exponential stability of the zero solution for (4.1). Notice that the proof of the exponential stability property is made by means of a Lyapunov functional.

**Lemma 4.3 (Lyapunov functional for the general case):** *Suppose that $\int_0^A |\varphi(a)| da < 1$. Then $0 \in X$ is globally exponentially stable for (4.1). Moreover, the functional $V: X \to \Re_+$ defined by $V(x) := \max_{a \in [0,A]} \left( \exp(-\sigma a) |x(-a)| \right)$, where $\sigma > 0$ is a constant that satisfies $\int_0^A |\varphi(a)| \exp(\sigma a) da < 1$, satisfies the differential inequality:*

$$\limsup_{h \to 0^+} \left( h^{-1} (V(x_{t+h}) - V(x_t)) \right) \leq -\sigma V(x_t), \text{ for all } t \geq 0 \quad (4.7)$$

*for every solution of (4.1).*

Lemma 4.3 is useful because we next construct Lyapunov functionals of the form used in Lemma 4.3. However, we are mostly interested in kernels $\varphi \in C^0([0, A]; \Re)$ with non-negative values that satisfy $\int_0^A \varphi(a) da = 1$. We show next that even for this specific case, it is possible to construct a Lyapunov functional on an invariant subspace of the state space $X = \left\{ x \in C^0([-A, 0]; \Re) : x(0) = \int_0^A \varphi(a) x(-a) da \right\}$. We next introduce a technical assumption.

**(H1)** *The function $\varphi \in C^0([0, A]; \Re)$ satisfies $\varphi(a) \geq 0$ for all $a \in [0, A]$ and $\int_0^A \varphi(a) da = 1$. Moreover, there exists $\lambda > 0$ such that $\int_0^A \left| \varphi(a) - r\lambda \int_a^A \varphi(s) ds \right| da < 1$, where $r := \left( \int_0^A a\varphi(a) da \right)^{-1}$.*

The following result provides the construction of a Lyapunov functional for system (4.1) under assumption (H1).

**Theorem 4.4 (Lyapunov functional for linear IDEs with special kernels):** *Consider system (4.1), where $\varphi \in C^0([0, A]; \Re_+)$ satisfies assumption (H1). Let $\lambda > 0$ be a real number for which $\int_0^A \left| \varphi(a) - r\lambda \int_a^A \varphi(s) ds \right| da < 1$, where $r := \left( \int_0^A a\varphi(a) da \right)^{-1}$. Define the functional $V: X \to \Re_+$ by means of the equation:*

$$V(x) := \max_{a \in [0,A]} \left( \exp(-\sigma a) |x(-a) - P(x)| \right) \quad (4.8)$$

*where $\sigma > 0$ is a real number for which $\int_0^A \left| \varphi(a) - r\lambda \int_a^A \varphi(s) ds \right| \exp(\sigma a) da < 1$ and $P: X \to \Re$ is the functional defined by (4.4). Then the following relations hold*

$$P(x_t) = P(x_0) \text{, for all } t \geq 0 \tag{4.9}$$

$$\limsup_{h \to 0^+}\left(h^{-1}(V(x_{t+h}) - V(x_t))\right) \leq -\sigma V(x_t)\text{, for all } t \geq 0 \tag{4.10}$$

*for every solution of (4.1).*

**Remark 4.5:** Theorem 4.4 is a Lyapunov-like version of the scalar strong ergodic theorem (compare with Corollary 4.2) for kernels that satisfy assumption (H1). Corollary 4.2 does not allow us to estimate the magnitude of the constant $\sigma > 0$ that determines the convergence rate. On the other hand, Theorem 4.4 allows us to estimate $\sigma > 0$: the Comparison Lemma on page 85 in [19] and differential inequality (4.10) guarantee that $V(x_t) \leq \exp(-\sigma t)V(x_0)$ for all $t \geq 0$ and for every solution of (4.1). Using (4.9), definition (4.8) and the previous estimate, we can guarantee that

$$\max_{a \in [0,A]}(|x(t-a) - P(x_t)|) = \max_{a \in [0,A]}(|x(t-a) - P(x_0)|) \leq \exp(-\sigma(t-A))\max_{a \in [0,A]}(|x(-a) - P(x_0)|)\text{, for all } t \geq 0$$

Therefore, bounds for $\sigma > 0$ can be computed in a straightforward way using the inequality $\int_0^A \left|\varphi(a) - r\lambda\int_a^A \varphi(s)ds\right|\exp(\sigma a)da < 1$ (e.g., an allowable value for $\sigma > 0$ is $-A^{-1}\ln\left(\int_0^A \left|\varphi(a) - r\lambda\int_a^A \varphi(s)ds\right|da\right)$).

Moreover, Corollary 4.2 does not provide a Lyapunov-like functional for equation (4.1). However, the cost of these features is the loss of generality: while Corollary 4.2 holds for all kernels $\varphi \in PC^1([0,A];\Re_+)$ that satisfy $\varphi(a) \geq 0$ for all $a \in [0,A]$ and $\int_0^A \varphi(a)da = 1$, Theorem 4.4 holds only for kernels that satisfy Assumption (H1).

Theorem 4.4 can allow us to guarantee exponential stability for the zero solution of (4.1), when the state evolves in certain invariant subsets of the state space. This is shown in the following result.

**Corollary 4.6 (Lyapunov functional for linear IDEs on invariant sets):** *Consider system (4.1), where $\varphi \in C^0([0,A];\Re_+)$ satisfies assumption (H1). Let $\lambda > 0$ be a real number for which $\int_0^A \left|\varphi(a) - r\lambda\int_a^A \varphi(s)ds\right|da < 1$, where $r := \left(\int_0^A a\varphi(a)da\right)^{-1}$. Let $P : X \to \Re$ be the functional defined by (4.4). Define the functional $W : X \to \Re_+$ by means of the equation:*

$$W(x) := \max_{a \in [0,A]}\left(\exp(-\sigma a)|x(-a)|\right) \tag{4.11}$$

*where $\sigma > 0$ is a real number for which $\int_0^A \left|\varphi(a) - r\lambda\int_a^A \varphi(s)ds\right|\exp(\sigma a)da < 1$. Let $S \subseteq X$ be a positively invariant set for system (4.1) and let $C : S' \to [\kappa, +\infty)$, where $\kappa > 0$ is a constant and $S' \subseteq C^0([-A,0];\Re)$ is an open set with $S \subset S'$, be a continuous functional that satisfies*

$$\limsup_{h \to 0^+}\left(h^{-1}(C(x_{t+h}) - C(x_t))\right) \leq 0 \tag{4.12}$$

*for every $t \geq 0$ and for every solution $x(t) \in \Re$ of (4.1) with $x_t \in S$. Then for every $x_0 \in S$ with $P(x_0) = 0$ and for every $b \in K_\infty \cap C^1([0,+\infty);\Re_+)$, the following hold for the solution $x(t) \in \Re$ of (4.1) with initial condition $x_0 \in S$:*

$$\limsup_{h \to 0^+}\left(h^{-1}(C(x_{t+h})b(W(x_{t+h})) - C(x_t)b(W(x_t)))\right) \leq -\sigma C(x_t)b'(W(x_t))W(x_t)\text{, for all } t \geq 0 \tag{4.13}$$

$$P(x_t) = 0\text{, for all } t \geq 0 \tag{4.14}$$

**Remark 4.7: (a)** The differential inequality (4.12) is equivalent to the assumption that the mapping $t \to C(x_t)$ is non-increasing.

**(b)** Using assumption (H1) and Lemma 4.1, we can guarantee that the mapping $t \to g_1(x_t)$ is non-decreasing and that the mapping $t \to g_2(x_t)$ is non-increasing for every solution of (4.1), where $g_1, g_2 : C^0([-A,0]; \Re) \to \Re$ are the continuous functionals $g_1(x) := \min_{a \in [0,A]}(x(-a))$ and $g_2(x) := \max_{a \in [0,A]}(x(-a))$. Indeed, Lemma 4.1 implies that for every solution of (4.1) it holds that:

$$g_1(x_0) \leq g_1(x_t) \leq g_2(x_t) \leq g_2(x_0), \text{ for all } t \geq 0 \tag{4.15}$$

Consequently, any set $S \subseteq X$ of the form $S = \{x \in X : \min_{a \in [0,A]}(x(-a)) > c_1\}$, $S = \{x \in X : \max_{a \in [0,A]}(x(-a)) < c_2\}$, $S = \{x \in X : c_1 < \min_{a \in [0,A]}(x(-a)) \leq \max_{a \in [0,A]}(x(-a)) < c_2\}$, where $c_1 < c_2$ are constants, is a positively invariant set for (4.1). Moreover, using the semigroup property for the solution of (4.1) and (4.15), we get

$$g_1(x_{t_2}) \leq g_1(x_{t_1}) \leq g_2(x_{t_1}) \leq g_2(x_{t_2}), \text{ for all } t_1 \geq t_2 \geq 0$$

The above inequality shows that the mapping $t \to g_1(x_t)$ is non-decreasing and that the mapping $t \to g_2(x_t)$ is non-increasing for every solution of (4.1).

## 5. Proofs of Main Results

We next turn our attention to the proof of Theorem 2.1. Throughout this section we use the notation:

$$q(x) := \min(D_{\max} - D^*, \max(-(D^* - D_{\min}), x)), \text{ for all } x \in \Re \tag{5.1}$$

Notice that $q(x)$ is a non-decreasing function, which satisfies the equation:

$$q(x) = sat(D^* + x) - D^*, \text{ for all } x \in \Re \tag{5.2}$$

Equation (5.2) and the fact $D^* \in (D_{\min}, D_{\max})$ imply the inequality

$$|q(x)| \leq \max(D_{\max} - D^*, D^* - D_{\min}), \text{ for all } x \in \Re \tag{5.3}$$

We also notice that the inequality

$$xq(x) \geq \min(1, D_{\max} - D^*, D^* - D_{\min}) \frac{x^2}{1+|x|} \text{ for all } x \in \Re \tag{5.4}$$

holds. Indeed, inequality (5.4) can be derived by using definition (5.1) and distinguishing three cases: (i) $-(D^* - D_{\min}) \leq x \leq D_{\max} - D^*$, (ii) $x > D_{\max} - D^*$, and (iii) $-(D^* - D_{\min}) > x$. For case (i) we get from (5.1) $xq(x) = x^2$ and since $1 \geq \frac{\min(1, D_{\max} - D^*, D^* - D_{\min})}{1+|x|}$, we conclude that (5.4) holds in this case. For case (ii) we get from (5.1) $xq(x) = (D_{\max} - D^*)|x|$ and since $D_{\max} - D^* \geq \min(1, D_{\max} - D^*, D^* - D_{\min})$, $|x| \geq \frac{x^2}{1+|x|}$, we conclude that (5.4) holds in this case. The proof is similar for case (iii).

The proof of Theorem 2.1 is based on the transformation shown in Figure 1 and on the following lemmas. Their proofs can be found in the Appendix.

**Lemma 5.1:** *Consider the control system*

$$\dot{\eta}(t) = D^* - D(t) \quad , \quad \psi(t) = \int_0^A \tilde{k}(a)\psi(t-a)da \quad , \quad (\eta(t),\psi(t)) \in \Re^2 \tag{5.5}$$

*where $A>0$ is a constant, $D^* \in (D_{min}, D_{max})$ is a constant, $D_{max} > D_{min} > 0$ are constants, $\tilde{k} \in C^0([0,A];\Re_+)$ satisfies assumption (A) and $\int_0^A \tilde{k}(a)da = 1$. The control system (5.5) is defined on the set $\Re \times S$, where*

$$S = \tilde{S} \cap \{\psi \in C^0([-A,0];\Re): P(\psi) = 0\}, \quad \tilde{S} = \left\{\psi \in C^0([-A,0];(-1,+\infty)): \psi(0) = \int_0^A \tilde{k}(a)\psi(-a)da\right\}$$

*and $P(\psi)$ is the linear functional $P(\psi) := r\int_0^A \psi(-a)\int_a^A \tilde{k}(s)dsda$ with $r := \left(\int_0^A a\tilde{k}(a)da\right)^{-1}$. The measured output of system (5.5) is given by the equation*

$$Y(t) = \eta(t) + \ln\left(1 + \int_0^A g(a)\psi(t-a)da\right) \tag{5.6}$$

*where the function $g \in C^0([0,A];\Re)$ satisfies $g(a) \geq 0$ for all $a \in [0,A]$ and $\int_0^A g(a)da = 1$. Consider the closed-loop system (5.5) with the dynamic feedback law given by*

$$\dot{z}_1(t) = z_2(t) - D(t) - l_1(z_1(t) - Y(t)) \quad , \quad \dot{z}_2(t) = -l_2(z_1(t) - Y(t)) \quad , \quad z(t) = (z_1(t), z_2(t))' \in \Re^2 \tag{5.7}$$

*and*

$$D(t) = sat(z_2(t) + \gamma Y(t)) \tag{5.8}$$

*where $l_1, l_2, \gamma > 0$ are constants. Let $p_1, p_2 > 0$ be a pair of constants satisfying $(2 + l_1 p_1 - 2l_2 p_2)^2 < 8l_1 p_1 - 4l_2 p_1^2$, $p_1^2 < 4p_2$. Then there exist sufficiently large constants $M, G > 0$ and sufficiently small constants $L, \sigma > 0$ such that for every constant $\beta \geq 0$ and for every $z_0 \in \Re^2$, $(\eta_0, \psi_0) \in \Re \times S$ the solution $(z(t), \eta(t), \psi_t) \in \Re^2 \times \Re \times S$ of the closed-loop system (5.5), (5.6) with (5.7), (5.8) and initial condition $(\eta(0), z(0)) = (\eta_0, z_0)$, $\psi(-a) = \psi_0(-a)$ for all $a \in [0,A]$, is unique, exists for all $t \geq 0$ and satisfies the differential inequality*

$$\limsup_{h \to 0^+}\left(h^{-1}\left(V(\eta(t+h), z(t+h), \psi_{t+h}) - V(\eta(t), z(t), \psi_t)\right)\right) \leq -L\frac{V(\eta(t), z(t), \psi_t)}{1 + \sqrt{V(\eta(t), z(t), \psi_t)}}, \text{ for all } t \geq 0$$

$$\tag{5.9}$$

*where $V: \Re \times \Re^2 \times C^0([-A,0];(-1,+\infty)) \to \Re_+$ is the continuous functional defined by:*

$$V(\eta, z, \psi) := \eta^2 + G\sqrt{Q(z, \eta, \psi)} + \beta Q(z, \eta, \psi) \tag{5.10}$$

$$Q(z, \eta, \psi) := (z_1 - \eta)^2 - p_1(z_1 - \eta)(z_2 - D^*) + p_2(z_2 - D^*)^2 + \frac{M}{2}\left(\frac{\max_{a \in [0,A]}(\exp(-\sigma a)|\psi(-a)|)}{1 + \min\left(0, \min_{a \in [0,A]}(\psi(-a))\right)}\right)^2 \tag{5.11}$$

**Lemma 5.2:** *Suppose that there exists a constant $L > 0$ such that the continuous function $\varphi: \Re_+ \to \Re_+$ satisfies:*

$$\limsup_{h \to 0^+}\left(h^{-1}(\varphi(t+h) - \varphi(t))\right) \leq -L\frac{\varphi(t)}{1 + \sqrt{\varphi(t)}}, \text{ for all } t \geq 0 \tag{5.12}$$

*Then the following estimate holds:*

$$\varphi(t) \leq \varphi(0)\exp(\max(0, \varphi(0) - 1))\exp\left(-\frac{L}{2}t\right), \text{ for all } t \geq 0 \tag{5.13}$$

We are now ready to provide the proof of Theorem 2.1.

**Proof of Theorem 2.1:** Define

$$\psi_0(-a) = \frac{f_0(a)}{f^*(a)\Pi(f_0)} - 1, \text{ for all } a \in [0, A] \quad (5.14)$$

$$\eta_0 = \ln(\Pi(f_0)) \quad (5.15)$$

It is straightforward to verify (using definitions (2.8), (5.14), equation (2.4) and the fact that $\tilde{k}(a) = k(a)\exp\left(-D^*a - \int_0^a \mu(s)ds\right)$ for all $a \in [0, A]$) that $P(\psi_0) = 0$, where $P(\psi) = r\int_0^A \psi(-a)\int_a^A \tilde{k}(s)ds\,da$ with $r := \left(\int_0^A a\tilde{k}(a)da\right)^{-1}$. Define

$$g(a) := p(a)f^*(a)/y^*, \text{ for all } a \in [0, A] \quad (5.16)$$

Notice that (2.4) and the fact that $M = y^*\left(\int_0^A p(a)\exp\left(-D^*a - \int_0^a \mu(s)ds\right)da\right)^{-1}$ imply that the function $g \in C^0([0, A]; \Re)$ satisfies $g(a) \geq 0$ for all $a \in [0, A]$ and $\int_0^A g(a)da = 1$.

Next consider the solution $(z(t), \eta(t), \psi_t) \in \Re^2 \times \Re \times S$ of the closed-loop system (5.5) with (5.7), (5.8) and initial condition $(\eta(0), z(0)) = (\eta_0, z_0)$, $\psi(-a) = \psi_0(-a)$ for all $a \in [0, A]$. Lemma 5.1 guarantees that the solution $(z(t), \eta(t), \psi_t) \in \Re^2 \times \Re \times S$ of the closed-loop system (5.5) with (5.7), (5.8) exists for all $t \geq 0$. The solution of the IDE $\psi(t) = \int_0^A \tilde{k}(a)\psi(t-a)da$ is $C^1$ on $(0, +\infty)$ since it coincides with the solution of the delay differential equation $\dot{\psi}(t) = \tilde{k}(0)\psi(t) - \tilde{k}(A)\psi(t-A) + \int_0^A \tilde{k}'(a)\psi(t-a)da$ with the same initial condition. Therefore by virtue of (5.14), the function defined by:

$$f(t, a) = (1 + \psi(t-a))f^*(a)\exp(\eta(t)) \text{ for } (t, a) \in \Re_+ \times [0, A] \quad (5.17)$$

is continuous and $f \in C^1(D_f; (0, +\infty))$, where $D_f = \{(t, a) \in (0, +\infty) \times (0, A) : (a-t) \notin B \cup \{0, A\}\}$ and $B \subseteq (0, A)$ is the finite (possibly empty) set where the derivative of $f_0 \in \tilde{X}$ is not defined or is not continuous. Since $\psi_t \in S$, it follows that $f[t] \in \tilde{X}$ for all $t \geq 0$, where $(f[t])(a) = f(t, a)$ for $a \in [0, A]$. Using (5.15), (5.17) and (5.5) we conclude that

$$\eta(t) = \ln(\Pi(f[t])), \text{ for all } t \geq 0 \quad (5.18)$$

Moreover, using (5.5), (5.6), (5.7), (5.8), (5.16), (5.17) and (5.18), we conclude that equations (2.5), (2.6), (2.7) hold for all $t \geq 0$ and equation (2.1) holds for all $(t, a) \in D_f$. Finally, by virtue of (5.14), (5.15) and (5.17), it follows that $z(0) = z_0$, $f(0, a) = f_0(a)$ for all $a \in [0, A]$.

Using (2.11), (2.12), (5.17), (5.18) we conclude that $V(\eta(t), z(t), \psi_t) = W(z(t), f[t])$ for all $t \geq 0$. Therefore, the differential inequality (5.9) implies the differential inequality (2.13).

Lemma 5.2 in conjunction with inequality (5.9), implies that the following estimate holds:

$$V(\eta(t), z(t), \psi_t) \leq V(\eta_0, z_0, \psi_0)\exp(\max(0, V(\eta_0, z_0, \psi_0) - 1))\exp\left(-\frac{L}{2}t\right), \text{ for all } t \geq 0 \quad (5.19)$$

Since $p_1, p_2 > 0$ is a pair of constants with $p_1^2 < 4p_2$, it follows that the quadratic form $A(e) := e_1^2 - p_1 e_1 e_2 + p_2 e_2^2$ is positive definite. Therefore, there exist constants $K_2 \geq K_1 > 0$, such that $K_1 |e|^2 \leq A(e) = e_1^2 - p_1 e_1 e_2 + p_2 e_2^2 \leq K_2 |e|^2$, for all $e \in \Re^2$. Using the previous inequality, (5.10), (5.11) and (5.19), we obtain the following estimates for all $t \geq 0$:

$$\eta^2(t) \leq V(\eta_0, z_0, \psi_0) \exp(\max(0, V(\eta_0, z_0, \psi_0) - 1)) \exp\left(-\frac{L}{2} t\right) \quad (5.20)$$

$$G \sqrt{\frac{M}{2}} \frac{\max_{a \in [0,A]} \left(\exp(-\sigma a) |\psi_t(-a)|\right)}{1 + \min\left(0, \min_{a \in [0,A]} (\psi_t(-a))\right)} \leq V(\eta_0, z_0, \psi_0) \exp(\max(0, V(\eta_0, z_0, \psi_0) - 1)) \exp\left(-\frac{L}{2} t\right) \quad (5.21)$$

$$G \sqrt{K_1} \max\left(|z_1(t) - \eta(t)|, |z_2(t) - D^*|\right) \leq V(\eta_0, z_0, \psi_0) \exp(\max(0, V(\eta_0, z_0, \psi_0) - 1)) \exp\left(-\frac{L}{2} t\right) \quad (5.22)$$

Estimates (5.20), (5.21), (5.22) and the fact that $M, G > 0$ are sufficiently large constants (with $G \sqrt{K_1} \geq 3$ and $G \sqrt{\frac{M}{2}} \geq 3$), imply the following estimate for all $t \geq 0$:

$$|\eta(t)| + \frac{\max_{a \in [0,A]} \left(\exp(-\sigma a) |\psi_t(-a)|\right)}{1 + \min\left(0, \min_{a \in [0,A]} (\psi_t(-a))\right)} + |z_1(t)| + |z_2(t) - D^*|$$

$$\leq \left(V(\eta_0, z_0, \psi_0) + \sqrt{V(\eta_0, z_0, \psi_0)}\right) \exp(\max(0, V(\eta_0, z_0, \psi_0) - 1)) \exp\left(-\frac{L}{4} t\right) \quad (5.23)$$

Using (5.17) and the fact that $|\ln(1+x)| \leq \frac{|x|}{1 + \min(x,0)}$ for all $x > -1$, we get the estimate for all $a \in [0, A]$:

$$\left|\ln\left(\frac{f(t,a)}{f^*(a)}\right)\right| \leq |\eta(t)| + \frac{|\psi(t-a)|}{1 + \min(0, \psi(t-a))} \leq |\eta(t)| + \frac{\exp(\sigma A) \max_{s \in [0,A]} \left(\exp(-\sigma s) |\psi_t(-s)|\right)}{1 + \min\left(0, \min_{s \in [0,A]} (\psi_t(-s))\right)}$$

which implies the following estimate for all $t \geq 0$:

$$\max_{a \in [0,A]} \left(\left|\ln\left(\frac{f(t,a)}{f^*(a)}\right)\right|\right) \leq |\eta(t)| + \frac{\exp(\sigma A) \max_{a \in [0,A]} \left(\exp(-\sigma a) |\psi_t(-a)|\right)}{1 + \min\left(0, \min_{a \in [0,A]} (\psi_t(-a))\right)} \quad (5.24)$$

Combining (5.23) and (5.24), we obtain the following estimate for all $t \geq 0$:

$$\max_{a \in [0,A]} \left(\left|\ln\left(\frac{f(t,a)}{f^*(a)}\right)\right|\right) + |z_1(t)| + |z_2(t) - D^*|$$

$$\leq \exp(\sigma A) \left(V(\eta_0, z_0, \psi_0) + \sqrt{V(\eta_0, z_0, \psi_0)}\right) \exp(\max(0, V(\eta_0, z_0, \psi_0) - 1)) \exp\left(-\frac{L}{4} t\right) \quad (5.25)$$

Taking into account (5.25), we conclude that the validity of (2.10) relies on showing that there exists a function $b \in K_\infty$ that satisfies the following inequality for all $f_0 \in \tilde{X}$ and $z_0 \in \Re^2$ and for $(\eta_0, \psi_0) \in \Re \times S$ satisfying (5.14), (5.15):

$$V(\eta_0, z_0, \psi_0) \leq b\left(\max_{a \in [0,A]} \left(\left|\ln\left(\frac{f_0(a)}{f^*(a)}\right)\right|\right) + |z_{1,0}| + |z_{2,0} - D^*|\right) \quad (5.26)$$

In order to show (5.26), and taking into account definitions (5.10), (5.11), it suffices to show there exist functions $b_1, b_2 \in K_\infty$ for which the following inequalities hold for all $f_0 \in \tilde{X}$ and for $(\eta_0, \psi_0) \in \Re \times S$ satisfying (5.14), (5.15):

$$|\eta_0| \leq b_1 \left( \max_{a\in[0,A]} \left( \left| \ln\left( \frac{f_0(a)}{f^*(a)} \right) \right| \right) \right), \quad \frac{\max_{a\in[0,A]}\left(\exp(-\sigma a)|\psi_0(-a)|\right)}{1+\min\left(0, \min_{a\in[0,A]}(\psi_0(-a))\right)} \leq b_2 \left( \max_{a\in[0,A]} \left( \left| \ln\left( \frac{f_0(a)}{f^*(a)} \right) \right| \right) \right) \quad (5.27)$$

In what follows we are using repeatedly the notation $v = \max_{a\in[0,A]}\left(\left|\ln\left(\frac{f_0(a)}{f^*(a)}\right)\right|\right)$ and the fact that

$$\exp(-v) \leq \min_{a\in[0,A]}\left(\frac{f_0(a)}{f^*(a)}\right) \leq \max_{a\in[0,A]}\left(\frac{f_0(a)}{f^*(a)}\right) \leq \exp(v) \quad (5.28)$$

Inequality (5.28) follows from the definition $v = \max_{a\in[0,A]}\left(\left|\ln\left(\frac{f_0(a)}{f^*(a)}\right)\right|\right)$ and the following implications

$$\left|\ln\left(\frac{f_0(a)}{f^*(a)}\right)\right| \leq v, \text{ for all } a \in [0, A] \Rightarrow$$

$$\Rightarrow -v \leq \ln\left(\frac{f_0(a)}{f^*(a)}\right) \leq v, \text{ for all } a \in [0, A] \Rightarrow$$

$$\Rightarrow \exp(-v) \leq \frac{f_0(a)}{f^*(a)} \leq \exp(v), \text{ for all } a \in [0, A]$$

Using (2.4) and (2.8) we get for all $f_0 \in \tilde{X}$:

$$\min_{a\in[0,A]}\left(\frac{f_0(a)}{f^*(a)}\right) \leq \Pi(f_0) \leq \max_{a\in[0,A]}\left(\frac{f_0(a)}{f^*(a)}\right) \quad (5.29)$$

Inequalities (5.29) is derived by means of definition (2.8), which directly implies

$$\Pi(f_0) \leq \frac{\int_0^A f^*(a)\left(\int_a^A k(s)\exp\left(\int_s^a (\mu(l)+D^*)dl\right)ds\right)da}{\int_0^A ak(a)f^*(a)da} \max_{a\in[0,A]}\left(\frac{f_0(a)}{f^*(a)}\right) = \Pi(f^*) \max_{a\in[0,A]}\left(\frac{f_0(a)}{f^*(a)}\right)$$

$$\Pi(f_0) \geq \frac{\int_0^A f^*(a)\left(\int_a^A k(s)\exp\left(\int_s^a (\mu(l)+D^*)dl\right)ds\right)da}{\int_0^A ak(a)f^*(a)da} \min_{a\in[0,A]}\left(\frac{f_0(a)}{f^*(a)}\right) = \Pi(f^*) \min_{a\in[0,A]}\left(\frac{f_0(a)}{f^*(a)}\right)$$

Moreover, by virtue of (2.4) and (2.8), we have $\Pi(f^*) = 1$, since

$$\Pi(f^*) = \frac{\int_0^A f^*(a)\left(\int_a^A k(s)\exp\left(\int_s^a (\mu(l)+D^*)dl\right)ds\right)da}{\int_0^A ak(a)f^*(a)da}$$

$$= \frac{M\int_0^A \exp\left(-D^*a - \int_0^a \mu(s)ds\right)\left(\int_a^A k(s)\exp\left(\int_s^a (\mu(l)+D^*)dl\right)ds\right)da}{M\int_0^A ak(a)\exp\left(-D^*a - \int_0^a \mu(s)ds\right)da}$$

$$= \frac{\int_0^A \left(\int_a^A k(s)\exp\left(-D^*s - \int_0^s \mu(l)dl\right)ds\right)da}{\int_0^A ak(a)\exp\left(-D^*a - \int_0^a \mu(s)ds\right)da} = 1$$

Notice that for the last equality above we have used integration by parts for the integral in the numerator. Combining (5.28), (5.29) and using (5.15) we get:

$$|\eta_0| \leq v \tag{5.30}$$

Consequently, (5.30) shows that the first inequality (5.27) holds with $b_1(s) = s$ for all $s \geq 0$. On the other hand, using (5.14), (5.15) we get for all $f_0 \in \tilde{X}$ and $a \in [0, A]$:

$$|\psi_0(-a)| \leq \exp(-\eta_0)\left|\frac{f_0(a)}{f^*(a)} - 1\right| + |\exp(-\eta_0) - 1| \text{ and } \psi_0(-a) \geq \exp(-\eta_0) \min_{a \in [0,A]}\left(\frac{f_0(a)}{f^*(a)}\right) - 1 \tag{5.31}$$

Using (5.28) and (5.30), we obtain for all $f_0 \in \tilde{X}$ and $a \in [0, A]$:

$$|\psi_0(-a)| \leq \exp(2v) - 1 \text{ and } \psi_0(-a) \geq \exp(-2v) - 1 \tag{5.32}$$

The following inequality is a direct consequence of (5.32):

$$\frac{\max_{a \in [0,A]}\left(\exp(-\sigma a)|\psi_0(-a)|\right)}{1 + \min\left(0, \min_{a \in [0,A]}(\psi_0(-a))\right)} \leq \exp(2v)(\exp(2v) - 1) \tag{5.33}$$

Consequently, (5.33) shows that the second inequality (5.27) holds with $b_2(s) = \exp(2s)(\exp(2s) - 1)$ for all $s \geq 0$. The proof is complete. ◁

The proof of Theorem 2.4 is based on the transformation shown in Figure 1 and on the following lemma. Its proof can be found in the Appendix.

**Lemma 5.3:** *Consider the control system (5.5) where $A > 0$ is a constant, $D^* \in (D_{\min}, D_{\max})$ is a constant, $D_{\max} > D_{\min} > 0$ are constants, $\tilde{k} \in C^0([0, A]; \Re_+)$ satisfies $\int_0^A \tilde{k}(a)da = 1$. The control system (5.5) is defined on the set $\Re \times S$, where*

$$S = \tilde{S} \cap \{\psi \in C^0([-A,0]; \Re) : P(\psi) = 0\}, \quad \tilde{S} = \left\{\psi \in C^0([-A,0]; (-1,+\infty)) : \psi(0) = \int_0^A \tilde{k}(a)\psi(-a)da\right\}$$

*and $P(\psi) := r\int_0^A \psi(-a)\int_a^A \tilde{k}(s)ds\,da$ with $r := \left(\int_0^A a\tilde{k}(a)da\right)^{-1}$ is a linear functional. The measured output of system (5.5) is given by the equation (5.6), where the function $g \in C^0([0, A]; \Re)$ satisfies $g(a) \geq 0$ for all $a \in [0, A]$ and $\int_0^A g(a)da = 1$. Consider the closed-loop system (5.5) with the dynamic feedback law given by*

$$D(t) = sat\left(D^* + T^{-1}Y(iT)\right), \text{ for all } t \in [iT, (i+1)T) \text{ and for all integers } i \geq 0 \tag{5.34}$$

*where $T > 0$ is a constant. Let $G: C^0([-A,0]; \Re) \to C^0([0, A]; \Re)$ be the operator defined by the relation $(Gv)(a) = v(-a)$ for all $a \in [0, A]$ for every $v \in C^0([-A,0]; \Re)$. Then there exist a constant $L > 0$ and a function $\kappa \in K_\infty$ such that for every $(\eta_0, \psi_0) \in \Re \times S$ with $(G\psi_0) \in PC^1([0, A]; \Re)$ the solution $(\eta(t), \psi_t) \in \Re \times S$ of the closed-loop system (5.5), (5.6) with (5.34) and initial condition $\eta(0) = \eta_0$, $\psi(-a) = \psi_0(-a)$ for all $a \in [0, A]$, is unique, exists for all $t \geq 0$ and satisfies the following estimate*

$$|\eta(t)| + \frac{\max_{-A \leq s \leq 0}(|\psi(t+s)|)}{1 + \min\left(0, \min_{-A \leq s \leq 0}(\psi(t+s))\right)} \leq \exp(-Lt)\kappa\left(|\eta_0| + \frac{\max_{-A \leq s \leq 0}(|\psi_0(s)|)}{1 + \min\left(0, \min_{-A \leq s \leq 0}(\psi_0(s))\right)}\right), \text{ for all } t \geq 0 \tag{5.35}$$

We are now ready to provide the proof of Theorem 2.4.

**Proof of Theorem 2.4:** Define $(\eta_0, \psi_0) \in \Re \times S$ by means of (5.14), (5.15). It is straightforward to verify (using definitions (2.8), (5.14), equation (2.4) and the fact that $\tilde{k}(a) = k(a)\exp\left(-D^*a - \int_0^a \mu(s)ds\right)$ for all $a \in [0, A]$) that $P(\psi_0) = 0$ and $(G\psi_0) \in PC^1([0, A]; \Re)$, where $P(\psi) = r\int_0^A \psi(-a)\int_a^A \tilde{k}(s)ds\,da$ with $r := \left(\int_0^A a\tilde{k}(a)da\right)^{-1}$ and $G: C^0([-A, 0]; \Re) \to C^0([0, A]; \Re)$ is the operator defined by the relation $(Gv)(a) = v(-a)$ for all $a \in [0, A]$ for every $v \in C^0([-A, 0]; \Re)$. Define $g$ by means of (5.16) and notice that (2.4) and the fact that $M = y^*\left(\int_0^A p(a)\exp\left(-D^*a - \int_0^a \mu(s)ds\right)da\right)^{-1}$ imply that the function $g \in C^0([0, A]; \Re)$ satisfies $g(a) \geq 0$ for all $a \in [0, A]$ and $\int_0^A g(a)da = 1$.

Next consider the solution $(\eta(t), \psi_t) \in \Re \times S$ of the closed-loop system (5.5), (5.6) with (5.34) and initial condition $\eta(0) = \eta_0$, $\psi(-a) = \psi_0(-a)$ for all $a \in [0, A]$. Lemma 5.3 guarantees that the solution $(\eta(t), \psi_t) \in \Re \times S$ of the closed-loop system (5.5), (5.6) with (5.34) exists for all $t \geq 0$. Moreover, there exist a constants $L > 0$ and a function $\kappa \in K_\infty$ such that for every $(\eta_0, \psi_0) \in \Re \times S$ with $(G\psi_0) \in PC^1([0, A]; \Re)$ the solution $(\eta(t), \psi_t) \in \Re \times S$ of the closed-loop system (5.5), (5.6) with (5.34) and initial condition $\eta(0) = \eta_0$, $\psi(-a) = \psi_0(-a)$ for all $a \in [0, A]$, satisfies estimate (5.35). The solution of the IDE $\psi(t) = \int_0^A \tilde{k}(a)\psi(t - a)da$ is $C^1$ on $(0, +\infty)$ since it coincides with the solution of the delay differential equation $\dot{\psi}(t) = \tilde{k}(0)\psi(t) - \tilde{k}(A)\psi(t - A) + \int_0^A \tilde{k}'(a)\psi(t - a)da$ with the same initial condition. Therefore by virtue of (5.14), the function defined by (5.17) is continuous and $f \in C^1(D_f; (0, +\infty))$, where $D_f = \{(t, a) \in (0, \tau) \times (0, A) : (a - t) \notin B \cup \{0, A\}, t \neq T[T^{-1}t]\}$ and $B \subseteq (0, A)$ is the finite (possibly empty) set where the derivative of $f_0 \in \tilde{X}$ is not defined or is not continuous. Since $\psi_t \in S$, it follows that $f[t] \in \tilde{X}$ for all $t \geq 0$, where $(f[t])(a) = f(t, a)$ for $a \in [0, A]$. Using (5.15), (5.17) and (5.5) we conclude that (5.18) holds. Moreover, using (5.5), (5.6), (5.34), (5.16), (5.17) and (5.18), we conclude that equations (2.5), (2.14) hold for all $t \geq 0$ and equation (2.1) holds for all $(t, a) \in D_f$. Finally, by virtue of (5.14), (5.15) and (5.17), it follows that $f(0, a) = f_0(a)$ for all $a \in [0, A]$.

Using (5.17) and the fact that $|\ln(1 + x)| \leq \dfrac{|x|}{1 + \min(x, 0)}$ for all $x > -1$, we get the estimate for all $t \geq 0$:

$$\max_{a \in [0, A]}\left(\left|\ln\left(\frac{f(t, a)}{f^*(a)}\right)\right|\right) \leq |\eta(t)| + \frac{\max_{a \in [0, A]}(|\psi_t(-a)|)}{1 + \min\left(0, \min_{a \in [0, A]}(\psi_t(-a))\right)} \quad (5.36)$$

Combining (5.35) and (5.36), we obtain the following estimate for all $t \geq 0$:

$$\max_{a \in [0, A]}\left(\left|\ln\left(\frac{f(t, a)}{f^*(a)}\right)\right|\right) \leq \exp(-Lt)\kappa\left(|\eta_0| + \frac{\max_{-A \leq s \leq 0}(|\psi_0(s)|)}{1 + \min\left(0, \min_{-A \leq s \leq 0}(\psi_0(s))\right)}\right) \quad (5.37)$$

Estimate (2.15) for certain $\rho \in K_\infty$ is a direct consequence and inequalities (5.27) for certain $b_1, b_2 \in K_\infty$. The proof is complete. ◁

## 6. Using a Reduced Order Observer

Instead of using the full-order observer (2.6) of the system $\dot{\eta}(t) = D^*(t) - D(t)$, $\dot{D}^*(t) = 0$, one can think of the possibility of using a reduced order observer that estimates the equilibrium value of the dilution rate $D^*$. Such a dynamic, output feedback law will be given by the equations:

$$\dot{z}(t) = -l_1 l_2^{-1} z(t) + l_1^2 l_2^{-1} \ln\left(\frac{y(t)}{y^*}\right) - l_1 D(t) \quad , \quad z(t) \in \Re \tag{6.1}$$

and

$$D(t) = sat\left(-l_2^{-1} z(t) + (\gamma + l_1 l_2^{-1}) \ln\left(\frac{y(t)}{y^*}\right)\right) \tag{6.2}$$

where $l_1, l_2, \gamma > 0$ are constants. In such a case, a solution of the initial-value problem (2.1), (2.2), (2.5) with (6.1), (6.2) with initial condition $(f_0, z_0) \in \tilde{X} \times \Re$, where $\tilde{X} = \left\{ f \in PC^1([0, A]; (0, +\infty)) : f(0) = \int_0^A k(a) f(a) da \right\}$, means a pair of mappings $f \in C^0([0, \tau) \times [0, A]; (0, +\infty))$, $z \in C^1([0, \tau); \Re)$, where $\tau > 0$, which satisfies the following properties:

(i) $f \in C^1(D_f; (0, +\infty))$, where $D_f = \{(t, a) \in (0, \tau) \times (0, A) : (a - t) \notin B \cup \{0, A\}\}$ and $B \subseteq (0, A)$ is the finite (possibly empty) set where the derivative of $f_0 \in \tilde{X}$ is not defined or is not continuous,

(ii) $f[t] \in \tilde{X}$ for all $t \in [0, \tau)$,

(iii) equations (2.5), (6.1), (6.2) hold for all $t \in [0, \tau)$,

(iv) equation $\frac{\partial f}{\partial t}(t, a) + \frac{\partial f}{\partial a}(t, a) = -(\mu(a) + D(t)) f(t, a)$ holds for all $(t, a) \in D_f$, and

(v) $z(0) = z_0$, $f(0, a) = f_0(a)$ for all $a \in [0, A]$.

The mapping $[0, \tau) \ni t \to (f[t], z(t)) \in \tilde{X} \times \Re$ is called the *solution of the closed-loop system (2.1), (2.2), (2.5) with (6.1), (6.2)* and initial condition $(f_0, z_0) \in \tilde{X} \times \Re$ defined for $t \in [0, \tau)$.

For the reduced-order observer case, we are in a position to prove, exactly in the same way of proving Theorem 2.1, the following result. Since its proof is almost identical to the proof of Theorem 2.1, it is omitted.

**Theorem 6.1 (Stabilization with a reduced order observer):** *Consider the age-structured chemostat model (2.1), (2.2) with $k \in PC^1([0, A]; \Re_+)$ under Assumption (A). Then for every $f_0 \in \tilde{X}$ and $z_0 \in \Re$ there exists a unique solution of the closed-loop (2.1), (2.2), (2.5) with (6.1), (6.2) and initial condition $(f_0, z_0) \in \tilde{X} \times \Re$. Furthermore, there exist a constant $L > 0$ and a function $\rho \in K_\infty$ such that for every $f_0 \in \tilde{X}$ and $z_0 \in \Re$ the unique solution of the closed-loop (2.1), (2.2), (2.5) with (6.1), (6.2) and initial condition $(f_0, z_0) \in \tilde{X} \times \Re^2$ is defined for all $t \geq 0$ and satisfies the following estimate:*

$$\max_{a \in [0, A]} \left( \left| \ln\left(\frac{f(t, a)}{f^*(a)}\right) \right| \right) + |z(t) + l_2 D^*| \leq \exp\left(-\frac{L}{4} t\right) \rho\left( \max_{a \in [0, A]} \left( \left| \ln\left(\frac{f_0(a)}{f^*(a)}\right) \right| \right) + |z_0 + l_2 D^*| \right) \quad , \text{ for all } t \geq 0 \tag{6.3}$$

*Moreover, the continuous functional $W : \Re \times C^0([0, A]; (0, +\infty)) \to \Re_+$ defined by:*

$$W(z, f) := \left(\ln(\Pi(f))\right)^2 + G\sqrt{Q(z,f)} + \beta Q(z,f) \tag{6.4}$$

where $\beta \geq 0$ is an arbitrary constant,

$$Q(z,f) := \left(z - l_1 \ln(\Pi(f)) + l_2 D^*\right)^2 + \frac{M}{2} \left( \frac{\max\limits_{a \in [0,A]}\left(\exp(-\sigma a)\left|\frac{f(a) - \Pi(f)f^*(a)}{f^*(a)}\right|\right)}{\min\left(\Pi(f), \min\limits_{a \in [0,A]}\left(\frac{f(a)}{f^*(a)}\right)\right)} \right)^2 \tag{6.5}$$

$\Pi : C^0([0,A]; \Re) \to \Re$ is given by (2.8), $\sigma > 0$ is a sufficiently small constant and $M, G > 0$ are sufficiently large constants, is a Lyapunov functional for the closed-loop system (2.1), (2.2), (2.5) with (6.1), (6.2), in the sense that every solution $(f[t], z(t)) \in \tilde{X} \times \Re$ of the closed-loop system (2.1), (2.2), (2.5) with (6.1), (6.2) satisfies the inequality:

$$\limsup_{h \to 0^+}\left(h^{-1}\left(W(z(t+h), f[t+h]) - W(z(t), f[t])\right)\right) \leq -L \frac{W(z(t), f[t])}{1 + \sqrt{W(z(t), f[t])}}, \text{ for all } t \geq 0 \tag{6.6}$$

The family of dynamic, bounded, output feedback laws (6.1), (6.2) presents the same features as the family (2.6), (2.7). The only difference lies in the dimension of the observer.

## 7. Simulations

To demonstrate the sampled-data control design from Theorem 2.4, three simulations were carried out. In each simulation we considered the case where

$$A = 2, \; \mu(a) \equiv \mu = 0.1, \; D^* = 1 \tag{7.1}$$

and the birth modulus is given by

$$k(a) = ag, \; a \in [0,1] \text{ and } k(a) = (2-a)g, \; a \in (1,2], \tag{7.2}$$

where $g > 0$ is the constant for which the Lotka-Sharpe condition (2.3) holds. The model is dimensionless (a dimensionless version of (2.1), (2.2) can be obtained by using appropriate scaling of all variables). After a simple calculation it can be found that the constant $g > 0$ is given by:

$$g = \frac{(\mu + D^*)^2}{\left(1 - e^{-(\mu + D^*)}\right)^2} = 2.718728 \tag{7.3}$$

The output is given by the equation

$$y(t) = \int_0^2 f(t, a)da, \; t \geq 0 \tag{7.4}$$

In other words, the output is the total concentration of the microorganism in the chemostat. The chosen equilibrium profile that has to be stabilized is the profile that is given by the equation:

$$f^*(a) = \exp\left(-(D^* + \mu)a\right), \text{ for } a \in [0,2] \tag{7.5}$$

The equilibrium value of the output is given by:

$$y^* = \frac{1 - \exp\left(-2(D^* + \mu)\right)}{D^* + \mu} = 0.808361 \tag{7.6}$$

Two sampled-data feedback laws were tested: the state feedback law given by

$$D(t) = D_i = sat\left(D^* + T^{-1} \ln\left(f(iT,0)/f^*(0)\right)\right),$$
for all $t \in [iT,(i+1)T)$ and for all integers $i \geq 0$ \hfill (7.7)

which is the sampled-data feedback law proposed in [18], and the output feedback law given by

$$D(t) = D_i = sat\left(D^* + T^{-1} \ln\left(y(iT)/y^*\right)\right),$$
for all $t \in [iT,(i+1)T)$ and for all integers $i \geq 0$ \hfill (7.8)

which is the sampled-data feedback law given by Theorem 2.4. For both feedback laws we chose:

$$T = 0.4, \quad D_{\min} = 0.5, \quad D_{\max} = 1.5 \hfill (7.9)$$

The following family of functions was used for initial conditions:

$$f_0(a) = b_0 - b_1 a + c\exp(-\theta a), \text{ for } a \in [0,2] \hfill (7.10)$$

where $b_0, c, \theta > 0$ are free parameters and the constant $b_1 \in \mathfrak{R}$ is chosen so that the condition $f_0(0) = \int_0^2 k(a) f_0(a) da$ holds. After some simple calculations we find that

$$b_1 = g^{-1}(g-1)b_0 + c\theta^{-2}\left(1 - e^{-\theta}\right)^2 - cg^{-1} \hfill (7.11)$$

However, we notice that not all parameters $b_0, c, \theta > 0$ can be used because the additional condition $\min_{a \in [0,2]}(f_0(a)) > 0$ must hold as well.

The simulations were made with the generation of a uniform grid of function values:

$$f(ih, jh), \text{ for } j = 0,1,...,50 \text{ and } i \geq 0 \hfill (7.12)$$

where $h = 0.04$. For $i = 0$ we had $f(0, jh) = f_0(jh)$ for $j = 0,1,...,50$. The calculation of the integrals $y(ih) = \int_0^2 f(ih,a) da$ and $f(ih,0) = \int_0^2 k(a) f(ih,a) da$ for every $i \geq 0$, was made numerically. However, since we wanted the numerical integrator to be able to evaluate <u>exactly</u> the integrals $y(ih) = \int_0^2 f(ih,a) da$ and $f(ih,0) = \int_0^2 k(a) f(ih,a) da$ for every $i \geq 0$ when $f[ih]$ is an exponential function, i.e., when $f(ih,a) = C\exp(\sigma a)$ for $a \in [0,2]$ and for certain constants $C, \sigma \in \mathfrak{R}$, we could not use a conventional numerical integration scheme like the trapezoid's rule or Simpson's rule. The reason for this demand to be able to evaluate <u>exactly</u> the integrals $y(ih) = \int_0^2 f(ih,a) da$ and $f(ih,0) = \int_0^2 k(a) f(ih,a) da$ for every $i \geq 0$ when $f[ih]$ is an exponential function is explained by the fact that the equilibrium profile (given by (7.5)) is an exponential function and we would like to avoid a steady-state error due to the error induced by the numerical integration. To this end, we used the following integration schemes:

$$\int_{jh}^{(j+1)h} f(ih,a)da \approx I_i(j) = h\frac{f(ih,(j+1)h) - f(ih,jh)}{\ln(f(ih,(j+1)h)) - \ln(f(ih,jh))} \quad for \quad f(ih,(j+1)h) \neq f(ih,jh)$$

$$\int_{jh}^{(j+1)h} f(ih,a)da \approx I_i(j) = hf(ih,jh) \quad for \quad f(ih,(j+1)h) = f(ih,jh)$$

$$for \quad j = 2,3,\ldots,49 \text{ and } i \geq 0 \tag{7.13}$$

$$\int_0^{2h} f(ih,a)da \approx I_i(2) = h\frac{f(ih,2h) - f^2(ih,h)/f(ih,2h)}{\ln(f(ih,2h)) - \ln(f(ih,h))} \quad for \quad f(ih,h) \neq f(ih,2h)$$

$$\int_0^{2h} f(ih,a)da \approx I_i(2) = 2hf(ih,h) \quad for \quad f(ih,h) = f(ih,2h) \tag{7.14}$$

$$\int_{jh}^{(j+1)h} af(ih,a)da \approx J_i(j) = h^2 \frac{f(ih,(j+1)h) + j(f(ih,(j+1)h) - f(ih,jh))}{\ln(f(ih,(j+1)h)) - \ln(f(ih,jh))}$$

$$-\frac{h^2(f(ih,(j+1)h) - f(ih,jh))}{(\ln(f(ih,(j+1)h)) - \ln(f(ih,jh)))^2} \quad for \quad f(ih,(j+1)h) \neq f(ih,jh)$$

$$\int_{jh}^{(j+1)h} af(ih,a)da \approx J_i(j) = \frac{2j+1}{2}h^2 f(ih,jh) \quad for \quad f(ih,(j+1)h) = f(ih,jh)$$

$$for \quad j = 2,3,\ldots,24 \text{ and } i \geq 0 \tag{7.15}$$

$$\int_0^{2h} af(ih,a)da \approx J_i(j) = \frac{2h^2 f(ih,2h)}{\ln(f(ih,2h)) - \ln(f(ih,h))}$$

$$-\frac{h^2\left(f(ih,2h) - \frac{f^2(ih,h)}{f(ih,2h)}\right)}{(\ln(f(ih,2h)) - \ln(f(ih,h)))^2} \quad for \quad f(ih,h) \neq f(ih,2h) \tag{7.16}$$

$$\int_0^{2h} af(ih,a)da \approx J_i(j) = 2h^2 f(ih,2h) \quad for \quad f(ih,h) \neq f(ih,2h)$$

$$\int_{jh}^{(j+1)h} (2-a)f(ih,a)da \approx K_i(j) = -h^2 \frac{f(ih,(j+1)h)}{\ln(f(ih,(j+1)h)/f(ih,jh))}$$

$$+\left(2 - jh + \frac{h}{\ln(f(ih,(j+1)h)/f(ih,jh))}\right)\frac{hf(ih,(j+1)h) - hf(ih,jh)}{\ln(f(ih,(j+1)h)/f(ih,jh))} \quad for \quad f(ih,(j+1)h) \neq f(ih,jh)$$

$$\int_{jh}^{(j+1)h} (2-a)f(ih,a)da \approx K_i(j) = \left(2 - \frac{2j+1}{2}h\right)hf(ih,jh) \quad for \quad f(ih,(j+1)h) = f(ih,jh)$$

$$for \quad j = 25,26,\ldots,49 \text{ and } i \geq 0 \tag{7.17}$$

The derivation of formulas (7.13), (7.14), (7.15), (7.16), (7.17) is based on the interpolation of the function $\tilde{f}_j(a) = C_j \exp(\sigma_j a)$ through the points $(jh, f(ih,jh))$ and $((j+1)h, f(ih,(j+1)h))$ for $j = 1,2,\ldots,49$. More specifically, we obtain for $j = 1,2,\ldots,49$:

$$\sigma_j = h^{-1}\ln(f(ih,(j+1)h)/f(ih,jh))$$
$$C_j = f(ih,jh)(f(ih,(j+1)h)/f(ih,jh))^{-j} \tag{7.18}$$

Based on the above interpolation, the exact integration formulas are used. For example, for the integral $\int_{jh}^{(j+1)h} af(ih,a)da$ for $j = 2,3,...,24$, we get when $\sigma_j = h^{-1}\ln(f(ih,(j+1)h)/f(ih,jh)) \neq 0$

$$\int_{jh}^{(j+1)h} af(ih,a)da \approx \int_{jh}^{(j+1)h} a\tilde{f}_j(a)da = C_j \int_{jh}^{(j+1)h} a\exp(\sigma_j a)da$$

$$= C_j \sigma_j^{-1} h \exp(\sigma_j jh)((j+1)\exp(\sigma_j h) - j) - C_j \sigma_j^{-2}(\exp(\sigma_j(j+1)h) - \exp(\sigma_j jh))$$

and when $\sigma_j = h^{-1}\ln(f(ih,(j+1)h)/f(ih,jh)) = 0$

$$\int_{jh}^{(j+1)h} af(ih,a)da \approx \int_{jh}^{(j+1)h} a\tilde{f}_j(a)da = C_j \int_{jh}^{(j+1)h} ada = \frac{h^2}{2} C_j (2j+1)$$

Combining the above formulas with the estimated values for $C_j, \sigma_j$ given by (7.18), we obtain formula (7.15). Similarly, we derive formulas (7.13), (7.14), (7.16) and (7.17). Notice that the formulas (7.13), (7.14), (7.15), (7.16), (7.17) allow the numerical evaluation of the integrals $y(ih) = \int_0^2 f(ih,a)da$ and $f(ih,0) = \int_0^2 k(a)f(ih,a)da$ for every $i \geq 0$ without knowledge of $f(ih,0)$.

Since the time step has been chosen to be equal to the discretization space step $h$, we are able to use the exact formula:

$$f((i+1)h, jh) = f(ih,(j-1)h)\exp(-(\mu + D_i)h), \text{ for } j = 1,2,...,50 \text{ and } i \geq 0 \quad (7.19)$$

Therefore, we are in a position to use the following algorithm for the simulation of the closed-loop system under the effect of the output feedback law (7.8).

<u>Algorithm:</u> Given $f(ih, jh)$, for $j = 1,...,50$ and certain $i \geq 0$ do the following:

→ Calculate $f(ih,0) \approx g \sum_{j=2}^{24} J_i(j) + g \sum_{j=25}^{49} K_i(j)$, where $J_i(j), K_i(j)$ are given by (7.15), (7.16), (7.17).

→ Calculate $y(ih) \approx \sum_{j=2}^{49} I_i(j)$, where $I_i(j)$ is given by (7.13), (7.14).

→ If $\frac{ih}{T}$ is an integer then set $D_i = \max(D_{\min}, \min(D_{\max}, D^* + T^{-1}\ln(y(ih)/y^*)))$; else set $D_i = D_{i-1}$.

→ Calculate $f((i+1)h, jh)$, for $j = 1,...,50$ using (7.19).

The above algorithm with obvious modifications was also used for the simulation of the open-loop system as well as for the simulation of the closed-loop system under the effect of the output feedback law (7.7).

In our first simulation, we used the parameter values $b_0 = 0.2$, $b_1 = 0.15184212$, $c = 0.8$, $\theta = 1$ in our initial conditions. In Figure 2, we plot the control values and the newborn individual values. We show the values for the open loop feedback $D(t) \equiv 1$, and for the state and output feedbacks from (7.7) and (7.8). Our simulation shows the efficacy of our control design.

In our second simulation, we changed the parameter values to $b_0 = 1$, $b_1 = 0.7592106$, $c = 4$, $\theta = 1$ and plotted the same values as before, in Figure 3. The responses for the output feedback law (7.7) and the output feedback law (7.8) are almost identical. The second simulation was made with an initial condition which is not close to the equilibrium profile (in the sense that it is an initial condition with very large initial population). The difference in the performance of the feedback controllers (7.7) and (7.8) cannot be distinguished.

In the final simulation, we tested the robustness of the controller with respect to errors in the choice of $D^*$ being used in the controllers. We chose the values $b_0 = 0.2$, $b_1 = 0.15184212$, $c = 0.8$, $\theta = 1$, but instead of (7.7) and (7.8), we applied the following controllers: the state feedback law

$$D(t) = D_i = sat\left(0.7 + T^{-1} \ln\left(f(iT,0)/f^*(0)\right)\right),$$
for all $t \in [iT, (i+1)T)$ and for all integers $i \geq 0$  (7.20)

and the output feedback law given by

$$D(t) = D_i = sat\left(0.7 + T^{-1} \ln\left(y(iT)/y^*\right)\right),$$
for all $t \in [iT, (i+1)T)$ and for all integers $i \geq 0$  (7.21)

We obtained in both cases: $\lim_{t \to +\infty} (f(t,0)) = 1.1275$ and $\lim_{t \to +\infty} (D(t)) = D^* = 1$. A -30% error in $D^*$ gives a +12.75% steady-state deviation from the desired value of the newborn individuals. See Figure 4. Notice that a constant error in $D^*$ is equivalent to an error in the set point since we have:

$$D(t) = D_i = sat\left(0.7 + T^{-1} \ln\left(f(iT,0)/f^*(0)\right)\right)$$
$$= sat\left(D^* + T^{-1} \ln\left(f(iT,0)/f^*(0)\right) - T^{-1} 0.12\right)$$
$$= sat\left(D^* + T^{-1} \ln\left(f(iT,0)/(1.1275 f^*(0))\right)\right)$$

for the state feedback case (7.7) and

$$D(t) = D_i = sat\left(0.7 + T^{-1} \ln\left(y(iT)/y^*\right)\right)$$
$$= sat\left(D^* + T^{-1} \ln\left(y(iT)/y^*\right) - T^{-1} 0.12\right)$$
$$= sat\left(D^* + T^{-1} \ln\left(y(iT)/(1.1275 y^*)\right)\right)$$

for the output feedback case (7.8).

## 8. Concluding Remarks

Age-structured chemostats present challenging control problems for first-order hyperbolic PDEs that require novel results. We studied the problem of stabilizing an equilibrium age profile in an age-structured chemostat, using the dilution rate as the control. We built a family of dynamic, bounded, output feedback laws with continuously adjusted input that ensures asymptotic stability under arbitrary physically meaningful initial conditions and does not require knowledge of the model. We also built a sampled-data, bounded, output feedback stabilizer which guarantees asymptotic stability under arbitrary physically meaningful initial conditions and requires only the knowledge of one parameter: the equilibrium value of the dilution rate. In addition, we provided a family of CLFs for the age-structured chemostat model. The construction of the CLF was based on novel stability results on linear IDEs, which are of independent interest. The newly developed results, provide a Lyapunov-like proof of the scalar, strong ergodic theorem for special cases of the integral kernel.

Since the growth of the microorganism may sometimes depend on the concentration of a limiting substrate, it would be useful to solve the stabilization problem for an enlarged system that has one PDE for the age distribution, coupled with one ODE for the substrate (as proposed in [29], in the context of studying limit cycles with constant dilution rates instead of a control). This is going to be the topic of our future research.

**Acknowledgements:** The authors would like to thank Michael Malisoff for his help in the initial stages of the writing process of the paper.

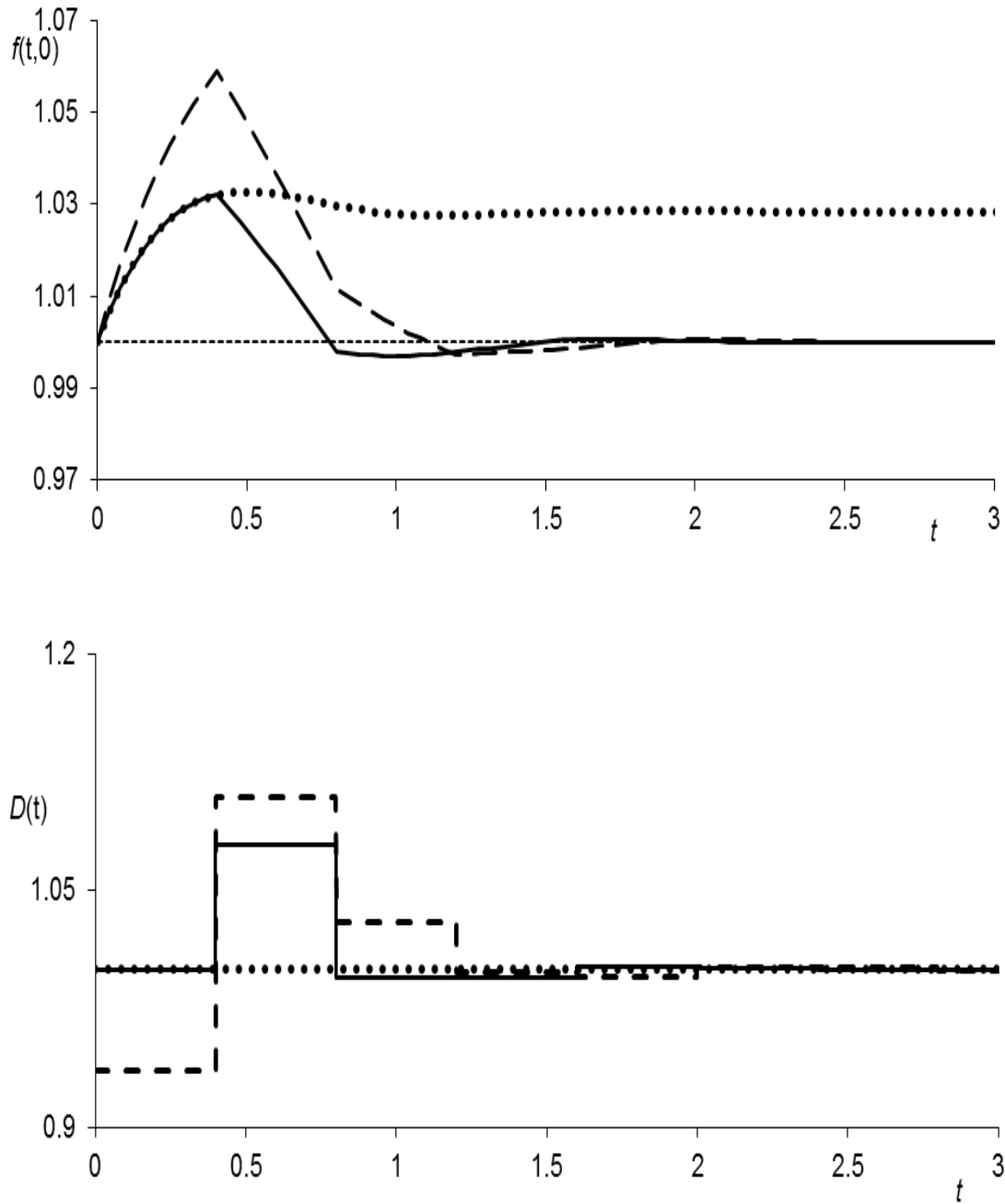

**Figure 2:** Simulation for the initial condition given by (7.10) with $b_0 = 0.2$, $b_1 = 0.15184212$, $c = 0.8$, $\theta = 1$. The upper part of the figure shows the response for the newborn individuals. The solid line with bullets is for the state feedback (7.7), the dashed line is for the output feedback (7.8) and the bulleted line is for the open-loop system with $D(t) \equiv 1$. The lower part of the figure shows the applied control action $D(t)$. Again, the solid line is for the state feedback (7.7), the dashed line is for the output feedback (7.8), while the bulleted line shows the equilibrium value ($D^* = 1$) of the dilution rate.

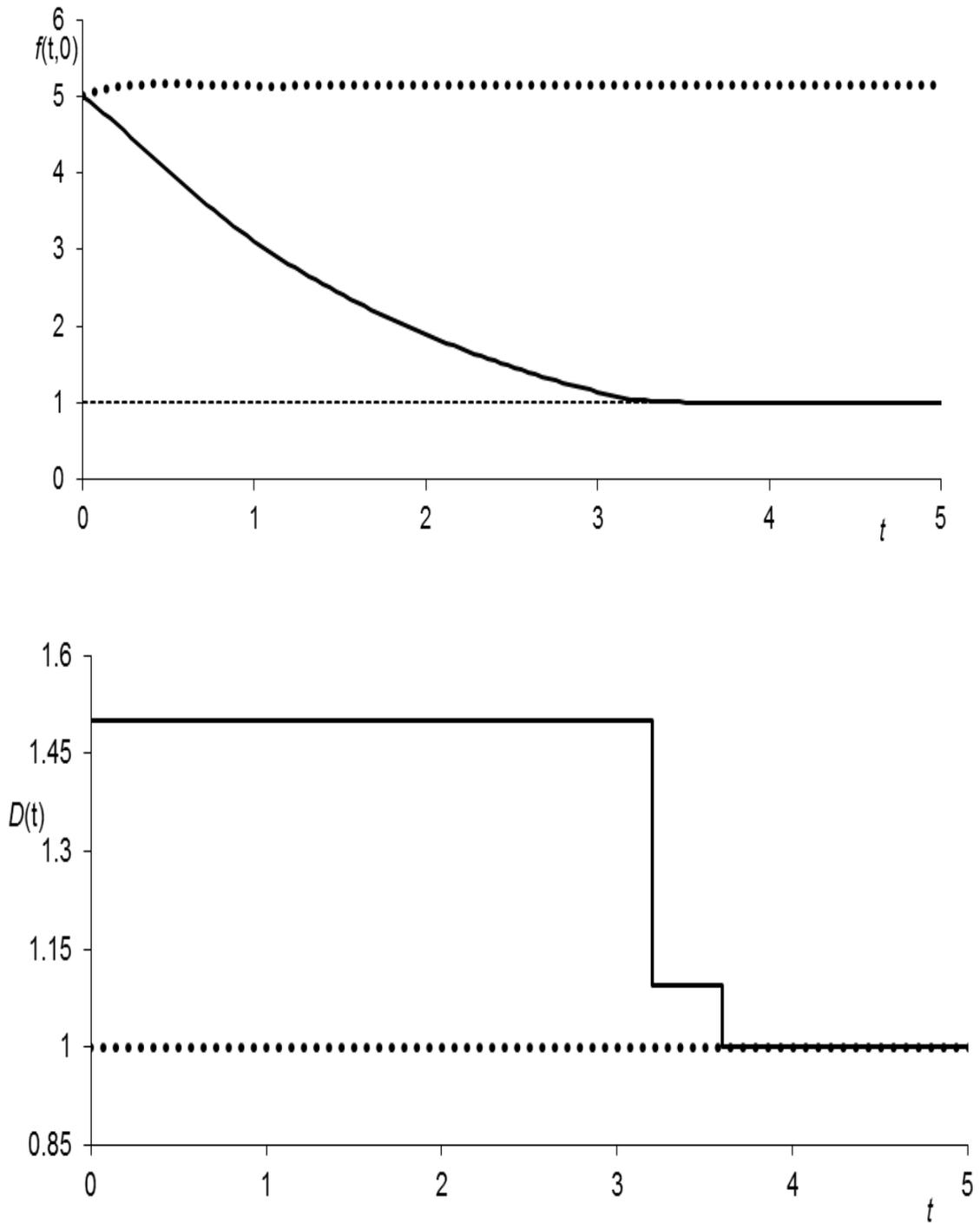

**Figure 3:** Simulation for the initial condition given by (7.10) with $b_0 = 1$, $b_1 = 0.7592106$, $c = 4$, $\theta = 1$. The upper part of the figure shows the response for the newborn individuals. The solid line is both for the state feedback (7.7) and for the output feedback (7.8) (identical) and the bulleted line is for the open-loop system with $D(t) \equiv 1$. The lower part of the figure shows the applied control action $D(t)$. Again, the solid line is both for the state feedback (7.7) and for the output feedback (7.8), while the bulleted line shows the equilibrium value ($D^* = 1$) of the dilution rate.

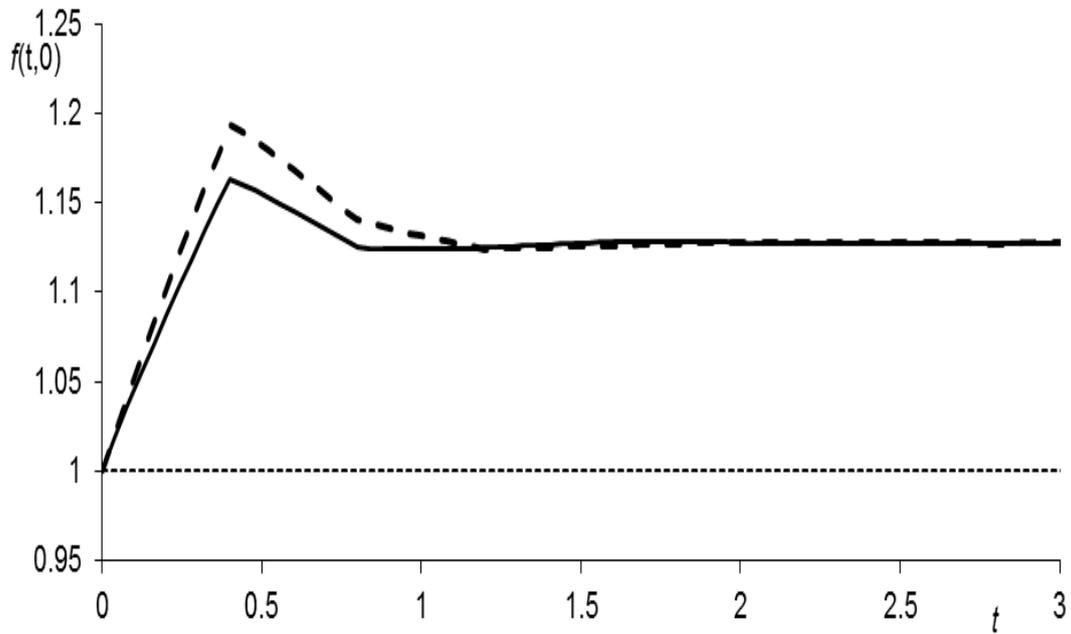

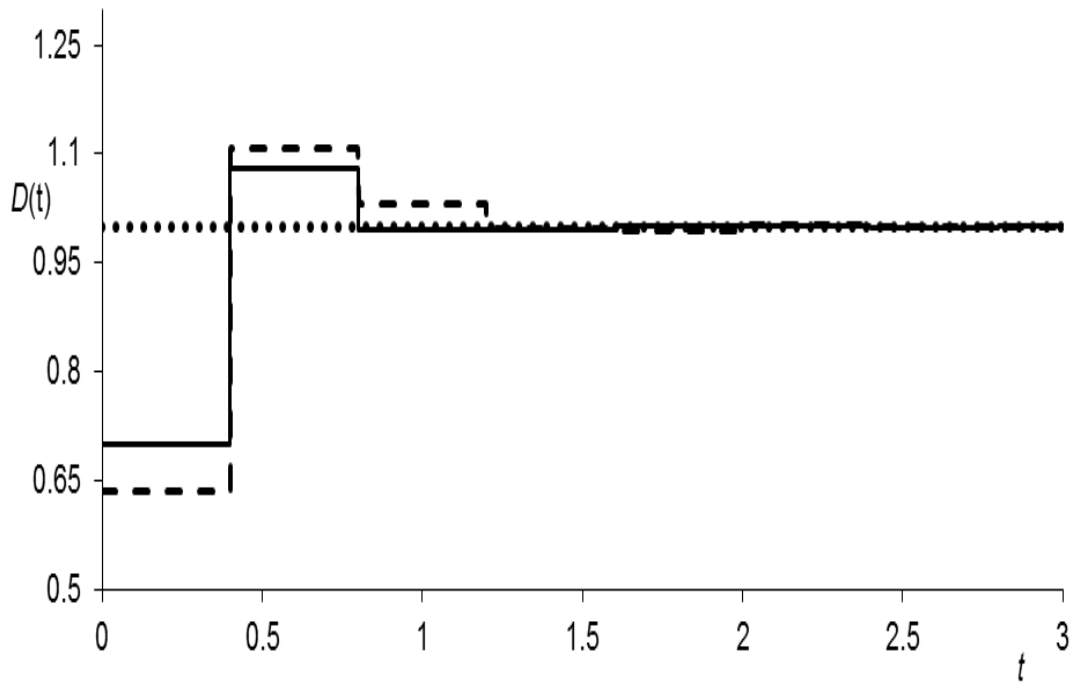

**Figure 4:** Control in presence of modeling errors (-30% error in $D^*$). Simulation for the initial condition given by (7.10) with $b_0 = 0.2$, $b_1 = 0.15184212$, $c = 0.8$, $\theta = 1$. The upper part of the figure shows the response for the newborn individuals. The solid line is for the state feedback (7.20), the dashed line is for the output feedback (7.21), while the dotted line shows the equilibrium value ($f^* = 1$) of the newborn individuals. The lower part of the figure shows the applied control action $D(t)$. Again, the solid line is for the state feedback (7.20), the dashed line is for the output feedback (7.21), while the bulleted line shows the equilibrium value ($D^* = 1$) of the dilution rate.

**References**
# References

[1] Bastin, G., and J.-M. Coron, "On Boundary Feedback Stabilization of Non-Uniform Linear 2x2 Hyperbolic Systems Over a Bounded Interval", *Systems and Control Letters*, 60, 2011, 900-906.

[2] Bernard, P. and M. Krstic, "Adaptive Output-Feedback Stabilization of Non-Local Hyperbolic PDEs", *Automatica*, 50, 2014, 2692-2699.

[3] Boucekkine, R., N. Hritonenko and Y. Yatsenko, *Optimal Control of Age-structured Populations in Economy, Demography, and the Environment* (Google eBook), 2013.

[4] Brauer, F. and C. Castillo-Chavez, *Mathematical Models in Population Biology and Epidemiology*, Springer-Verlag, New York, 2001.

[5] Charlesworth, B., *Evolution in Age-Structured Populations*, 2nd Edition, Cambridge University Press, 1994.

[6] Coron, J.-M., R. Vazquez, M. Krstic, and G. Bastin, "Local Exponential H2 Stabilization of a 2x2 Quasilinear Hyperbolic System Using Backstepping", *SIAM Journal of Control and Optimization*, 51, 2013, 2005-2035.

[7] Di Meglio, F., R. Vazquez, and M. Krstic, "Stabilization of a System of $n + 1$ Coupled First-Order Hyperbolic Linear PDEs with a Single Boundary Input," *IEEE Transactions on Automatic Control*, 58, 2013, 3097-3111.

[8] Feichtinger, G., G. Tragler and V. M. Veliov, "Optimality Conditions for Age-Structured Control Systems", *Journal of Mathematical Analysis and Applications*, 288(1), 2003, 47–68.

[9] Gouze, J. L. and G. Robledo, "Robust Control for an Uncertain Chemostat Model", *International Journal of Robust and Nonlinear Control*, 16(A.19), 2006, 133-155.

[10] Hale, J. K. and S. M. V. Lunel, *Introduction to Functional Differential Equations*, Springer-Verlag, New York, 1993.

[11] Inaba, H., "A Semigroup Approach to the Strong Ergodic Theorem of the Multistate Stable Population Process", *Mathematical Population Studies*, 1(4.1), 1988, 49-77.

[12] Inaba, H., "Asymptotic Properties of the Inhomogeneous Lotka-Von Foerster System", *Mathematical Population Studies*, 1(A.19), 1988, 247-264.

[13] Karafyllis, I., C. Kravaris, L. Syrou and G. Lyberatos, "A Vector Lyapunov Function Characterization of Input-to-State Stability with Application to Robust Global Stabilization of the Chemostat", *European Journal of Control*, 14(4.1), 2008, 47-61.

[14] Karafyllis, I., C. Kravaris and N. Kalogerakis, "Relaxed Lyapunov Criteria for Robust Global Stabilization of Nonlinear Systems", *International Journal of Control*, 82(7.11), 2009, 2077-2094.

[15] Karafyllis, I., and Z.-P. Jiang, "A New Small-Gain Theorem with an Application to the Stabilization of the Chemostat", *International Journal of Robust and Nonlinear Control*, 22(7.14), 2012, 1602–1630.

[16] Karafyllis, I. and M. Krstic, "On the Relation of Delay Equations to First-Order Hyperbolic Partial Differential Equations", *ESAIM Control, Optimisation and Calculus of Variations*, 20(A.19), 2014, 894 - 923.

[17] Karafyllis, I., M. Malisoff and M. Krstic, "Ergodic Theorem for Stabilization of a Hyperbolic PDE Inspired by Age-Structured Chemostat", arXiv:1501.04321 [math.OC].

[18] Karafyllis, I., M. Malisoff and M. Krstic, "Sampled-Data Feedback Stabilization of Age-Structured Chemostat Models", *Proceedings of the American Control Conference 2015*, Chicago, IL, U.S.A., pp. 4549-4554.

[19] Khalil, H. K., *Nonlinear Systems*, 2nd Edition, Prentice-Hall, 1996.

[20] Krstic, M., and A. Smyshlyaev, "Backstepping Boundary Control for First-Order Hyperbolic PDEs and Application to Systems With Actuator and Sensor Delays", *Systems and Control Letters*, 57(A.24), 2008, 750-758.

[21] Mazenc, F., M. Malisoff, J. Harmand, "Stabilization and Robustness Analysis for a Chemostat Model with Two Species and Monod Growth Rates via a Lyapunov Approach", *Proceedings of the 46th IEEE Conference on Decision and Control*, New Orleans, 2007.

[22] Mazenc, F., M. Malisoff and J. Harmand, "Further Results on Stabilization of Periodic Trajectories for a Chemostat with Two Species", *IEEE Transactions on Automatic Control*, 53(4.1), 2008, 66-74.
[23] Melchor-Aguilar, D., "Exponential Stability of Some Linear Continuous Time Difference Systems", *Systems and Control Letters*, 61, 2012, 62–68.
[24] Pazy, A., *Semigroups of Linear Operators and Applications to Partial Differential Equations*, Springer-Verlag, New York, 1983.
[25] Rao N. S., and E. O. Roxin, "Controlled Growth of Competing Species", *SIAM Journal on Applied Mathematics*, 50(3), 1990, 853-864.
[26] Rundnicki, R. and M. C. Mackey, "Asymptotic Similarity and Malthusian Growth in Autonomous and Nonautonomous Populations", *Journal of Mathematical Analysis and Applications*, 187, 1994, 548-566.
[27] Smith, H. and P. Waltman, *The Theory of the Chemostat. Dynamics of Microbial Competition*, Cambridge Studies in Mathematical Biology, 13, Cambridge University Press: Cambridge, 1995.
[28] Sun, B., "Optimal Control of Age-Structured Population Dynamics for Spread of Universally Fatal Diseases II", *Applicable Analysis: An International Journal*, 93(A.23), 2014, 1730-1744.
[29] Toth, D., and M. Kot, "Limit Cycles in a Chemostat Model for a Single Species with Age Structure", *Mathematical Biosciences*, 202, 2006 194–217.

# Appendix

**Proof of Proposition 2.3:** Define for $\lambda \in [0, r^{-1}\varepsilon]$:

$$g(\lambda) := \int_0^A \left| \tilde{k}(a) - r\lambda \int_a^A \tilde{k}(s)ds \right| da \tag{A.1}$$

Since $T := \sup\{a \in [0, A] : \tilde{k}(a) > 0\}$, we have:

$$g(\lambda) = \int_0^T \left| \tilde{k}(a) - r\lambda \int_a^A \tilde{k}(s)ds \right| da \text{ and}$$

$$\int_0^A a\tilde{k}(a)da = \int_0^T a\tilde{k}(a)da = \int_0^A \int_a^A \tilde{k}(s)ds\,da = \int_0^T \int_a^A \tilde{k}(s)ds\,da = \int_0^T \int_a^T \tilde{k}(s)ds\,da \tag{A.2}$$

Define the Lebesgue measurable sets:

$$S^+ = \left\{ a \in [0, T] : \tilde{k}(a) > r\lambda \int_a^A \tilde{k}(s)ds \right\} \quad S^- = \left\{ a \in [0, T] : \tilde{k}(a) \leq r\lambda \int_a^A \tilde{k}(s)ds \right\} \tag{A.3}$$

and the integrals

$$I^+ = \int_{S^+} \tilde{k}(a)da \quad J^+ = \int_{S^+} \left( \int_a^A \tilde{k}(s)ds \right) da \tag{A.4}$$

Notice that $S^+ \cup S^- = [0, T]$ and consequently, equations (A.2) in conjunction with the fact that $r = \left( \int_0^A a\tilde{k}(a)da \right)^{-1}$, implies that:

$$g(\lambda) = \int_{S^+} \left( \tilde{k}(a) - r\lambda \int_a^A \tilde{k}(s)ds \right) da + \int_{S^-} \left( r\lambda \int_a^A \tilde{k}(s)ds - \tilde{k}(a) \right) da$$

$$= I^+ - r\lambda J^+ + \left( \int_0^T \left( r\lambda \int_a^A \tilde{k}(s)ds - \tilde{k}(a) \right) da - \int_{S^+} \left( r\lambda \int_a^A \tilde{k}(s)ds - \tilde{k}(a) \right) da \right) \tag{A.5}$$

$$= I^+ - r\lambda J^+ + r\lambda(r^{-1} - J^+) - (1 - I^+) = 2I^+ - 1 + \lambda(1 - 2rJ^+)$$

Since $S^+ \cup S^- = [0,T]$ and $\int_a^A \tilde{k}(s)ds \leq 1$ for all $a \in [0,A]$, we get:

$$r^{-1} = \int_0^T \left( \int_a^A \tilde{k}(s)ds \right) da = J^+ + \int_{S^-} \left( \int_a^A \tilde{k}(s)ds \right) da \leq |S^-| + J^+ \tag{A.6}$$

Moreover, since $\lambda \in [0, r^{-1}\varepsilon]$ and $\int_a^A \tilde{k}(s)ds \leq 1$ for all $a \in [0,A]$, it follows that $S^- \subseteq S_\varepsilon$. Therefore, we obtain from (A.6):

$$2 \leq 2r|S_\varepsilon| + 2rJ^+$$

or equivalently

$$1 - 2rJ^+ \leq 2r|S_\varepsilon| - 1 \tag{A.7}$$

Definition (A.4) and the fact that $\tilde{k}(a) \geq 0$ for all $a \in [0,A]$ with $\int_0^A \tilde{k}(a)da = 1$ implies that $I^+ \leq 1$. Combining the previous inequality with (A.7) and (A.5) we obtain the desired inequality $\int_0^A \left| \tilde{k}(a) - r\lambda \int_a^A \tilde{k}(s)ds \right| da \leq 1 - \lambda(1 - 2r|S_\varepsilon|)$ for all $\lambda \in [0, r^{-1}\varepsilon]$. The proof is complete. ◁

**Proof of Lemma 4.1:** Local existence and uniqueness for every initial condition $\xi \in L^\infty([-A,0]; \Re)$ is guaranteed by Theorem 2.1 in [16].

Define for all $t \geq 0$ for which the solution of (4.1) exists:

$$V(t) = \sup_{-A \leq a < 0} (x(t+a)), \quad W(t) = \inf_{-A \leq a < 0} (x(t+a)) \tag{A.8}$$

Let $q > 0$ and $t \geq 0$ be sufficiently small so that the solution exists on $[t, t+q)$. We get from definition (A.8) and equation (4.1):

$$\begin{aligned}
V(t+q) &= \sup_{-A \leq a < 0} (x(t+q+a)) \\
&= \sup_{q-A \leq s < q} (x(t+s)) = \max\left( \sup_{q-A \leq s < 0} (x(t+s)), \sup_{0 \leq s < q} (x(t+s)) \right) \\
&\leq \max\left( V(t), \sup_{0 \leq s < q} \left( \int_0^A \varphi(a) x(t+s-a) da \right) \right) \\
&= \max\left( V(t), \sup_{0 \leq s < q} \left( \int_\delta^A \varphi(a) x(t+s-a) da + \int_0^\delta \varphi(a) x(t+s-a) da \right) \right) \\
&\leq \max\left( V(t), \sup_{0 \leq s < q} \left( \sup_{s-A \leq l < s-\delta} (x(t+l)) \int_\delta^A \varphi(a) da + \sup_{s-\delta \leq l < s} (x(t+l)) \int_0^\delta \varphi(a) da \right) \right) \\
&\leq \max\left( V(t), \sup_{-A \leq l < q-\delta} (x(t+l)) \int_\delta^A \varphi(a) da + \sup_{-\delta \leq l < q} (x(t+l)) \int_0^\delta \varphi(a) da \right)
\end{aligned} \tag{A.9}$$

Using the fact that $L := \int_0^A \varphi(a)da \geq 1$ and assuming that $q \leq \min(\delta, A-\delta)$, we obtain from (A.9):

$$V(t+q) \leq \max(V(t), V(t)(L-c) + cV(t+q)) \tag{A.10}$$

Using the fact that $\delta > 0$ is a constant with $c := \int_0^\delta \varphi(a) da < 1$ and the fact that $L := \int_0^A \varphi(a) da \geq 1$, we distinguish the following cases: (i) $V(t) \leq V(t)(L-c) + cV(t+q)$ and in this case (A.10) implies that $V(t+q) \leq \left(\frac{L-c}{1-c}\right) V(t)$, (ii) $V(t) > V(t)(L-c) + cV(t+q)$ and in this case (A.10) implies that $V(t+q) \leq V(t)$. Therefore, in any case we get:

$$V(t+q) \leq \max\left(\frac{L-c}{1-c} V(t), V(t)\right) \tag{A.11}$$

Similarly, we get from definition (A.8), equation (4.1) and the fact that $q \leq \min(\delta, A-\delta)$:

$$\begin{aligned}
W(t+q) &= \inf_{-A \leq a < 0} (x(t+q+a)) \\
&= \inf_{q-A \leq s < q} (x(t+s)) = \min\left(\inf_{q-A \leq s < 0} (x(t+s)), \inf_{0 \leq s < q} (x(t+s))\right) \\
&\geq \min\left(W(t), \inf_{0 \leq s < q}\left(\int_0^A \varphi(a) x(t+s-a) da\right)\right) \\
&= \min\left(W(t), \inf_{0 \leq s < q}\left(\int_\delta^A \varphi(a) x(t+s-a) da + \int_0^\delta \varphi(a) x(t+s-a) da\right)\right) \\
&\geq \min\left(W(t), \inf_{0 \leq s < q}\left(\inf_{s-A \leq l < s-\delta}(x(t+l))\int_\delta^A \varphi(a) da + \inf_{s-\delta \leq l < s}(x(t+l))\int_0^\delta \varphi(a) da\right)\right) \\
&\geq \min\left(W(t), (L-c) \inf_{-A \leq l < q-\delta}(x(t+l)) + c \inf_{-\delta \leq l < q}(x(t+l))\right) \\
&\geq \min(W(t), W(t)(L-c) + cW(t+q))
\end{aligned} \tag{A.12}$$

Using the fact that $\delta > 0$ is a constant with $c := \int_0^\delta \varphi(a) da < 1$ and the fact that $L := \int_0^A \varphi(a) da \geq 1$, we obtain (again by distinguishing cases) from (A.12):

$$W(t+q) \geq \min\left(W(t), \left(\frac{L-c}{1-c}\right) W(t)\right) \tag{A.13}$$

It follows from (A.11) and (A.13) that the solution of (4.1) is bounded on $[t, t+q)$ for $q \leq \min(\delta, A-\delta)$. A standard contradiction argument in conjunction with Theorem 2.1 in [16] implies that the solution exists for all $t \geq 0$. Indeed, a finite maximal existence time $t_{\max} < +\infty$ for the solution in conjunction with Theorem 2.1 in [16] would imply that $\limsup_{t \to t_{\max}^-} V(t) = +\infty$ or $\liminf_{t \to t_{\max}^-} W(t) = -\infty$. Using induction and (A.11), (A.13), we are in a position to show that:

$$\min\left(W(0), \left(\frac{L-c}{1-c}\right)^i W(0)\right) \leq W(ih) \leq V(ih) \leq \max\left(V(0), \left(\frac{L-c}{1-c}\right)^i V(0)\right) \tag{A.14}$$

for all integers $0 \leq i \leq [h^{-1} t_{\max}]$, where $h = \min(\delta, A-\delta)$. Moreover, using the fact that $L \geq 1$, (A.14), (A.11) and (A.13) with $t = h[h^{-1} t_{\max}]$ for the case $t_{\max} \neq h[h^{-1} t_{\max}]$ or $t = t_{\max} - h$ for the case $t_{\max} = h[h^{-1} t_{\max}]$ and arbitrary $q \in [0, t_{\max} - t) \subseteq [0, h]$, we get

$$\min\left(W(0), \left(\frac{L-c}{1-c}\right)^{1+[h^{-1} t_{\max}]} W(0)\right) \leq \sup_{0 \leq t < t_{\max}} W(t) \leq \sup_{0 \leq t < t_{\max}} V(t) \leq \max\left(V(0), \left(\frac{L-c}{1-c}\right)^{1+[h^{-1} t_{\max}]} V(0)\right)$$

which contradicts the assertion that $\limsup_{t \to t_{\max}^-} V(t) = +\infty$ or $\liminf_{t \to t_{\max}^-} W(t) = -\infty$.

Inequality (4.2) is a direct consequence of definitions (A.8), the fact that $L \geq 1$ and inequalities (A.11), (A.13), (A.14). The proof is complete. ◁

**Proof of Corollary 4.2:** Since (4.1) holds and since $1 = \int_0^A \varphi(a) da$, we get:

$$x(t) - P(x_0) = \int_0^A \varphi(a)(x(t-a) - P(x_0)) da, \text{ for all } t \geq 0 \tag{A.15}$$

Let $\tilde{K}, \varepsilon > 0$ be the constants involved in (4.5). Using (4.5) and (A.15), we get for all $t \geq 0$:

$$|x(t) - P(x_0)| \leq \max_{0 \leq a \leq A}(\varphi(a)) \int_0^A |x(t-a) - P(x_0)| da \leq \max_{0 \leq a \leq A}(\varphi(a)) \tilde{K} \exp(-\varepsilon t) \int_0^A |x(-a)| da \tag{A.16}$$

It follows from (A.16) that the following estimate holds for all $t \geq A$:

$$\max_{-A \leq \theta \leq 0}(|x(t+\theta) - P(x_0)|) \leq \max_{0 \leq a \leq A}(\varphi(a)) \tilde{K} \exp(-\varepsilon(t-A)) A \max_{-A \leq a \leq 0}(|x_0(a)|) \tag{A.17}$$

When $0 \leq t < A$, we get from (A.16):

$$\max_{-t \leq \theta \leq 0}(|x(t+\theta) - P(x_0)|) = \max_{0 \leq s \leq t}(|x(s) - P(x_0)|) \leq \max_{0 \leq s \leq t}\left(\max_{0 \leq a \leq A}(\varphi(a)) \tilde{K} \exp(-\varepsilon s) \int_0^A |x(-a) da|\right)$$

$$= \max_{0 \leq a \leq A}(\varphi(a)) \tilde{K} \int_0^A |x_0(-a) da| \leq \max_{0 \leq a \leq A}(\varphi(a)) \tilde{K} A \max_{-A \leq a \leq 0}(|x_0(a)|)$$

Moreover, using definition (4.4) and the fact that $r := \left(\int_0^A a\varphi(a) da\right)^{-1}$ (which imply that $|P(x_0)| \leq \max_{-A \leq a \leq 0}(|x_0(a)|)$ for all $x_0 \in C^0([-A,0]; \Re)$), we get for $0 \leq t < A$:

$$\max_{-A \leq \theta \leq -t}(|x(t+\theta) - P(x_0)|) = \max_{t-A \leq s \leq 0}(|x(s) - P(x_0)|) \leq \max_{-A \leq s \leq 0}(|x(s) - P(x_0)|)$$

$$= \max_{-A \leq s \leq 0}(|x_0(s) - P(x_0)|) \leq \max_{-A \leq s \leq 0}(|x_0(s)| + |P(x_0)|) = \max_{-A \leq s \leq 0}(|x_0(s)|) + |P(x_0)| \leq 2 \max_{-A \leq a \leq 0}(|x_0(a)|)$$

The two above inequalities give for $0 \leq t < A$:

$$\max_{-A \leq \theta \leq 0}(|x(t+\theta) - P(x_0)|) \leq \max\left(\max_{-t \leq \theta \leq 0}(|x(t+\theta) - P(x_0)|), \max_{-A \leq \theta \leq -t}(|x(t+\theta) - P(x_0)|)\right)$$

$$\leq \max\left(\max_{0 \leq a \leq A}(\varphi(a)) \tilde{K} A, 2\right) \max_{-A \leq a \leq 0}(|x_0(a)|) \leq \exp(-\varepsilon(t-A)) \max\left(\max_{0 \leq a \leq A}(\varphi(a)) \tilde{K} A, 2\right) \max_{-A \leq a \leq 0}(|x_0(a)|) \tag{A.18}$$

Combining (A.17) and (A.18), we conclude that estimate (4.6) holds with $\sigma = \varepsilon$ and $M = \exp(\varepsilon A) \max\left(\max_{0 \leq a \leq A}(\varphi(a)) \tilde{K} A, 2\right)$. The proof is complete. ◁

**Proof of Lemma 4.3:** Notice that since $x \in C^0([-A,+\infty);\Re)$, it follows that the mapping $\Re_+ \ni t \to V(x_t)$ is continuous. We have for all $t,h \geq 0$ with $h < A$:

$$V(x_{t+h}) = \max_{a \in [0,A]} \left(\exp(-\sigma a)|x(t+h-a)|\right) = \max_{t+h-A \leq s \leq t+h} \left(\exp(-\sigma(t+h-s))|x(s)|\right)$$
$$= \max\left(\max_{t+h-A \leq s \leq t} \left(\exp(-\sigma(t+h-s))|x(s)|\right), \max_{t \leq s \leq t+h} \left(\exp(-\sigma(t+h-s))|x(s)|\right)\right)$$
$$= \exp(-\sigma h) \max\left(\max_{0 \leq a \leq A-h} \left(\exp(-\sigma a)|x(t-a)|\right), \max_{t \leq s \leq t+h} \left(\exp(-\sigma(t-s))|x(s)|\right)\right)$$
$$\leq \exp(-\sigma h) \max\left(V(x_t), \max_{t \leq s \leq t+h} \left(\exp(-\sigma(t-s))|x(s)|\right)\right) \quad (A.19)$$

Using (4.1) and (A.19) we obtain:

$$\exp(\sigma(t+h))V(x_{t+h}) \leq \max\left(\exp(\sigma t)V(x_t), \max_{t \leq s \leq t+h}\left(\exp(\sigma s)\left|\int_0^A \varphi(a)x(s-a)da\right|\right)\right)$$
$$\leq \max\left(\exp(\sigma t)V(x_t), \int_0^A |\varphi(a)|\exp(\sigma a)da \max_{t \leq s \leq t+h}\left(\exp(\sigma s) \max_{0 \leq a \leq A}\left(\exp(-\sigma a)|x(s-a)|\right)\right)\right) \quad (A.20)$$
$$= \max\left(\exp(\sigma t)V(x_t), \int_0^A |\varphi(a)|\exp(\sigma a)da \max_{t \leq s \leq t+h}\left(\exp(\sigma s)V(x_s)\right)\right)$$

Consequently, we obtain from (A.20) for all $t,h \geq 0$ with $h < A$:

$$\max_{t \leq s \leq t+h}\left(\exp(\sigma s)V(x_s)\right) \leq \max\left(\exp(\sigma t)V(x_t), \int_0^A |\varphi(a)|\exp(\sigma a)da \max_{t \leq s \leq t+h}\left(\exp(\sigma s)V(x_s)\right)\right) \quad (A.21)$$

Since $\int_0^A |\varphi(a)|\exp(\sigma a)da < 1$, we obtain from (A.21) for all $t,h \geq 0$ with $h < A$:

$$\max_{t \leq s \leq t+h}\left(\exp(\sigma s)V(x_s)\right) \leq \exp(\sigma t)V(x_t) \quad (A.22)$$

Indeed, the proof of (A.22) follows from distinguishing the cases (i) $\exp(\sigma t)V(x_t) \geq \int_0^A |\varphi(a)|\exp(\sigma a)da \max_{t \leq s \leq t+h}\left(\exp(\sigma s)V(x_s)\right)$, and (ii) $\exp(\sigma t)V(x_t) < \int_0^A |\varphi(a)|\exp(\sigma a)da \max_{t \leq s \leq t+h}\left(\exp(\sigma s)V(x_s)\right)$.

Case (ii) leads to a contradiction, since in this case (A.21) in conjunction with $\int_0^A |\varphi(a)|\exp(\sigma a)da < 1$ implies $\max_{t \leq s \leq t+h}\left(\exp(\sigma s)V(x_s)\right) = 0$, which contradicts the assumption $\exp(\sigma t)V(x_t) < \int_0^A |\varphi(a)|\exp(\sigma a)da \max_{t \leq s \leq t+h}\left(\exp(\sigma s)V(x_s)\right)$. Therefore, we obtain from (A.22) for all $t,h \geq 0$ with $h < A$:

$$V(x_{t+h}) \leq \exp(-\sigma h)V(x_t) \quad (A.23)$$

It follows from (A.23) for all $t \geq 0$ and $h \in (0,A)$:

$$h^{-1}(V(x_{t+h}) - V(x_t)) \leq h^{-1}(\exp(-\sigma h) - 1)V(x_t) \quad (A.24)$$

Letting $h \to 0^+$ and using (A.24), we obtain (4.7). The proof is complete. ◁

**Proof of Theorem 4.4:** Notice that since $x \in C^0([-A,+\infty);\Re)$, it follows that the mappings $\Re_+ \ni t \to V(x_t)$, $\Re_+ \ni t \to P(x_t)$ are continuous. Moreover, definition (4.4) implies that

$$P(x_t) = r\int_0^A x(t-a) \int_a^A \varphi(s)ds\,da = r\int_{t-A}^t x(w) \int_{t-w}^A \varphi(s)ds\,dw$$

It follows from Leibniz's rule that the mapping $(0,+\infty) \ni t \to P(x_t)$ is continuously differentiable and its derivative satisfies

$$\frac{d}{dt}P(x_t) = rx(t) - r\int_{t-A}^t x(w)\varphi(t-w)dw = r\left(x(t) - \int_0^A x(t-a)\varphi(a)da\right), \text{ for all } t \geq 0$$

Notice that for the derivation of the above equality we have used the fact that $\int_0^A \varphi(a)da = 1$. Using (4.1) we can conclude that (4.9) holds. Next, define

$$y(t) = x(t) - P(x_0), \text{ for all } t \geq -A \tag{A.25}$$

Using (4.1), definition (A.25) and the fact that $\int_0^A \varphi(a)da = 1$, we obtain:

$$y(t) = \int_0^A \varphi(a) y(t-a)da, \text{ for all } t \geq 0 \tag{A.26}$$

Moreover, it follows from definition (4.4) and (4.9) that:

$$0 = P(x_t) - P(x_0) = r\int_0^A x(t-a) \int_a^A \varphi(s)ds\,da - P(x_0), \text{ for all } t \geq 0 \tag{A.27}$$

Since $r\int_0^A\int_a^A \varphi(s)ds\,da = r\int_0^A a\varphi(a)da = 1$, we obtain from (A.25) and (A.27):

$$0 = \int_0^A y(t-a)\left(r\int_a^A \varphi(s)ds\right)da, \text{ for all } t \geq 0 \tag{A.28}$$

Combining (A.26) and (A.28), we get:

$$y(t) = \int_0^A y(t-a)\left(\varphi(a) - \lambda r\int_a^A \varphi(s)ds\right)da, \text{ for all } t \geq 0 \tag{A.29}$$

where $\lambda > 0$ is the real number for which $\int_0^A \left|\varphi(a) - r\lambda\int_a^A \varphi(s)ds\right|da < 1$. Therefore, we are in a position to apply Lemma 4.3 for the solution of (A.29). More specifically, we get:

$$\limsup_{h \to 0^+}\left(h^{-1}(W(y_{t+h}) - W(y_t))\right) \leq -\sigma W(y_t), \text{ for all } t \geq 0 \tag{A.30}$$

where $W(y_t) := \max_{a\in[0,A]}\left(\exp(-\sigma a)|y(t-a)|\right)$ and $\sigma > 0$ is a real number for which $\int_0^A \left|\varphi(a) - r\lambda\int_a^A \varphi(s)ds\right|\exp(\sigma a)da < 1$. Finally, we notice that definitions (4.8), (A.25) and equality (4.9) imply the following equalities:

$$V(x_t) = \max_{a\in[0,A]}\left(\exp(-\sigma a)|x(t-a) - P(x_t)|\right) = \max_{a\in[0,A]}\left(\exp(-\sigma a)|x(t-a) - P(x_0)|\right)$$
$$= \max_{a\in[0,A]}\left(\exp(-\sigma a)|y(t-a)|\right) = W(y_t) \tag{A.31}$$

The differential inequality (4.10) is a direct consequence of equation (A.31) and inequality (A.30). The proof is complete. ◁

**Proof of Corollary 4.6:** Working as in the proof of Theorem 4.4, we show that (4.14) holds. Since $C: S' \to [\kappa, +\infty)$ is a continuous functional and $x \in C^0([-A, +\infty); \Re)$, it follows that the mappings $\mathfrak{R}_+ \ni t \to C(x_t)$, $\mathfrak{R}_+ \ni t \to W(x_t)$ are continuous. Applying Theorem 4.4 and taking into account the fact that $P(x_t) = 0$ for all $t \geq 0$, we obtain:

$$\limsup_{h \to 0^+}\left(h^{-1}\left(W(x_{t+h}) - W(x_t)\right)\right) \leq -\sigma W(x_t), \text{ for all } t \geq 0 \tag{A.32}$$

Let $x(t) \in \mathfrak{R}$ be a solution of (4.1). The differential inequalities (4.12), (A.32) imply that the mappings $t \to C(x_t)$, $t \to W(x_t)$ are non-increasing. Consequently, we get $C(x_{t+h})b(W(x_{t+h})) \leq C(x_t)b(W(x_{t+h}))$ for all $h \geq 0$, which implies

$$h^{-1}\left(C(x_{t+h})b(W(x_{t+h})) - C(x_t)b(W(x_t))\right) \leq h^{-1}C(x_t)\left(b(W(x_{t+h})) - b(W(x_t))\right), \text{ for all } h \geq 0$$

By virtue of the mean value theorem, we obtain the existence of $s \in (0,1)$ such that

$$b(W(x_{t+h})) - b(W(x_t)) = \left(W(x_{t+h}) - W(x_t)\right)b'\left(sW(x_t) + (1-s)W(x_{t+h})\right).$$

Using the fact that the mapping $t \to W(x_t)$ is non-increasing and combining the above relations, we obtain

$$\begin{aligned} &h^{-1}\left(C(x_{t+h})b(W(x_{t+h})) - C(x_t)b(W(x_t))\right) \\ &\leq C(x_t)h^{-1}\left(W(x_{t+h}) - W(x_t)\right) \min_{s \in [0,1]}\left(b'\left(sW(x_t) + (1-s)W(x_{t+h})\right)\right) \end{aligned} \tag{A.33}$$

for all sufficiently small $h > 0$. The differential inequality (4.13) is a direct consequence of inequalities (A.32), (A.33) and the fact that the mapping $t \to W(x_t)$ is continuous. The proof is complete. ◁

**Proof of Lemma 5.1:** By virtue of Remark 4.7 and Corollary 4.6, for every $\psi_0 \in S$ the solution of the IDE $\psi(t) = \int_0^A \tilde{k}(a)\psi(t-a)da$ exists for all $t \geq 0$, is unique and satisfies $\psi_t \in S$ for all $t \geq 0$. More specifically, using Lemma 4.1, we can guarantee that the solution $\psi_t \in C^0([-A,+\infty);(-1,+\infty))$ of the IDE $\psi(t) = \int_0^A \tilde{k}(a)\psi(t-a)da$ satisfies

$$\inf_{t \geq -A}(\psi(t)) \geq \min_{-A \leq s \leq 0}(\psi_0(s)) > -1, \text{ for all } t \geq 0 \tag{A.34}$$

Indeed, since $\tilde{k} \in C^0([0,A];\mathfrak{R}_+)$ with $\int_0^A \tilde{k}(a)da = 1$, it follows that all assumptions of Lemma 4.1 hold. Therefore, we get from (4.2) with $L = 1$ and arbitrary $c \in (0,1)$ that the inequality

$$\inf_{-A \leq a < 0}(\psi_0(a)) \leq \inf_{-A \leq a < 0}(\psi(t+a)) \leq \sup_{-A \leq a < 0}(\psi(t+a)) \leq \sup_{-A \leq a < 0}(\psi_0(a))$$

holds for all $t \geq 0$. Inequality (A.34) is a direct consequence of continuity of $\psi_0 \in S$ and the above inequality.

Given the facts that $g \in C^0([0,A];\mathfrak{R})$ satisfies $g(a) \geq 0$ for all $a \in [0,A]$ with $\int_0^A g(a)da = 1$ and that the solution $\psi_t \in C^0([-A,+\infty);(-1,+\infty))$ of the IDE $\psi(t) = \int_0^A \tilde{k}(a)\psi(t-a)da$ satisfies (A.34), we are in a position to guarantee that the mapping $\mathfrak{R}_+ \ni t \to v(t) \in \mathfrak{R}$ defined by

$$v(t) = \ln\left(1 + \int_0^A g(a)\psi(t-a)da\right), \text{ for all } t \geq 0 \quad (A.35)$$

is well-defined and is a continuous mapping. It follows that for every $z_0 \in \Re^2$, $\eta_0 \in \Re$ the solution of the system of differential equations $\dot\eta(t) = D^* - sat(z_2(t) + \gamma\eta(t) + \gamma v(t))$, $\dot z_1(t) = z_2(t) - sat(z_2(t) + \gamma\eta(t) + \gamma v(t)) - l_1(z_1(t) - \eta(t) - v(t))$, $\dot z_2(t) = -l_2(z_1(t) - \eta(t) - v(t))$ exists locally and is unique. Moreover, due to the fact that the right hand side of the differential equations satisfies a linear growth condition, it follows that the solution $(z(t), \eta(t)) \in \Re^2 \times \Re$ of the system of differential equations $\dot\eta(t) = D^* - sat(z_2(t) + \gamma\eta(t) + \gamma v(t))$, $\dot z_1(t) = z_2(t) - sat(z_2(t) + \gamma\eta(t) + \gamma v(t)) - l_1(z_1(t) - \eta(t) - v(t))$, $\dot z_2(t) = -l_2(z_1(t) - \eta(t) - v(t))$ exists for all $t \geq 0$. Due to definition (A.35) and equations (5.6), (5.8), we are in a position to conclude that the constructed mappings coincide with a solution $(z(t), \eta(t), \psi_t) \in \Re^2 \times \Re \times S$ of the closed-loop system (5.5) with (5.7), (5.8) with initial condition $z_0 \in \Re^2$, $(\eta_0, \psi_0) \in \Re \times S$. Uniqueness of solution of the closed-loop system (5.5) with (5.7), (5.8) is a direct consequence of the above procedure of the construction of the solution.

Define

$$\begin{aligned} e_1(t) &= z_1(t) - \eta(t) \\ e_2(t) &= z_2(t) - D^* \end{aligned}, \text{ for all } t \geq 0 \quad (A.36)$$

and notice that equations (5.5), (5.6), (5.7), (A.35) and definition (A.36) allow us to conclude that the following differential equations hold for all $t \geq 0$:

$$\begin{aligned} \frac{d}{dt}\left(e_1^2(t) - p_1 e_1(t)e_2(t) + p_2 e_2^2(t)\right) &= -(2l_1 - l_2 p_1)e_1^2(t) + (2 + l_1 p_1 - 2p_2 l_2)e_1(t)e_2(t) - p_1 e_2^2(t) \\ &\quad + ((2l_1 - p_1 l_2)e_1(t) + (2p_2 l_2 - p_1 l_1)e_2(t))v(t) \end{aligned} \quad (A.37)$$

$$\begin{aligned} \dot e_1(t) &= e_2(t) - l_1 e_1(t) + l_1 v(t) \\ \dot e_2(t) &= -l_2 e_1(t) + l_2 v(t) \end{aligned} \quad (A.38)$$

Since the inequalities $l_1, l_2 > 0$, $p_1, p_2 > 0$, $(2 + l_1 p_1 - 2l_2 p_2)^2 < 8l_1 p_1 - 4l_2 p_1^2$, $p_1^2 < 4p_2$ hold, it follows that the quadratic forms $A(e) := e_1^2 - p_1 e_1 e_2 + p_2 e_2^2$, $B(e) := (2l_1 - l_2 p_1)e_1^2 - (2 + l_1 p_1 - 2p_2 l_2)e_1 e_2 + p_1 e_2^2$ are positive definite (recall Remark 2.2(vi)). It follows that there exist constants $K_2 \geq K_1 > 0$, $K_4 \geq K_3 > 0$ such that

$$K_1|e|^2 \leq A(e) = e_1^2 - p_1 e_1 e_2 + p_2 e_2^2 \leq K_2|e|^2, \text{ for all } e \in \Re^2 \quad (A.39)$$

$$K_3|e|^2 \leq B(e) = (2l_1 - l_2 p_1)e_1^2 - (2 + l_1 p_1 - 2p_2 l_2)e_1 e_2 + p_1 e_2^2 \leq K_4|e|^2, \text{ for all } e \in \Re^2 \quad (A.40)$$

Using (A.37), (A.39), (A.40) and the inequality $((2l_1 - p_1 l_2)e_1 + (2p_2 l_2 - p_1 l_1)e_2)v \leq \frac{K_3}{2}|e|^2 + \frac{(2l_1 - p_1 l_2)^2 + (2p_2 l_2 - p_1 l_1)^2}{2K_3}|v|^2$ (which holds for all $e \in \Re^2$ and $v \in \Re$), we conclude that there exist constants $\mu, c > 0$ such that the following differential inequality holds for all $t \geq 0$:

$$\frac{d}{dt}\left(e_1^2(t) - p_1 e_1(t)e_2(t) + p_2 e_2^2(t)\right) \leq -2\mu\left(e_1^2(t) - p_1 e_1(t)e_2(t) + p_2 e_2^2(t)\right) + c|v(t)|^2 \quad (A.41)$$

Since $\tilde{k} \in C^0([0,A]; \Re_+)$ satisfies assumption (H1) and since the mapping $t \to g_1(\psi_t)$ is non-decreasing for every solution of the IDE $\psi(t) = \int_0^A \tilde{k}(a)\psi(t-a)da$, where $g_1$ is the continuous functional $g_1(\psi) := \min_{a \in [0,A]}(\psi(-a))$ (recall Remark 4.7), it follows that the mapping $t \to C(\psi_t) := \left(1 + \min\left(0, \min_{a \in [0,A]}(\psi_t(-a))\right)\right)^{-2}$ is non-increasing. Using Remark 4.7 and Corollary 4.6 with $b(s) = Ms^2/2$, where $M > 0$ is an arbitrary constant, we conclude that

$$\limsup_{h \to 0^+}\left( h^{-1}\left( \frac{M}{2}\left( \frac{\max_{a \in [0,A]}\left(\exp(-\sigma a)|\psi_{t+h}(-a)|\right)}{1 + \min\left(0, \min_{a \in [0,A]}(\psi_{t+h}(-a))\right)}\right)^2 - \frac{M}{2}\left( \frac{\max_{a \in [0,A]}\left(\exp(-\sigma a)|\psi_t(-a)|\right)}{1 + \min\left(0, \min_{a \in [0,A]}(\psi_t(-a))\right)}\right)^2 \right)\right)$$

, for all $t \geq 0$  (A.42)

$$\leq -\sigma M \left( \frac{\max_{a \in [0,A]}\left(\exp(-\sigma a)|\psi_t(-a)|\right)}{1 + \min\left(0, \min_{a \in [0,A]}(\psi_t(-a))\right)}\right)^2$$

where $\sigma > 0$ is a real number for which $\int_0^A \left|\tilde{k}(a) - r\lambda \int_a^A \tilde{k}(s)ds\right| \exp(\sigma a)da < 1$, $\lambda > 0$ is a real number for which $\int_0^A \left|\tilde{k}(a) - r\lambda \int_a^A \tilde{k}(s)ds\right| da < 1$ and $r := \left(\int_0^A a\tilde{k}(a)da\right)^{-1}$.

Since $|\ln(1+x)| \leq \frac{|x|}{1 + \min(x,0)}$ for all $x > -1$ and using the facts that $g(a) \geq 0$ for all $a \in [0,A]$ and $\int_0^A g(a)da = 1$, we obtain from (A.35):

$$|v(t)| \leq \frac{\exp(\sigma A) \max_{a \in [0,A]}\left(|\psi_t(-a)|\exp(-\sigma a)\right)}{1 + \min\left(0, \min_{a \in [0,A]}(\psi_t(-a))\right)}, \text{ for all } t \geq 0 \quad (A.43)$$

Using (A.41), (A.42), (A.43) we obtain the following differential inequality

$$\limsup_{h \to 0^+}\left(h^{-1}\left(\tilde{Q}(e(t+h),\psi_{t+h}) - \tilde{Q}(e(t),\psi_t)\right)\right)$$

$$\leq -2\mu\left(e_1^2(t) - p_1 e_1(t)e_2(t) + p_2 e_2^2(t)\right) - (\sigma M - c\exp(2\sigma A))\left(\frac{\max_{a \in [0,A]}\left(\exp(-\sigma a)|\psi_t(-a)|\right)}{1 + \min\left(0, \min_{a \in [0,A]}(\psi_t(-a))\right)}\right)^2, \text{ for all } t \geq 0 \quad (A.44)$$

where

$$\tilde{Q}(e,\psi) := e_1^2 - p_1 e_1 e_2 + p_2 e_2^2 + \frac{M}{2}\left(\frac{\max_{a \in [0,A]}\left(\exp(-\sigma a)|\psi(-a)|\right)}{1 + \min\left(0, \min_{a \in [0,A]}(\psi(-a))\right)}\right)^2 \quad (A.45)$$

Selecting $M > \sigma^{-1}c\exp(2\sigma A)$, we obtain from (A.44) and definition (A.45):

$$\limsup_{h \to 0^+}\left(h^{-1}\left(\tilde{Q}(e(t+h),\psi_{t+h}) - \tilde{Q}(e(t),\psi_t)\right)\right) \leq -2\tilde{\mu}\tilde{Q}(e(t),\psi_t), \text{ for all } t \geq 0 \quad (A.46)$$

where $\tilde{\mu} = \min\left(\mu, \frac{\sigma M - c\exp(2\sigma A)}{2M}\right)$.

Suppose that $\tilde{Q}(e(t),\psi_t) > 0$. Since the mapping $t \to \tilde{Q}(e(t),\psi_t)$ is continuous, it follows that $\tilde{Q}(e(t+h),\psi_{t+h})$ for all sufficiently small $h > 0$. The differential inequality (A.46) implies that the mapping $t \to \tilde{Q}(e(t),\psi_t)$ is non-increasing. Consequently, by virtue of the mean value theorem, we obtain:

$$h^{-1}\left(\sqrt{\tilde{Q}(e(t+h),\psi_{t+h}))} - \sqrt{\tilde{Q}(e(t),\psi_t))}\right) \leq \frac{\tilde{Q}(e(t+h),\psi_{t+h})) - \tilde{Q}(e(t),\psi_t))}{2h\sqrt{\tilde{Q}(e(t),\psi_t))}} \quad (A.47)$$

for all sufficiently small $h > 0$. Therefore, using (A.46), (A.47), we obtain the differential inequality

$$\limsup_{h \to 0^+}\left(h^{-1}\left(\sqrt{\tilde{Q}(e(t+h),\psi_{t+h})} - \sqrt{\tilde{Q}(e(t),\psi_t)}\right)\right) \leq -\tilde{\mu}\sqrt{\tilde{Q}(e(t),\psi_t)} \quad (A.48)$$

for all $t \geq 0$ with $\tilde{Q}(e(t),\psi_t) > 0$. On the other hand, using (A.39) and (A.45), we are in a position to conclude that $e(t) = 0$, $\psi_t = 0$ when $\tilde{Q}(e(t),\psi_t) = 0$. In this case and using (5.5), we conclude that $\psi_{t+h} = 0$ for all $h \geq 0$ (the unique solution of $\psi(t) = \int_0^A \tilde{k}(a)\psi(t-a)da$). Therefore, in this case (A.35) implies that $v(t+h) = 0$ for all $h \geq 0$. Finally, (A.38) and the facts that $e(t) = 0$ and $v(t+h) = 0$ for all $h \geq 0$, imply that $e(t+h) = 0$ for all $h \geq 0$, when $\tilde{Q}(e(t),\psi_t) = 0$. Definition (A.45) allows us to conclude that $\tilde{Q}(e(t+h),\psi_{t+h}) = 0$ for all $h \geq 0$, when $\tilde{Q}(e(t),\psi_t) = 0$. Therefore, the differential inequality (A.48) holds for all $t \geq 0$.

Finally, using (5.2), (5.5), (5.6), (5.8), (A.35), (A.36) we get for all $t \geq 0$:

$$\frac{d}{dt}\left(\eta^2(t)\right) = 2\eta(t)\left(D^* - sat(z_2(t) + \gamma\eta(t) + \gamma v(t))\right) = -2\eta(t)q(e_2(t) + \gamma\eta(t) + \gamma v(t)) \quad (A.49)$$

We distinguish the following cases:

-Case 1: $|e_2(t) + \gamma v(t)| \leq \frac{\gamma}{2}|\eta(t)|$ and $\eta(t) \geq 0$. In this case, we have $\frac{\gamma}{2}\eta(t) \leq e_2(t) + \gamma v(t) + \gamma\eta(t) \leq \frac{3\gamma}{2}\eta(t)$. Using the fact that the function $q$ defined by (5.1) is non-decreasing and the previous inequality we obtain from (A.49) that $\frac{d}{dt}\left(\eta^2(t)\right) \leq -2\eta(t)q\left(\frac{\gamma}{2}\eta(t)\right)$.

-Case 2: $|e_2(t) + \gamma v(t)| \leq \frac{\gamma}{2}|\eta(t)|$ and $\eta(t) < 0$. In this case, we have $\frac{3\gamma}{2}\eta(t) \leq e_2(t) + \gamma v(t) + \gamma\eta(t) \leq \frac{\gamma}{2}\eta(t)$. Using the fact that the function $q$ defined by (5.1) is non-decreasing and the previous inequality we obtain from (A.49) that $\frac{d}{dt}\left(\eta^2(t)\right) \leq -2\eta(t)q\left(\frac{\gamma}{2}\eta(t)\right)$.

-Case 3: $|e_2(t) + \gamma v(t)| > \frac{\gamma}{2}|\eta(t)|$. Inequality (5.3) implies that $|q(e_2(t) + \gamma\eta(t) + \gamma v(t))| \leq \max(D_{\max} - D^*, D^* - D_{\min})$ and $\left|q\left(\frac{\gamma}{2}\eta(t)\right)\right| \leq \max(D_{\max} - D^*, D^* - D_{\min})$. Consequently, we obtain from (A.49) that

$$\frac{d}{dt}\left(\eta^2(t)\right) + 2\eta(t)q\left(\frac{\gamma}{2}\eta(t)\right) \leq 4|\eta(t)|\max(D_{\max} - D^*, D^* - D_{\min})$$

$$\leq \frac{8}{\gamma}|e_2(t) + \gamma v(t)|\max(D_{\max} - D^*, D^* - D_{\min})$$

Combining all the above three cases, we conclude that

$$\frac{d}{dt}\left(\eta^2(t)\right) \le -2\eta(t)q\left(\frac{\gamma}{2}\eta(t)\right) + R|e_2(t) + \gamma v(t)|, \text{ for all } t \ge 0 \tag{A.50}$$

where $R := 8\gamma^{-1}\max(D_{max} - D^*, D^* - D_{min})$. Combining (A.46), (A.48) with (A.50) and using the triangle inequality, we obtain for every $G, \beta \ge 0$:

$$\limsup_{h \to 0^+}\left(h^{-1}\left(\eta^2(t+h) + G\sqrt{\tilde{Q}(e(t+h),\psi_{t+h})} + \beta\tilde{Q}(e(t+h),\psi_{t+h}) - \eta^2(t) - G\sqrt{\tilde{Q}(e(t),\psi_t)} - \beta\tilde{Q}(e(t),\psi_t)\right)\right)$$
$$\le -2\eta(t)q\left(\frac{\gamma}{2}\eta(t)\right) + R|e_2(t)| + R\gamma|v(t)| - \tilde{\mu}G\sqrt{\tilde{Q}(e(t),\psi_t)} - 2\tilde{\mu}\beta\tilde{Q}(e(t),\psi_t)$$
$$\text{for all } t \ge 0 \tag{A.51}$$

Using (A.39), (A.43) and definition (A.45), we obtain:

$$\max\left(\sqrt{K_1}|e_2(t)|, \exp(-\sigma A)\sqrt{\frac{M}{2}}|v(t)|\right) \le \sqrt{\tilde{Q}(e(t),\psi_t)}, \text{ for all } t \ge 0 \tag{A.52}$$

Using (A.51) and (A.52), we obtain for every $G, \beta \ge 0$:

$$\limsup_{h \to 0^+}\left(h^{-1}\left(\eta^2(t+h) + G\sqrt{\tilde{Q}(e(t+h),\psi_{t+h})} + \beta\tilde{Q}(e(t+h),\psi_{t+h}) - \eta^2(t) - G\sqrt{\tilde{Q}(e(t),\psi_t)} - \beta\tilde{Q}(e(t),\psi_t)\right)\right)$$
$$\le -2\eta(t)q\left(\frac{\gamma}{2}\eta(t)\right) - \left(\tilde{\mu}G - \frac{R}{\sqrt{K_1}} - \frac{R\gamma\sqrt{2}}{\sqrt{M}}\exp(\sigma A)\right)\sqrt{\tilde{Q}(e(t),\psi_t)} - 2\tilde{\mu}\beta\tilde{Q}(e(t),\psi_t)$$
$$\text{for all } t \ge 0 \tag{A.53}$$

Therefore, we obtain from (A.53), (5.4) and definitions (5.10), (5.11), (A.36), (A.45) for $G > \frac{R}{\tilde{\mu}\sqrt{K_1}} + \frac{R\gamma\sqrt{2}}{\tilde{\mu}\sqrt{M}}\exp(\sigma A)$ the differential inequality (5.9) with $L := \min\left(\tilde{\mu} - \frac{R}{G\sqrt{K_1}} - \frac{R\gamma\sqrt{2}}{G\sqrt{M}}\exp(\sigma A), \min(2,\gamma)\min(1, D_{max} - D^*, D^* - D_{min})\right)$. The proof is complete. ◁

**Proof of Lemma 5.2:** First we notice that the differential inequality (5.12) shows that $\varphi: \Re_+ \to \Re_+$ is non-increasing. We also make the following claim.

Claim: $\varphi(t) \le 1$ for all $t \ge T$, where $T = 2L^{-1}\max(0, \varphi(0) - 1)$.

If $\varphi(0) \le 1$ then the claim holds by virtue of the fact that $\varphi: \Re_+ \to \Re_+$ is non-increasing.

If $\varphi(0) > 1$ then the proof of the claim is made by contradiction. Suppose that there exists $t \ge T$ with $\varphi(t) > 1$. Since $\varphi: \Re_+ \to \Re_+$ is non-increasing, it follows that $\varphi(\tau) > 1$ for all $\tau \in [0, t]$. Consequently, we obtain from (5.12):

$$\limsup_{h \to 0^+}\left(h^{-1}(\varphi(\tau+h) - \varphi(\tau))\right) \le -L\frac{\varphi(\tau)}{1+\sqrt{\varphi(\tau)}} \le -\frac{L}{2}, \text{ for all } \tau \in [0, t] \tag{A.54}$$

Using the Comparison Lemma on page 85 in [19] and (A.54), we obtain $\varphi(\tau) \le \varphi(0) - \frac{L}{2}\tau$ for all $\tau \in [0, t]$. Since $\varphi(T) \le \varphi(0) - \frac{L}{2}T \le 1$, we obtain a contradiction.

Since $\varphi(t) \le 1$ for all $t \ge T$, we obtain from (5.12):

$$\limsup_{h \to 0^+}\left(h^{-1}\left(\varphi(t+h)-\varphi(t)\right)\right) \leq -\frac{L}{2}\varphi(t), \text{ for all } t \geq T \qquad (A.55)$$

Using the Comparison Lemma on page 85 in [19] and (A.55), we obtain $\varphi(t) \leq \varphi(T)\exp\left(-\frac{L}{2}(t-T)\right)$ for all $t \geq T$. Using the fact that $T = 2L^{-1}\max(0, \varphi(0)-1)$ (which implies the fact that $T=0$ when $\varphi(0) \leq 1$ and $\varphi(T) \leq 1$ when $\varphi(0) > 1$), we obtain the estimate (5.13) for all $t \geq T$. Since $\varphi: \Re_+ \to \Re_+$ is non-increasing and satisfies $\varphi(t) \leq \varphi(0)$ for all $t \in [0,T]$, we conclude that (5.13) holds for all $t \geq 0$. The proof is complete. ◁

**Proof of Lemma 5.3:** By virtue of Remark 4.7 and Corollary 4.6, for every $\psi_0 \in S$ the solution of the IDE $\psi(t) = \int_0^A \tilde{k}(a)\psi(t-a)da$ exists for all $t \geq 0$, is unique and satisfies $\psi_t \in S$ for all $t \geq 0$. More specifically, using Lemma 4.1, we can guarantee that the solution $\psi_t \in C^0([-A,+\infty);(-1,+\infty))$ of the IDE $\psi(t) = \int_0^A \tilde{k}(a)\psi(t-a)da$ satisfies (A.34). Working as in the proof of Theorem 4.4, we can also show that

$$P(\psi_t) = 0, \text{ for all } t \geq 0 \qquad (A.56)$$

Given the facts that $g \in C^0([0,A];\Re)$ satisfies $g(a) \geq 0$ for all $a \in [0,A]$ with $\int_0^A g(a)da = 1$ and that the solution $\psi_t \in C^0([-A,+\infty);(-1,+\infty))$ of the IDE $\psi(t) = \int_0^A \tilde{k}(a)\psi(t-a)da$ satisfies (A.34), we are in a position to guarantee that the mapping $\Re_+ \ni t \to v(t) \in \Re$ defined by (A.35) is well-defined and is a continuous mapping. It follows that for every $\eta_0 \in \Re$ the solution of the differential equation $\dot{\eta}(t) = D^* - sat\left(D^* + T^{-1}\eta(iT) + T^{-1}v(iT)\right)$, for all integers $i \geq 0$ and $t \in [iT, (i+1)T)$ with initial condition $\eta(0) = \eta_0$ exists locally and is unique. Moreover, due to the fact that the right hand side of the differential equation is bounded, it follows that the solution $\eta(t) \in \Re$ of the differential equation $\dot{\eta}(t) = D^* - sat\left(D^* + T^{-1}\eta(iT) + T^{-1}v(iT)\right)$, for all integers $i \geq 0$ and $t \in [iT, (i+1)T)$ with initial condition $\eta(0) = \eta_0$ exists for all $t \geq 0$. Due to definition (A.35) and equations (5.6), (5.34), we are in a position to conclude that the constructed mappings coincide with a solution $(\eta(t), \psi_t) \in \Re \times S$ of the closed-loop system (5.5), (5.6) with (5.34) with initial condition $(\eta_0, \psi_0) \in \Re \times S$. Uniqueness of solution of the closed-loop system (5.5), (5.6) with (5.34) is a direct consequence of the above procedure of the construction of the solution.

Using Corollary 4.2 with $\varphi(a) = \tilde{k}(a)$ for all $a \in [0,A]$, we conclude that there exist constants $M, \sigma > 0$ such that for every $\psi_0 \in S$ with $(G\psi_0) \in PC^1([0,A];\Re)$, the unique solution of the IDE $\psi(t) = \int_0^A \tilde{k}(a)\psi(t-a)da$ with initial condition $\psi(-a) = \psi_0(-a)$, for all $a \in [0,A]$ satisfies the following estimate for all $t \geq 0$:

$$\max_{-A \leq \theta \leq 0}\left(|\psi(t+\theta)|\right) \leq M\exp(-\sigma t)\max_{-A \leq a \leq 0}\left(|\psi_0(a)|\right) \qquad (A.57)$$

It follows from (5.5), (5.6), (A.35) and (5.34) that the following equation holds for all integers $i \geq 0$ and $t \in [iT, (i+1)T)$:

$$\eta(t) = \eta(iT) - (D_i - D^*)(t - iT) \qquad (A.58)$$

where

$$D_i = sat\left(D^* + T^{-1}\eta(iT) + T^{-1}v(iT)\right), \text{ for all integers } i \geq 0 \tag{A.59}$$

We next show the following claim.

<u>Claim 1:</u> The following inequality holds for all integers $i \geq 0$:

$$|\eta((i+1)T)| \leq |\eta(iT)| - \min(|\eta(iT)|, 2\delta) + 2|v(iT)| \tag{A.60}$$

where $\delta := \frac{1}{2}\min\left((D_{\max} - D^*)T, (D^* - D_{\min})T\right) > 0$.

<u>Proof of Claim 1:</u> We distinguish the following cases.

Case 1: $D_{\min} \leq D^* + T^{-1}\eta(iT) + T^{-1}v(iT) \leq D_{\max}$.
Definition (A.59) implies that $D_i = D^* + T^{-1}\eta(iT) + T^{-1}v(iT)$. Using (A.58), we get $\eta((i+1)T) = -v(iT)$, which directly implies (A.60).

Case 2: $D^* + T^{-1}\eta(iT) + T^{-1}v(iT) < D_{\min}$.
Definition (A.59) implies that $D_i = D_{\min}$. Using (A.58), we get $\eta((i+1)T) = \eta(iT) - (D_{\min} - D^*)T$. The inequality $D^* + T^{-1}\eta(iT) + T^{-1}v(iT) < D_{\min}$ implies that $-(D_{\min} - D^*)T + \eta(iT) + v(iT) < 0$. Thus, we get:

$$\begin{aligned}|\eta((i+1)T)| &= |\eta(iT) - (D_{\min} - D^*)T + v(iT) - v(iT)| \\ &\leq |\eta(iT) - (D_{\min} - D^*)T + v(iT)| + |v(iT)| \\ &= -\eta(iT) + (D_{\min} - D^*)T - v(iT) + |v(iT)| \\ &\leq |\eta(iT)| - (D^* - D_{\min})T + 2|v(iT)|\end{aligned}$$

The above inequality in conjunction with the fact that $\delta := \frac{1}{2}\min\left((D_{\max} - D^*)T, (D^* - D_{\min})T\right) > 0$ implies that (A.60) holds.

Case 3: $D^* + T^{-1}\eta(iT) + T^{-1}v(iT) > D_{\max}$.
Definition (A.59) implies that $D_i = D_{\max}$. Using (A.58) we get $\eta((i+1)T) = \eta(iT) - (D_{\max} - D^*)T$. The inequality $D^* + T^{-1}\eta(iT) + T^{-1}v(iT) > D_{\max}$ implies that $-(D_{\max} - D^*)T + \eta(iT) + v(iT) > 0$. Thus, we get:

$$\begin{aligned}|\eta((i+1)T)| &= |\eta(iT) - (D_{\max} - D^*)T + v(iT) - v(iT)| \\ &\leq |\eta(iT) - (D_{\max} - D^*)T + v(iT)| + |v(iT)| \\ &= \eta(iT) - (D_{\max} - D^*)T + v(iT) + |v(iT)| \\ &\leq |\eta(iT)| - (D_{\max} - D^*)T + 2|v(iT)|\end{aligned}$$

The above inequality in conjunction with the fact that $\delta := \frac{1}{2}\min\left((D_{\max} - D^*)T, (D^* - D_{\min})T\right) > 0$ implies that (A.60) holds.

The proof of Claim 1 is complete. ◁

<u>Claim 2:</u> The following inequalities hold for all integers $i \geq 0$:

$$\eta(iT) \geq \min\left(0, \eta_0 + i(D^* - D_{\min})T\right) + \min_{k=0,\ldots i}\left(\min(0, -v(kT))\right)$$
$$\eta(iT) \leq \max\left(0, \eta_0 - i(D_{\max} - D^*)T\right) + \max_{k=0,\ldots i}\left(\max(0, -v(kT))\right)$$
(A.61)

<u>Proof of Claim 2:</u> The proof of inequalities (A.61) is made by induction. First notice that both inequalities (A.61) hold for $i = 0$. Next assume that inequalities (A.61) hold for certain integer $i \geq 0$. We distinguish the following cases.

Case 1: $D_{\min} \leq D^* + T^{-1}\eta(iT) + T^{-1}v(iT) \leq D_{\max}$.
Definition (A.59) implies that $D_i = D^* + T^{-1}\eta(iT) + T^{-1}v(iT)$. Using (A.58), we get $\eta((i+1)T) = -v(iT)$. Consequently, we get:

$$\eta((i+1)T) = -v(iT) \leq \max(0, -v(iT)) \leq \max_{k=0,\ldots i+1}\left(\max(0, -v(kT))\right)$$
$$\leq \max\left(0, \eta_0 - (i+1)(D_{\max} - D^*)T\right) + \max_{k=0,\ldots i+1}\left(\max(0, -v(kT))\right)$$

which directly implies the second inequality (A.61) with $i+1$ in place of $i \geq 0$. Similarly, we obtain the first inequality (A.61) with $i+1$ in place of $i \geq 0$.

Case 2: $D^* + T^{-1}\eta(iT) + T^{-1}v(iT) < D_{\min}$.
Definition (A.59) implies that $D_i = D_{\min}$. Using (A.58) we get $\eta((i+1)T) = \eta(iT) - (D_{\min} - D^*)T$. Consequently, we get from (A.61):

$$\eta((i+1)T) = \eta(iT) + (D^* - D_{\min})T$$
$$\geq \min\left(0, \eta_0 + i(D^* - D_{\min})T\right) + (D^* - D_{\min})T + \min_{k=0,\ldots i}\left(\min(0, -v(kT))\right)$$
$$\geq \min\left((D^* - D_{\min})T, \eta_0 + (i+1)(D^* - D_{\min})T\right) + \min_{k=0,\ldots i+1}\left(\min(0, -v(kT))\right)$$
$$\geq \min\left(0, \eta_0 + (i+1)(D^* - D_{\min})T\right) + \min_{k=0,\ldots i+1}\left(\min(0, -v(kT))\right)$$

which is the first inequality (A.61) with $i+1$ in place of $i \geq 0$. Furthermore, the inequality $D^* + T^{-1}\eta(iT) + T^{-1}v(iT) < D_{\min}$ implies that $-(D_{\min} - D^*)T + \eta(iT) < -v(iT)$. Consequently, we get:

$$\eta((i+1)T) = \eta(iT) + (D^* - D_{\min})T \leq -v(iT) \leq \max(0, -v(iT)) \leq \max_{k=0,\ldots i+1}\left(\max(0, -v(kT))\right)$$
$$\leq \max\left(0, \eta_0 - (i+1)(D_{\max} - D^*)T\right) + \max_{k=0,\ldots i+1}\left(\max(0, -v(kT))\right)$$

which is the second inequality (A.61) with $i+1$ in place of $i \geq 0$.

Case 3: $D^* + T^{-1}\eta(iT) + T^{-1}v(iT) > D_{\max}$.
Definition (A.59) implies that $D_i = D_{\max}$. Using (A.58) we get $\eta((i+1)T) = \eta(iT) - (D_{\max} - D^*)T$. Consequently, we get from (A.61):

$$\eta((i+1)T) = \eta(iT) - (D_{max} - D^*)T$$
$$\leq \max\left(0, \eta_0 - i(D_{max} - D^*)T\right) - (D_{max} - D^*)T + \max_{k=0,\ldots i}\left(\max(0, -v(kT))\right)$$
$$\leq \max\left(-(D_{max} - D^*)T, \eta_0 - (i+1)(D_{max} - D^*)T\right) + \max_{k=0,\ldots i+1}\left(\max(0, -v(kT))\right)$$
$$\leq \max\left(0, \eta_0 - (i+1)(D_{max} - D^*)T\right) + \max_{k=0,\ldots i+1}\left(\max(0, -v(kT))\right)$$

which is the second inequality (A.61) with $i+1$ in place of $i \geq 0$. Furthermore, the inequality $D^* + T^{-1}\eta(iT) + T^{-1}v(iT) > D_{max}$ implies that $-(D_{max} - D^*)T + \eta(iT) > -v(iT)$. Consequently, we get:

$$\eta((i+1)T) = \eta(iT) - (D_{max} - D^*)T \geq -v(iT) \geq \min(0, -v(iT)) \geq \min_{k=0,\ldots i+1}\left(\min(0, -v(kT))\right)$$
$$\geq \min\left(0, \eta_0 + (i+1)(D^* - D_{min})T\right) + \min_{k=0,\ldots i+1}\left(\min(0, -v(kT))\right)$$

which is the first inequality (A.61) with $i+1$ in place of $i \geq 0$.

The proof of Claim 2 is complete. ◁

We next show the following claim.

<u>Claim 3:</u> The following inequalities hold for all $t \geq 0$:

$$\min(0, \eta_0) + 2 \min_{k=0,\ldots [t/T]}\left(\min(0, -v(kT))\right) \leq \eta(t) \leq \max(0, \eta_0) + 2 \max_{k=0,\ldots [t/T]}\left(\max(0, -v(kT))\right) \quad (A.62)$$

$$|\eta(t)| \leq |\eta([t/T]T)| + |v([t/T]T)| \quad (A.63)$$

<u>Proof of Claim 3:</u> Let arbitrary $t \geq 0$ and define $i = [t/T]$. Notice that the definition $i = [t/T]$ implies the inclusion $t \in [iT, (i+1)T)$. We distinguish the following cases.

Case 1: $D_{min} \leq D^* + T^{-1}\eta(iT) + T^{-1}v(iT) \leq D_{max}$.

Definition (A.59) implies that $D_i = D^* + T^{-1}\eta(iT) + T^{-1}v(iT)$. Using (A.58), we get for all $s \in [iT, (i+1)T]$:

$$\eta(s) = \left(1 - (s - iT)T^{-1}\right)\eta(iT) - (s - iT)T^{-1}v(iT)$$

The above equality in conjunction with the facts that $0 \leq 1 - (s-iT)T^{-1} \leq 1$, $0 \leq (s-iT)T^{-1} \leq 1$ and inequality (A.61) gives estimates (A.62), (A.63). More specifically, the above inequality in conjunction with the facts that $0 \leq 1 - (s-iT)T^{-1} \leq 1$, $0 \leq (s-iT)T^{-1} \leq 1$ implies for all $s \in [iT, (i+1)T]$

$$|\eta(s)| \leq \left(1 - (s-iT)T^{-1}\right)|\eta(iT)| + (s-iT)T^{-1}|v(iT)| \leq |\eta(iT)| + |v(iT)|$$

and since $t \in [iT, (i+1)T)$, the above inequality shows that (A.63) holds in this case. Moreover, the equation $\eta(s) = \left(1 - (s-iT)T^{-1}\right)\eta(iT) - (s-iT)T^{-1}v(iT)$ gives for all $s \in [iT, (i+1)T]$

$$\eta(s) \leq \left(1 - (s-iT)T^{-1}\right)\max(0, \eta(iT)) + (s-iT)T^{-1}\max(0, -v(iT))$$

which in conjunction with the facts that $0 \leq 1 - (s-iT)T^{-1} \leq 1$, $0 \leq (s-iT)T^{-1} \leq 1$ implies for all $s \in [iT, (i+1)T]$

$$\eta(s) \leq \max(0, \eta(iT)) + \max(0, -v(iT))$$

On the other hand, inequality (A.61) gives

$$\max(0,\eta(iT)) \leq \max\left(0,\eta_0 - i(D_{\max} - D^*)T\right) + \max_{k=0,\ldots i}\left(\max(0,-v(kT))\right)$$
$$\leq \max(0,\eta_0) + \max_{k=0,\ldots i}\left(\max(0,-v(kT))\right)$$

Combining the two above inequalities, we get for all $s \in [iT,(i+1)T]$

$$\eta(s) \leq \max(0,\eta_0) + 2\max_{k=0,\ldots i}\left(\max(0,-v(kT))\right)$$

and since $t \in [iT,(i+1)T]$, the above inequality shows that the right inequality (A.62) holds in this case. The left inequality (A.62) is proved in the same way.

Case 2: $D^* + T^{-1}\eta(iT) + T^{-1}v(iT) < D_{\min}$.
Definition (A.59) implies that $D_i = D_{\min}$. Using (A.58) we get $\eta(s) = \eta(iT) - (D_{\min} - D^*)(s - iT)$ for all $s \in [iT,(i+1)T]$. Therefore, we get for all $s \in [iT,(i+1)T]$

$$\eta(s) \geq \eta(iT)$$

which combined with (A.61) gives for all $s \in [iT,(i+1)T]$

$$\eta(s) \geq \min\left(0,\eta_0 + i(D^* - D_{\min})T\right) + \min_{k=0,\ldots i}\left(\min(0,-v(kT))\right)$$
$$\geq \min(0,\eta_0) + \min_{k=0,\ldots i}\left(\min(0,-v(kT))\right)$$
$$\geq \min(0,\eta_0) + 2\min_{k=0,\ldots i}\left(\min(0,-v(kT))\right)$$

Since $t \in [iT,(i+1)T]$, the above inequality shows that the left inequality (A.62) holds in this case. Furthermore, the inequality $D^* + T^{-1}\eta(iT) + T^{-1}v(iT) < D_{\min}$ implies that $-(D_{\min} - D^*)T + \eta(iT) < -v(iT)$. Since $\eta(s) = \eta(iT) - (D_{\min} - D^*)(s - iT)$ for all $s \in [iT,(i+1)T]$, we get for all $s \in [iT,(i+1)T]$:

$$\eta(s) = \eta(iT) - (D_{\min} - D^*)(s - iT) = \left(1 - (s-iT)T^{-1}\right)\eta(iT) + (s-iT)T^{-1}\left(\eta(iT) - (D_{\min} - D^*)T\right)$$
$$\leq \left(1 - (s-iT)T^{-1}\right)\eta(iT) - (s-iT)T^{-1}v(iT)$$
$$\leq \left(1 - (s-iT)T^{-1}\right)\max(0,\eta(iT)) + (s-iT)T^{-1}\max(0,-v(iT))$$

The above equality in conjunction with (A.61) and the facts that $0 \leq 1-(s-iT)T^{-1} \leq 1$, $0 \leq (s-iT)T^{-1} \leq 1$ shows that the right inequality (A.62) holds (exactly as in Case 1). Finally, we notice that the following inequalities hold for all $s \in [iT,(i+1)T]$

$$\eta(iT) \leq \eta(s) \leq \left(1-(s-iT)T^{-1}\right)\max(0,\eta(iT)) + (s-iT)T^{-1}\max(0,-v(iT))$$

which combined with the facts that $0 \leq 1-(s-iT)T^{-1} \leq 1$, $0 \leq (s-iT)T^{-1} \leq 1$ implies for all $s \in [iT,(i+1)T]$

$$-|\eta(iT)| - |v(iT)| \leq \eta(s) \leq |\eta(iT)| + |v(iT)|$$

Since $t \in [iT,(i+1)T]$, the above inequality shows that inequality (A.63) holds in this case.

Case 3: $D^* + T^{-1}\eta(iT) + T^{-1}v(iT) > D_{\max}$.
Definition (A.59) implies that $D_i = D_{\max}$. Using (A.58) we get $\eta(s) = \eta(iT) - (D_{\max} - D^*)(s - iT)$ for all $s \in [iT,(i+1)T]$. Therefore, we get for all $s \in [iT,(i+1)T]$

$$\eta(s) \leq \eta(iT)$$

which combined with (A.61) gives for all $s \in [iT, (i+1)T]$

$$\eta(s) \leq \max\left(0, \eta_0 - i(D_{max} - D^*)T\right) + \max_{k=0,\ldots i}\left(\max(0, -v(kT))\right)$$
$$\leq \max(0, \eta_0) + \max_{k=0,\ldots i}\left(\max(0, -v(kT))\right)$$
$$\leq \max(0, \eta_0) + 2\max_{k=0,\ldots i}\left(\max(0, -v(kT))\right)$$

Since $t \in [iT, (i+1)T)$, the above inequality shows that the right inequality (A.62) holds in this case.

Furthermore, the inequality $D^* + T^{-1}\eta(iT) + T^{-1}v(iT) > D_{max}$ implies that $-(D_{max} - D^*)T + \eta(iT) > -v(iT)$. Since $\eta(s) = \eta(iT) - (D_{max} - D^*)(s - iT)$ for all $s \in [iT, (i+1)T]$, we get for all $s \in [iT, (i+1)T]$:

$$\eta(s) = \eta(iT) - (D_{max} - D^*)(s - iT)$$
$$= \left(1 - (s - iT)T^{-1}\right)\eta(iT) + (s - iT)T^{-1}\left(\eta(iT) - (D_{max} - D^*)T\right)$$
$$\geq \left(1 - (s - iT)T^{-1}\right)\eta(iT) - (s - iT)T^{-1}v(iT)$$

The above equality in conjunction with the facts that $0 \leq 1 - (s - iT)T^{-1} \leq 1$, $0 \leq (s - iT)T^{-1} \leq 1$ gives for all $s \in [iT, (i+1)T]$:

$$\eta(s) \geq \left(1 - (s - iT)T^{-1}\right)\min(0, \eta(iT)) + (s - iT)T^{-1}\min(0, -v(iT))$$
$$\geq \min(0, \eta(iT)) + \min(0, -v(iT))$$

On the other hand, inequality (A.61) gives

$$\min(0, \eta(iT)) \geq \min\left(0, \eta_0 + i(D^* - D_{min})T\right) + \min_{k=0,\ldots i}\left(\min(0, -v(kT))\right)$$
$$\geq \min(0, \eta_0) + \min_{k=0,\ldots i}\left(\min(0, -v(kT))\right)$$

Combining the two above inequalities, we get for all $s \in [iT, (i+1)T]$

$$\eta(s) \geq \min(0, \eta_0) + 2\min_{k=0,\ldots i}\left(\min(0, -v(kT))\right)$$

and since $t \in [iT, (i+1)T)$, the above inequality shows that the left inequality (A.62) holds in this case.

Finally, we notice that the following inequalities hold for all $s \in [iT, (i+1)T]$

$$\eta(iT) \geq \eta(s) \geq \left(1 - (s - iT)T^{-1}\right)\min(0, \eta(iT)) + (s - iT)T^{-1}\min(0, -v(iT))$$

which combined with the facts that $0 \leq 1 - (s - iT)T^{-1} \leq 1$, $0 \leq (s - iT)T^{-1} \leq 1$ implies for all $s \in [iT, (i+1)T]$

$$-|\eta(iT)| - |v(iT)| \leq \eta(s) \leq |\eta(iT)| + |v(iT)|$$

Since $t \in [iT, (i+1)T)$, the above inequality shows that inequality (A.63) holds in this case.

The proof of Claim 3 is complete. ◁

Using (A.34), (A.35), (A.57) and the fact that $|\ln(1+x)| \leq \frac{|x|}{1+\min(x,0)}$ for all $x > -1$, we get the estimate for all $t \geq 0$:

$$|v(t)| = \left|\ln\left(1 + \int_0^A g(a)\psi(t-a)da\right)\right| \leq \frac{\left|\int_0^A g(a)\psi(t-a)da\right|}{1+\min\left(0, \int_0^A g(a)\psi(t-a)da\right)} \leq \frac{\max_{-A \leq s \leq 0}(|\psi(t+s)|)}{1+\min\left(0, \min_{-A \leq s \leq 0}(\psi(t+s))\right)}$$

$$\leq \frac{M \exp(-\sigma t) \max_{-A \leq s \leq 0}(|\psi_0(s)|)}{1+\min\left(0, \min_{-A \leq s \leq 0}(\psi_0(s))\right)} \quad (A.64)$$

Let $j \geq 0$ be an integer with $v(kT) \leq \delta/2$ for all $k \geq j$. We next show that the following inequality holds for all $i \geq j$:

$$\exp(|\eta((i+1)T)|) - 1 \leq \exp(-\delta)(\exp(|\eta(iT)|) - 1) + \exp(2|v(iT)|) - 1 \quad (A.65)$$

Indeed, when $|\eta(iT)| \leq 2\delta$, we get from (A.60) that $|\eta((i+1)T)| \leq 2|v(iT)|$, which directly implies (A.65). On the other hand, when $|\eta(iT)| > 2\delta$, we get from (A.60) that $|\eta((i+1)T)| \leq |\eta(iT)| - 2\delta + 2|v(iT)|$. The previous inequality, in conjunction with the fact that $v(iT) \leq \delta/2$ for all $i \geq j$ gives:

$$\exp(|\eta((i+1)T)|) - 1 \leq \exp(|\eta(iT)| - 2\delta + 2|v(iT)|) - 1$$
$$= \exp(2|v(iT)|) - 1 + \exp(|\eta(iT)| - 2\delta + 2|v(iT)|) - \exp(2|v(iT)|)$$
$$= \exp(2|v(iT)|) - 1 + \exp(2|v(iT)|)(\exp(|\eta(iT)| - 2\delta) - 1)$$
$$= \exp(2|v(iT)|) - 1 + \exp(2|v(iT)| - 2\delta)(\exp(|\eta(iT)|) - 1 + 1 - \exp(2\delta))$$
$$\leq \exp(2|v(iT)|) - 1 + \exp(2|v(iT)| - 2\delta)(\exp(|\eta(iT)|) - 1)$$
$$\leq \exp(2|v(iT)|) - 1 + \exp(-\delta)(\exp(|\eta(iT)|) - 1)$$

Consequently, (A.65) holds for all $i \geq j$. Using (A.65) and induction, we are in a position to prove the following inequality for all $i > j$:

$$\exp(|\eta(iT)|) - 1 \leq \exp(-\delta(i-j))(\exp(|\eta(jT)|) - 1) + \sum_{l=j}^{i-1} \exp(-\delta(i-1-l))(\exp(2|v(lT)|) - 1) \quad (A.66)$$

More specifically, inequality (A.66) follows from the definition of the sequence $\zeta_i := \exp(|\eta(iT)|) - 1$ and the fact that inequality (A.65) gives $\zeta_{i+1} \leq \exp(-\delta)\zeta_i + \exp(2|v(iT)|) - 1$ for all $i \geq j$. Using induction, we can prove the formula $\zeta_i \leq \exp(-\delta(i-j))\zeta_j + \sum_{l=j}^{i-1} \exp(-\delta(i-1-l))(\exp(2|v(lT)|) - 1)$ for all $i > j$, which directly implies (A.66) for all $i > j$.

Using the fact that $x \leq \exp(x) - 1 \leq x\exp(x)$ for all $x \geq 0$, we obtain from (A.66) and (A.64) the following inequality for all $i > j$:

$$|\eta(iT)| \leq \exp(-\delta(i-j))\exp(|\eta(jT)|)|\eta(jT)|$$
$$+ 2\exp(\delta)\sum_{l=j}^{i-1}\exp(2|v(lT)|)\exp(-\delta(i-l))|v(lT)|$$
$$\leq \exp(-\tilde{\delta}(i-j))\exp(|\eta(jT)|)|\eta(jT)| \quad (A.67)$$
$$+ \frac{2M\exp(\delta)\max_{-A\leq s\leq 0}(|\psi_0(s)|)}{1+\min\left(0, \min_{-A\leq s\leq 0}(\psi_0(s))\right)}\exp\left(\frac{2M\max_{-A\leq s\leq 0}(|\psi_0(s)|)}{1+\min\left(0, \min_{-A\leq s\leq 0}(\psi_0(s))\right)}\right)\sum_{l=j}^{i-1}\exp(-\tilde{\delta}(i-l))\exp(-\sigma lT)$$

where $\tilde{\delta} = \min(\delta, \sigma T)$. Since $\tilde{\delta} \leq \sigma T$, it follows that $\exp(-\tilde{\delta}(i-l))\exp(-\sigma lT) \leq \exp(-\tilde{\delta}i)$ for all $l = j,\ldots,i-1$ and thus we obtain from (A.67) the following inequality for all $i > j$:

$$|\eta(iT)| \leq \exp(-\tilde{\delta}(i-j))\exp(|\eta(jT)|)|\eta(jT)|$$
$$+ (i-j)\frac{2M\exp(\delta)\max_{-A\leq s\leq 0}(|\psi_0(s)|)}{1+\min\left(0, \min_{-A\leq s\leq 0}(\psi_0(s))\right)}\exp\left(\frac{2M\max_{-A\leq s\leq 0}(|\psi_0(s)|)}{1+\min\left(0, \min_{-A\leq s\leq 0}(\psi_0(s))\right)}\right)\exp(-\tilde{\delta}i) \quad (A.68)$$

Using (A.62) which implies $|\eta(t)| \leq |\eta_0| + 2\max_{k=0,\ldots[t/T]}(|v(kT)|)$ for all $t \geq 0$ in conjunction with (A.64), we obtain:

$$|\eta(t)| \leq |\eta_0| + \frac{2M\max_{-A\leq s\leq 0}(|\psi_0(s)|)}{1+\min\left(0, \min_{-A\leq s\leq 0}(\psi_0(s))\right)}, \text{ for all } t \geq 0 \quad (A.69)$$

Notice that (A.68) holds for $i = j$ as well and consequently, (A.68) holds for all $i \geq j$. Since $j \geq 0$ is an integer with $v(kT) \leq \delta/2$ for all $k \geq j$, it follows from (A.64) that $j \geq 0$ may be selected as the smallest integer that satisfies $j \geq \frac{1}{\sigma T}\ln\left(\frac{2M\delta^{-1}\max_{-A\leq s\leq 0}(|\psi_0(s)|)}{1+\min\left(0, \min_{-A\leq s\leq 0}(\psi_0(s))\right)}\right)$. The fact that $\tilde{\delta} = \min(\delta, \sigma T)$ in conjunction with (A.68), (A.69) and the fact that $(i-j)\exp(-\tilde{\delta}i/2) \leq i\exp(-\tilde{\delta}i/2) \leq \frac{2}{\tilde{\delta}}\exp(-1)$ for all integers $i \geq j \geq 0$, we get the following inequality for all $i \geq j$:

$$|\eta(iT)| \leq \exp(-\tilde{\delta}i/2)J\left(|\eta_0| + \frac{2M\max_{-A\leq s\leq 0}(|\psi_0(s)|)}{1+\min\left(0, \min_{-A\leq s\leq 0}(\psi_0(s))\right)}\right) \quad (A.70)$$

where $J(s) := s\left(\exp(\tilde{\delta})\delta^{-1}\max(\delta,s) + 2\tilde{\delta}^{-1}\exp(\delta-1)\right)\exp(s)$ for all $s \geq 0$. Using (A.69) and the fact that $j \geq 0$ is the smallest integer that satisfies $j \geq \frac{1}{\sigma T}\ln\left(\frac{2M\delta^{-1}\max_{-A\leq s\leq 0}(|\psi_0(s)|)}{1+\min\left(0, \min_{-A\leq s\leq 0}(\psi_0(s))\right)}\right)$, we can guarantee that (A.70) holds for all $i \geq 0$. Using (A.63), (A.57), (A.34), (A.70) and the fact that $-[t/T] \leq 1 - t/T$, we obtain (5.35) with $\kappa(s) := \exp(\tilde{\delta}/2)\left(J((1+2M)s) + (1+2M)s\right)$ and $L := \tilde{\delta}T^{-1}/2$.

The proof is complete. ◁